\documentclass[11pt]{amsart}

\usepackage[utf8]{inputenc}
\usepackage[T1]{fontenc}
\usepackage[english]{babel}
\usepackage{amsthm, times, xypic}
\usepackage{amscd}
\usepackage{amsfonts}         
\usepackage{amsmath}
\usepackage{amssymb}
\usepackage{hyperref}
\usepackage{enumerate}
\usepackage{tikz}
\usepackage{tikz-cd}
\usepackage{mathrsfs}  

\newcommand{\fref}[1]{\hyperref[{#1}]{\ref*{#1}}}

\newcommand{\B}{\mathbb{B}}
\newcommand{\Laz}{\mathbb{L}}
\newcommand{\A}{\mathbb{A}}
\newcommand{\Proj}{\mathbb{P}}
\newcommand{\R}{\mathbb{R}}
\newcommand{\C}{\mathbb{C}}

\newcommand{\OO}{\mathcal{O}}

\newcommand{\Z}{\mathbb{Z}}

\newcommand{\Li}{\mathscr{L}}

\newcommand{\op}{\mathrm{op}}
\newcommand{\ch}{\mathrm{CH}}

\newcommand{\coim}{\mathrm{\coim}}

\newcommand{\Hom}{\mathrm{Hom}}
\newcommand{\Spec}{\mathrm{Spec}}
\newcommand{\Gr}{\mathrm{Gr}}
\newcommand{\bl}{\mathrm{Bl}}
\newcommand{\wtil}{\widetilde}
\newcommand{\QCoh}{\mathrm{QCoh}}
\newcommand{\Perf}{\mathrm{Perf}}
\newcommand{\Sym}{\mathrm{Sym}}
\newcommand{\LSym}{\mathrm{LSym}}
\newcommand{\sSpec}{\cS pec}
\newcommand{\abs}[1]{\lvert #1 \rvert}
\newcommand{\iHom}{\underline{\mathrm{Hom}}}

\newtheorem{theo}{Tplottin ubuntuheorem}[section]
\theoremstyle{plain}
\newtheorem{thm}[theo]{Theorem}
\newtheorem{lem}[theo]{Lemma}
\newtheorem{prop}[theo]{Proposition}
\newtheorem{cor}[theo]{Corollary}
\newtheorem*{thm*}{Theorem}
\newtheorem*{lem*}{Lemma}
\newtheorem*{prop*}{Proposition}
\newtheorem*{cor*}{Corollary}

\theoremstyle{definition}
\newtheorem{defn}[theo]{Definition}

\newtheorem{cons}[theo]{Construction}
\newtheorem{ex}[theo]{Example}

\theoremstyle{remark}
\newtheorem{rem}[theo]{Remark}
\newtheorem{war}[theo]{Warning}

\def\alp{\alpha}
\def\be{\beta}

\def\Cal{\mathscr}
\def\ga{\gamma}

\def \bB{\mathbb B}

\def \bOB{\mathbb {OB}}
\def \bM{\mathbb M}
\def \bOM{\mathbb {OM}}

\def\bZ{\mathbb Z}

\def \bI{\mathbb I}
\def\op{\operatorname}
\def \cF{\mathcal F}
\def \cE{\mathcal E}
\def \cC{\mathcal C}

\def \cS{\mathcal S}
\def \cN{\mathcal N}

\numberwithin{equation}{section}
\title{Bivariant algebraic cobordism with 
bundles}
\author{Toni Annala}

\address{Department of Mathematics,
University of British Columbia, 
Room 121, 1984 Mathematics Road
Vancouver, BC
Canada V6T 1Z2}
\email
{tannala@math.ubc.ca}

\author{Shoji Yokura}

\address{Department of Mathematics and Computer Science,
Graduate School of Science and Engineering\\ Kagoshima University,
21-35 Korimoto 1-Chome,
Kagoshima 890-0065, Japan}
\email
{yokura@sci.kagoshima-u.ac.jp}

\date{}
\begin{document}

\maketitle

\begin{abstract} The purpose of this paper is to study an extended version of bivariant derived algebraic cobordism where the cycles carry a vector bundle on the source as additional data. We show that, over a field of characteristic 0, this extends the analogous homological theory of Lee and Pandharipande constructed earlier. We then proceed to study in detail the restricted theory where only rank $1$ vector bundles are allowed, and prove a weak version of projective bundle formula for bivariant cobordism. Since the proof of this theorem works very generally, we introduce precobordism theories over arbitrary Noetherian rings of finite Krull dimension as a reasonable class of theories where the proof can be carried out, and prove some of their basic properties. These results can be considered as the first steps towards a Levine-Morel style algebraic cobordism over a base ring that is not a field of characteristic 0.

\end{abstract}

\tableofcontents

\section{Introduction}

Algebraic cobordism is the most general oriented cohomology theory in algebraic geometry, and therefore the cobordism ring of a scheme $X$ contains a lot of geometric information about $X$. For example, two other fundamental invariants --- the Chow ring of $X$ and the $K$-theory of vector bundles on $X$ --- can be easily recovered from the algebraic cobordism ring of $X$. Moreover, the ring admits a purely geometric description: its elements are precisely the cobordism classes of nice enough varieties over $X$. Such a description is very interesting even in the case when $X = \Spec(k)$ is a spectrum of a field, and it has been used previously to prove conjectures in enumerative geometry (degree 0 Donaldson-Thomas conjectures, see \cite{LP}). More recently, Haution studies in \cite{Hau} the fixed loci of involutions of varieties using a ring closely related to the cobordism ring of $\Spec(k)$. For another recent application of algebraic cobordism, although in a slightly different spirit, see \cite{SeSe}, where Sechin and Semenov study algebraic groups using Morava $K$-theories, which are invariants constructed from algebraic cobordism by adding relations.

However, much like the theory of Chow rings, algebraic cobordism is still mostly hypothetical. Following the pioneering work of Voevodsky, Levine--Morel in their foundational book \cite{LM} were able to construct satisfactory algebraic cobordism rings for varieties smooth over a field of characteristic 0. The reason for the restriction to characteristic 0 is the liberal use of resolution of singularities and weak factorization in \cite{LM}, which seem to be essential for most of the proofs. On the other hand, the reason for the smoothness assumption is the same as in intersection theory: without it, it becomes very hard to define an intersection product. In \cite{An}, building on \cite{LS}, the first named author was able to overcome the smoothness assumptions (while staying in characteristic 0) using derived algebraic geometry. The purpose of this paper is to continue this work. More precisely, the initial goal was to perform a similar extension of the algebraic bordism of varieties with vector bundles of Lee-Pandharipande \cite{LeeP}. However, quite surprisingly, many of the proofs did not use the caracteristic 0 assumption in any essential way, allowing us to get nontrivial results over a more general base ring.

Our constructions are in fact much more general than just a construction of a cohomology theory. The cohomology rings are just a small part of a larger bivariant theory, which also contains a corresponding homology theory (e.g. Chow groups, $K$-theory of coherent sheaves, algebraic bordism), and forces these theories to behave extremely well. Recall that a \emph{bivariant theory} $\B^*$ assigns an Abelian group $\B^*(X \to Y)$ for every morphism between quasi-projective (derived) schemes, and natural operations associated to compositions (\emph{bivariant product}), factorizations through projective morphisms (\emph{bivariant pushforwards}) and (homotopy) Cartesian squares (\emph{bivariant pullbacks}). One can recover the corresponding homological and cohomological theories as
$$\B_*(X) := \B^{-*}(X \to pt)$$
and
$$\B^*(X) := \B^{*}(X \xrightarrow{\mathrm{Id}} X)$$
respectively. The bivariant formalism is recalled in detail in Section \fref{FM-BT}.

\subsection*{Summary of main results}

In Section \fref{Char0ConstructionSection}, working over a field of characteristic 0, we construct the bivariant cobordism $\Omega^{*,*}$ with vector bundles, generalizing both the algebraic cobordism of bundles $\omega_{*,*}$ 
constructed in \cite{LeeP} 
and the bivariant derived algebraic cobordism $\Omega^*$ constructed in \cite{An}. The bivariant group $\Omega^{d,r}(X \to Y)$ is generated by cycles of form
$$[V \to X, E],$$
where the morphism $V \to X$ is projective, the composition $V \to Y$ is quasi-smooth of relative dimension $-d$ and $E$ is a rank $r$ vector bundle on $V$. The relations are similar to those used to construct $\Omega^*$ in \cite{An}. We then have natural isomorphisms
$$\Omega^{*,*}(X \to pt) \cong \omega_{-*, *}(tX) \quad (\text{Theorem \fref{LeePIso}})$$
and 
$$\Omega^{*,0}(X \to Y) \cong \Omega^*(X \to Y) \quad (\fref{CobOfBundlesExtendsCob})$$
for all quasi-projective derived schemes $X$ and $Y$, where $tX$ denotes the classical truncation.

In Section \fref{CobordismOfLines} we study the structure of the restricted theory $\Omega^{*,1}$, henceforth called \emph{bivariant} \emph{cobordism with line bundles}, and show that \emph{over an arbitrary Noetherian base ring of finite Krull dimension} it has, in a sense, a nice basis over $\Omega^*$ (Theorems \fref{AlgebraicCobordismOfBGm} and \fref{KunnethFormulaForBGm}). This is a fundamental computation, the first of its kind in the study of derived cobordism theories, and the other results in Section \fref{CobordismOfLines} --- as well as several results of \cite{An4} --- are based on it. In order to work in the aforementioned generality, we initiate the study of \emph{bivariant precobordism theories} (and associated theories of line bundles) in Section \fref{PrecobordismSect}, which give a very general class of bivariant theories for which the results of Section \fref{CobordismOfLines} hold. The computation of the structure of $\Omega^{*,1}$ is done in Section \fref{StructOfPrecobOfLinesSect}. We then show in Section \fref{WPBFSect} that the bivariant groups $\Omega^*(\Proj^n \times X \to Y)$ embed in a natural way into $\Omega^{*,1}(X \to Y)$ allowing us to express $\Omega^*(\Proj^n \times X \to Y)$ additively as a direct sum of $n+1$ copies of $\Omega^*(X \to Y)$ (the \emph{weak projective bundle formula} --- Theorem \fref{WPBF}). The results are completely new whenever the main theorem of \cite{LS} does not apply, i.e.:
\begin{itemize}
\item when $Y$ is not smooth;
\item when the base ring $A$ is not a field of characteristic 0.
\end{itemize}
For a simple example, take $\Omega^*(\Proj^n \to pt)$ over a field of positive characteristic.

The computation of precobordism rings of products of projective spaces allow us to conclude that the behavior of the first Chern class in tensor products is controlled by a formal group law (Theorem \fref{PrecobordismFGL}). This is nontrivial, since the definition of a bivariant precobordism theory (Definition \fref{PrecobordismDef}) does not explicitly enforce a formal group law to hold (unlike in the construction of $\Omega^*$ in \cite{An}). Compare this to Proposition 8.3 of \cite{LP}. This also shows that we have a natural morphism $\Laz \to \Omega^*(pt)$, where $\Laz$ is the Lazard ring (Corollary \fref{LazardMapsToPrecobordism}).

As a final remark, we note that the results of Section \fref{CobordismOfLines} are stated in terms of general precobordism theories $\B^*$. This means that they will hold for the universal precobordism theory $\underline{\Omega}^*$, and, more importantly, for any bivariant theory obtained from $\underline\Omega^*$ by adding relations (i.e., taking a quotient by a \emph{bivariant ideal}, see Definition \fref{BI}). This is a crucial point: not only does one expect there to be many interesting bivariant theories receiving a surjective Grothendieck transformation from $\underline{\Omega}^*$, but also it is likely that the ``correct'' bivariant algebraic cobordism $\Omega^*$ is a quotient of $\underline\Omega^*$ by further relations.

\subsection*{Related work}
Algebraic cobordism $\op{MGL}^{*,*}(X)$ (now called higher algebraic cobordism) was first introduced by V. Voevodsky in the context of motivic homotopy theory and was used in his proof of the Milnor conjecture \cite{Vo1, Vo2, Vo3}. Later, in an attempt to better understand this higher algebraic cobordism, M. Levine and F. Morel constructed another algebraic cobordism 
$\Omega_*(X)$ in terms of cobordism cycles (of the form $[Y \xrightarrow f X; L_1, L_2, \cdots, L_r]$ where $L_i$'s being line bundles over the source variety $Y$) and some relations, as the universal oriented cohomology theory. To be a bit more precise, they first defined \emph{an oriented Borel--Moore functor with products} satisfying 12 conditions (D1) - (D4) and (A1) - (A8). Such a functor $\mathcal Z_*$ was constructed by using cobordism cycles.  Secondly they defined \emph{an oriented Borel--Moore functor with products of geometric type} by further imposing, on the oriented Borel--Moore functor with products, three axioms (Dim) (dimension axiom), (Sect) (section axiom) and (FGL) (formal group law axiom), which correspond to \emph{``of geometric type''}. They constructed such a functor $\mathcal Z_*/ \mathcal R$ by imposing relations $\mathcal R$ corresponding to these three axioms on the functor $\mathcal Z_*$. The functor  $\mathcal Z_*/\mathcal R$ is nothing but Levine--Morel's algebraic cobordism $\Omega_*$. In \cite{Le} M. Levine showed that there is an isomorphism $\Omega^*(X) \cong \op{MGL}^{2*,*}(X)$ for smooth $X$.

In \cite{LP} M. Levine and R. Pandharipande constructed what they call \emph{double point cobordism} $\omega_*(X)$ more geometrically and more simply, using \emph{double point relation}, which is similar to rational equivalence relation to define Chow group and they showed the isomorphism $\Omega_*(X) \cong \omega_*(X)$. They consider the set $\mathcal M(X)$ of isomorphism classes $[Y \xrightarrow f X]$ of projective morphisms $f:Y \to X$ with smooth $Y$ and $\mathcal M_*(X)$ is a monoid under disjoint union of domains and is graded by the dimension of $Y$ over the ground field $k$. Furthermore let $\mathcal M_*(X)^{+}$ denote the group completion of the monoid $\mathcal M_*(X)$ and let $\mathcal R_*(X) \subset \mathcal M_*(X)$ be the subgroup generated by all the double point relations over $X$. Then Levine--Pandharipande's double point cobordism $\omega_*(X)$ is defined to be the quotient $\frac{\mathcal M_*(X)}{\mathcal R_*(X)}$. A crucial and important difference from the construction of Levine--Morel's algebraic cobordism is that they consider $[Y \xrightarrow f X]$ without line bundles $L_i$'s, which are key ingredients of Levine--Morel'c cobordism cycles $[Y \xrightarrow f X; L_1, L_2, \cdots, L_r]$.

In \cite{LS} P. Lowrey and T. Sch\"urg extended Levine--Morel's algebraic cobordism to derived algebraic schemes and called it \emph{derived algebraic cobordism}, denoted by $d\Omega_*(X)$, and they showed the isomorphism $d\Omega_*(X) \cong \Omega_*(tX)$ for any derived quasi-projective scheme $X$. Above, $tX$ is the \emph{classical truncation} of $X$.

 In \cite{FM} W. Fulton and R. MacPherson have introduced \emph{Bivariant Theory} with an aim to deal with Riemann--Roch type theorems or formulas and to unify them. In this paper a bivariant theory is this Fulton--MacPherson's bivariant theory. With the aim to construct a bivariant-theoretic version $\Omega^*(X \xrightarrow f Y)$ of Levine--Morel's algebraic cobordism so that for a map to a point $\pi_X:X \to pt$ the covariant part $\Omega^*(X \xrightarrow {\pi_X} pt)$ is isomorphic to Levine--Morel's algebraic cobordism $\Omega_*(X)$, the second named author introduced \emph{an oriented bivariant theory} and \emph{a universal bivariant theory} in \cite{Yo1} (cf. \cite{Yo2}). In \cite{An}, using the construction of the above universal bivariant theory,  the first named author has constructed a bivariant-theoretic version of Lowrey--Sch\"urg's derived algebraic cobordism $d\Omega^*(X \xrightarrow f Y)$, called \emph{bivariant derived algebraic cobordism}, in such a way that for $Y=pt$ is a point, $d\Omega^*(X \xrightarrow {\pi_X} pt) \cong d\Omega_*(X)$, thus it follows from Lowrey--Sch\"urg's theorem above that 
$d\Omega^*(X \xrightarrow {\pi_X} pt) \cong \Omega_*(X)$.

In \cite{LP} Y.-P. Lee and R. Pandharipande extended the construction of Levine--Pandharipande's double point cobordism $\omega_*(X)$ to a double graded double point cobordism $\omega_{*,*}(X)$, considering 
similar double point relations on the group completion of the monoid generated by the isomorphism classes of $[Y \xrightarrow f X; E]$ instead of $[Y \xrightarrow f X]$ (see above). Here $E$ is a complex vector bundle over the source $Y$ and the second grade $*$ of $\omega_{*,*}(X)$ refers to the rank of this vector bundle $E$. This cobordism $\omega_{*,*}(X)$ is nothing but what they call \emph{algebraic cobordism of bundles on varieties}.


\subsection*{Structure of the paper}

The structure of the paper is as follows: In \S 2 we give an introduction to derived algebraic geometry with more emphasis on results or facts which we need in later sections. As a rule, only new results are given proofs. In \S 3 we give a quick introduction to Fulton--MacPherson's bivariant theory, since we sometimes refer to the seven axioms required on their bivariant theory. We also recall the notion of \emph{bivariant ideal}, which is introduced in \cite{An}, and the universal bivariant theory \cite{Yo1}. In \S4 we recall oriented bivariant theory and a universal oriented bivariant theory \cite{Yo1}, since we refer to them in proofs in later sections. The main results are obtained in \S5 and \S6.


\subsection*{Future work}

The study of bivariant precobordism theories is continued in \cite{An4}, where Chern classes are constructed for bivariant precobordism theories, allowing the generalization of most of the results of \cite{An} over general Noetherian base ring $A$ (Conner-Floyd, Grothendieck-Riemann-Roch). This, together with the computation of precobordism groups of line bundles in Section \fref{CobordismOfLines}, is then used to prove the validity of the projective bundle formula for general projective bundles $\Proj(E)$. Determining the structure of the bivariant theory $\Omega^{*,*}$ over $\Omega^*$ reduces, with some effort, to the case of line bundles, and will appear elsewhere. Using similar ideas as in here and in \cite{An4}, the structure of $\Omega^{*,r}$ can be used to compute the precobordism groups of Grassmannians of $r$-planes and bundles thereof.

\subsection*{Acknowledgements}

The authors would like to thank David Rydh and Adeel Khan for answering questions about derived blow ups, as well as Pavel Sechin and several other people for spotting a mistake. The first author would also like to thank his advisor Kalle Karu for multiple discussions. The first author is supported by Vilho, Yrj\"o and Kalle V\"ais\"al\"a Foundation of the Finnish Academy of Science and Letters. The second author is partially supported by JSPS KAKENHI Grant Numbers JP16H03936 and JP19K03468.

\section{$\infty$-categories and derived algebraic geometry}

Before we can give an introduction to derived algebraic geometry, we must ask ourselves a question: \emph{what is homotopy theory}? Or more specifically: what is the correct categorical structure underlying the theory of \emph{homotopy types}, i.e.,  spaces (CW complexes) considered up to homotopy equivalence? The problem, which is essentially about localization, appears to be a very innocent one: can you not just consider the category obtained from spaces by making all homotopy equivalences isomorphisms? One quickly notices that the category obtained this way is the \emph{homotopy category of spaces}: the category whose objects are spaces, and whose morphisms are homotopy classes of continuous maps.


Unfortunately, this does turn out to be insufficient as a categorical framework. Moreover, the shortcomings should be familiar to many algebraic geometers, although in a disguised form, for their (partial) resolution led to the now classical theory of stacks. A stack is supposed to be nothing more exotic than a sheaf of groupoids. For example, consider line bundles on a scheme $X$. We can associate for every open $U \subset X$ the groupoid of line bundles on $U$, and such an assignment clearly determines some sort of a presheaf of groupoids on $X$. The formalization of the sheaf condition is slightly more subtle: for example, to get a line bundle $\Li$ on $X$, it is \emph{not} enough to consider an open cover $(U_i)$ on $X$ and line bundles $\Li_i$ on $U_i$ that become pairwise isomorphic when restricted to the intersections $U_i \cap U_j$. We must also \emph{choose} isomorphisms $\Phi_{ji} : \Li_i \vert_{U_i \cap U_j} \to \Li_j \vert_{U_i \cap U_j}$ for all $i,j$ satisfying the  \emph{cocycle condition}
\begin{equation*}
\Phi_{ki} =  \Phi_{kj} \circ \Phi_{ji}
\end{equation*}
as maps $\Li_i \vert_{U_i \cap U_j \cap U_k} \to \Li_k \vert_{U_i \cap U_j \cap U_k}$. This is the general form of the sheaf condition for stacks on a topological space.

In order to see that the above problem fits in the framework of homotopy theory, we need to recall two facts:
\begin{itemize}
\item every groupoid is equivalent to the fundamental groupoid of a space that has no higher homotopy groups (a \emph{1-truncated space}); 
\item the homotopy type of a 1-truncated space is completely determined by its fundamental groupoid. 
\end{itemize}
Hence we may consider the above stack of line bundles as a presheaf $\cF$ taking values in 1-truncated spaces. Moreover, the cocycle condition for descent is replaced by its homotopical analogue saying that the space $\cF(X)$ of line bundles on $X$ should be obtained as the \emph{homotopy limit} over a certain diagram containing spaces $\cF(U_J)$, where $J \subset I$ and $U_J = \bigcap_{j \in J} U_j$. More precisely, the space of line bundles on $X$ should be naturally homotopy equivalent to the space (equipped with the compact-open topology) consisting of the data:
\begin{enumerate}
\item points $x_i \in \cF(U_i)$;
\item paths $\phi_{ji}$ from the image of $x_i$ to the image of $x_j$ in $\cF(U_{ij})$;
\item homotopies from the image path $\phi_{kj}\phi_{ji}$ to $\phi_{ki}$ in $\cF(U_{ijk})$ (filling the triangle);
\item fillings of the hollow tetrahedron formed by the images of the triangles obtained in (3) in $\cF(U_{ijkl})$
\item $\cdots$
\end{enumerate}
As $\cF$ takes values in 1-truncated spaces, we note that if a hollow triangle as in (3) can be filled, then all the different fillings must be homotopic (otherwise we would have a nontrivial element in the second homotopy group), and similarly that all the higher dimensional hollow tetrahedra can always be filled in a unique way (up to homotopy). Hence, passing to the fundamental groupoids, one obtains an equivalence between the classical theory of stacks and the more homotopical theory of sheaves of 1-truncated spaces. Note that without the assumption that $\cF$ takes values in 1-truncated spaces, none of the conditions (1), (2), (3)... become obsolete in general.

We can also see why passing to the homotopy category is a bad idea. Indeed, up to homotopy, a space only has one point for each path component, and hence trying to obtain the homotopy equivalence class of line bundles on $X$ as a limit over the homotopy equivalence classes of line bundles on the open sets produces nonsense: it would be the same as saying that a collection of line bundles on $U_i$ glue uniquely as soon as they are isomorphic on the pairwise intersections, which we already noted was wrong. 

\subsection{$\infty$-categories}
  
If the homotopy category is not the right framework for doing homotopy theory, then what is? As we already noticed, the problem is that we 
need to remember not only that certain maps are homotopic, but we also have to remember these homotopies, higher homotopies between homotopies, and so on. In other words, we need to keep track of the \emph{homotopy coherence data}. The modern way of doing this is using the theory of $\infty$-categories.

An \emph{$\infty$-category}, or more precisely an \emph{$(\infty,1)$-category}, is essentially a category enriched over spaces up to coherent homotopy equivalences. They were born out of the preceding notion of \emph{Quillen model categories}, and were (still are) a useful gadget for dealing with questions such as localization of categories, and computing homotopy (co)limits. We are not going to try to give a very comprehensive survey to the theory of $\infty$-categories --- the interested reader can consult the first chapter of \cite{Lur1} for a concise introduction of 50 pages. The following is essentially the Definition 1.1.2.4 of \emph{loc. cit.}

\begin{defn}
An \emph{$\infty$-category} is a simplicial set satisfying the weak Kan condition. An \emph{$\infty$-groupoid} is a simplicial set satisfying the full Kan condition. The latter notion should be understood as a combinatorial model for the \emph{homotopy type} of a nice enough topological space (a space having the homotopy type of a CW-complex).
\end{defn}

If $\cC$ is an $\infty$-category, then one regards its $0$-simplices as its \emph{objects}, its $1$-simplices as \emph{morphisms} between its endpoints, and the higher dimensional simplices give rise to compositions and (higher) homotopies between morphisms. Given two objects $x, y$ of $\cC$, one may associate to them in a natural way the \emph{mapping space} $\Hom_\cC(x,y)$ which is a Kan-complex. 

The most fundamental example of an $\infty$-category is the $\infty$-category $\cS$ of \emph{spaces}. The objects of $\cS$ are Kan complexes, and the mapping spaces $\Hom_\cS(x,y)$ are naturally equivalent to the internal mapping object $\iHom(\abs x, \abs y )$ of (compactly generated weakly Hausdorff) spaces, where $\abs{-}$ denotes the \emph{geometric realization of simplicial sets}. It is known that, in a precise sense, $\cS$ captures all the information of topological spaces considered up to weak homotopy equivalence.

Of course, in order to do meaningful category theory, one also needs to know the definition of a functor.
\begin{defn}
Given two $\infty$-categories $\mathcal C$, $\mathcal D$ a \emph{functor} $F: \mathcal{C} \to \mathcal{D}$ is a morphism of simplicial sets. 
\end{defn}
Hence, in order to give a functor between $\infty$-categories, it is not enough to know where objects and morphisms go: one must also know everything about the higher dimensional simplices. This makes it sometimes combinatorially very challenging to write down functors by hand and often the easiest way to proceed is to invoke some sort of a universal property.

A functor $F: \mathcal{C} \to \mathcal{D}$ naturally induces morphisms of mapping spaces
$$\Hom_{\cC}(x,y) \to \Hom_{\mathcal{D}}(F(x),F(y))$$
which allows us to make the following definition.

\begin{defn}
A functor $F: \mathcal{C} \to \mathcal{D}$ is \emph{fully faithful} if for all pairs $(x,y)$ of objects of $\cC$, the induced morphism
$$\Hom_\cC(x,y) \to \Hom_{\mathcal{D}}(F(x), F(y))$$
is an equivalence.
\end{defn}
As an example, given a subset $C'$ of objects of $\cC$ one may form the \emph{full subcategory} $\cC'$ of $\cC$ which is the simplicial subset of $\cC$ containign only simplices whose $0$-simplices are in $C'$. The natural inclusion $\cC' \to \cC$ of simplicial sets is fully faithful as a functor.

There is an $\infty$-categorical version of initial and final objects.
\begin{defn}
Let $\cC$ be an $\infty$-category and $x$ an object of $\cC$. Then $x$ is
\begin{enumerate}
\item \emph{initial} if for every object $y$ of $\cC$, the space of morphisms $\Hom_\cC(x,y)$ is contractible;
\item \emph{final} if for every object $y$ of $\cC$, the space of morphisms $\Hom_\cC(y,x)$ is contractible;
\item \emph{zero object} if $x$ is both initial and final.
\end{enumerate} 
\end{defn} 
We can also talk about \emph{overcategories} and \emph{undercategories} (see \cite{Lur1} Section 1.2.9). We will only define overcategories, as the notion of an undercategory is dual to it.
\begin{defn}
Suppose $\cC$ is an $\infty$-category, and $x$ is an object of $\cC$. We can now construct an $\infty$-category $\cC_{-/x}$ of \emph{objects of $\cC$ over $x$} whose 
\begin{enumerate}
\item objects are morphisms $\pi_y: y \to x$ in $\cC$;
\item morphisms between $\pi_y: y \to x$ and $\pi_{y'}: y' \to x$ are $2$-simplices in $\cC$ of form
\begin{center}
\begin{tikzcd}
y \arrow[->]{r} \arrow[->]{rd}[swap]{\pi_y} & y' \arrow[->]{d}{\pi_{y'}} \\ 
 & x; 
\end{tikzcd}
\end{center}
\item generally $n$-simplices are $(n+1)$-simplices in $\cC$, whose endpoint is $x$.
\end{enumerate}
There is a natural \emph{forgetful functor} $\cC_{-/x} \to \cC$ forgetting the structure morphisms.
\end{defn}
The above definition is merely a special case: more generally, given a diagram in an $\infty$-category, one may form the $\infty$-category of objects over/under the diagram. A final object in such an overcategory, if it exists is by definition the \emph{$\infty$-categorical limit} of the corresponding diagram. Dually, one defines the \emph{$\infty$-categorical colimit} as the initial object in an undercategory. Rather than explaining the general definition of (co)limits in full detail (which is not terribly complicated, see Section 1.2.13 of \cite{Lur1}), we will give the universal properties in the special cases we are going to use in the article.

\begin{rem}
The $\infty$-categorical overcategories can be more complicated than their classical counterparts. Consider for example the definition of a morphism over $x$. Indeed, the above definition may be translated into the following: a morphism between $\pi_y: y \to x$ and $\pi_{y'}: y' \to x$ is a morphism $f: y \to y'$ together with a path (a \emph{homotopy}) $\alpha$ from $\pi_{y'} \circ f$ to $\pi_y$ in the space of morphisms $\Hom_{\cC}(y \to x)$. Even if there was, up to homotopy, only one morphism $f: y \to y'$, the mapping space $\Hom_{\cC_{-/x}}(y,y')$ may still have multiple components.

For a simple example, consider the space of endomorphisms of a point $pt \to S^1$ in the $\infty$-category $\cS_{-/S^1}$ of spaces over the circle (for concreteness, let $pt$ be included as the point $(1,0) \in S^1$). Recall that the internal mapping space $\iHom(pt, S^1)$ is homeomorphic to $S^1$. As there is only one morphism $pt \to pt$, we see that the only data in a $S^1$-morphism $pt \to pt$ is a loop in $S^1$ based at $(1,0)$, and therefore $\Hom_{\cS_{-/S^1}}(pt,pt)$ is naturally equivalent to the space of loops $\Omega S^1$ in $S^1$ based at $(1,0)$. As $\Omega S^1 \simeq \Z$, we see that, even up to homotopy, there are infinitely many maps $pt \to pt$ over $S^1$.
\end{rem}

\subsection{Derived algebraic geometry --- basic definitions}

Derived algebraic geometry is obtained from algebraic geometry by replacing commutative rings with simplicial commutative algebras. There is a forgetful functor from the $\infty$-category of simplicial commutative rings (and from virtually any other $\infty$-category of derived algebraic objects, for example chain complexes) to the $\infty$-category of spaces sending a simplicial commutative ring to its \emph{underlying simplicial set} (or its \emph{underlying space}).
  
\begin{defn}
A \emph{derived scheme} $X$ is a topological space $X_\mathrm{top}$ equipped with a hypercomplete sheaf $\OO_X$ of simplicial commutative rings such that:
\begin{enumerate}
\item the \emph{truncation} $tX = (X_\mathrm{top}, \pi_0(\OO_X))$ is a scheme in the usual sense; 
\item the higher homotopy sheaves $\pi_i(\OO_X)$ (which descend to sheaves on the truncation) are quasi-coherent.
\end{enumerate}
\end{defn}

\begin{rem}
In the above definition, the underlying topological space $X_\mathrm{top}$ is \emph{not} considered up to any kind of homotopy equivalence. In derived algebraic geometry, we only consider the ring of functions up to homotopy, not the space on which these functions are defined on.
\end{rem}

\begin{rem}
It is known that whenever the underlying space $X_\mathrm{top}$ is Noetherian and finite dimensional, then the hypercompleteness assumption is automatically satisfied. This is the reason why we restrict our attention to Noetherian rings that in addition have finite Krull dimension.
\end{rem}

\begin{rem}
The condition of being a sheaf can be checked on the level of the underlying spaces (or on the underlying simplicial abelian group/connective chain complex), and it is essentially the same condition as in the introduction of this section. Namely, given an open set $U \subset X$ and an open cover $(V_i)$ of $U$, we want the space of sections on $U$ to be naturally equivalent to the homotopy limit of the spaces of a certain diagram containing the spaces of global sections on $V_i$ and on their intersections.

A classical example of an $\infty$-categorical sheaf from algebraic geometry is given by the \emph{hypercohomology} of a complex of coherent sheaves on a scheme $X$. A complex of sheaves $\cF_\bullet$ can be regarded as a presheaf on the underlying topological space $X_\mathrm{top}$ taking values in a suitable $\infty$-category of chain complexes, but it might fail to be a sheaf. Fortunately, due to the coherence assumption the sheaf condition is known to hold when the complex is restricted to affine schemes, and one can use this fact to understand the sections of the $\infty$-categorical sheafification of $\cF_\bullet$ on arbitrary open sets $U \subset X$: just choose an affine open cover $(V_i)$  of $U$ and use the sheaf property to compute the sections on $U$ as a homotopy limit. The cohomology groups of the sections on $U$ recover the \emph{hypercohomology} groups of $\cF\vert_U$, and they can be computed for example by using a spectral sequence or an injective resolution.

As a special case, one sees how to obtain a derived scheme from a \emph{dg-scheme} in the sense of Ciocan-Fontanine and Kapranov (see \cite{CFK}). A \emph{dg-scheme} is a topological space equipped with a differential graded sheaf of algebras satisfying some extra conditions. The structure sheaf can be regarded as a presheaf of derived $k$-algebras (differential graded $k$-algebras), and one obtains a derived scheme by applying the homotopical sheafification functor. As above, this process does not affect the sections on affine opens.
\end{rem}

Of course, in order to do derived algebraic geometry, it is not enough to know only the objects of the theory, but also the morphisms.

\begin{defn}
A \emph{morphism} $f: X \to Y$ between derived schemes is a continuous map $f: X_\mathrm{top} \to Y_\mathrm{top}$ of topological spaces and a map $f^\sharp:\OO_Y \to f_* \OO_X$ of the structure sheaves so that $(f, \pi_0 f^\sharp)$ defines a map of schemes.
\end{defn}
\begin{rem}
The pushforward $f_*$ of sheaves is defined in exactly the same way that it is defined in classical theory of sheaves.
\end{rem}
\begin{rem}
The ordinary category of classical schemes embeds fully faithfully to the $\infty$-category of derived schemes. Hence derived algebraic geometry truly is an extension of classical algebraic geometry.
\end{rem}

As in classical algebraic geometry, there are many important subclasses of morphisms. We will record the ones relevant to the article below.

\begin{defn}\label{ClosedEmbeddingDef}
A morphism $f: X \to Y$ of derived schemes is a \emph{closed embedding} if its truncation $tf : tX \to tY$ is a closed embedding in the classical sense.
\end{defn}

\begin{defn}\label{ProperDef}
A morphism $f: X \to Y$ of derived schemes is \emph{proper} if the truncation $tf : tX \to tY$ is proper in the classical sense. 
\end{defn}

\begin{defn}\label{QProjDef}
A morphism $f: X \to Y$ of derived schemes is \emph{quasi-projective} if it factors as a composition of a closed embedding $i: X \hookrightarrow U \times Y$ and a projection $p_2: U \times Y \to Y$, where $U$ is an open subscheme of a projective space $\Proj^n$. 

A quasi-projective morphism $f: X \to Y$ that is also proper is called \emph{projective}. This is equivalent to the existence of a factorization as a composition of a closed embedding $i: X \hookrightarrow \Proj^n \times Y$ and a projection $p_2: \Proj^n \times Y \to Y$. 
\end{defn}

\begin{war}
Unlike the first two classes of morphism, the quasi-projectivity of a morphism of derived schemes can \emph{not} be checked from the truncation. Indeed, counter examples for this were constructed in \cite{An3}: there exists derived schemes whose truncations are projective hypersurfaces, but which still fail to have any non constant maps to a projective space.
\end{war}

In the main part of the article, we will restrict out attention to quasi-projective derived schemes, for there are many technical benefits in working with them. Traditionally there are two reasons for this. First of all, a morphism between quasi-projective derived schemes is itself quasi-projective, and can therefore factored as a composition of a closed embedding and a smooth morphism. This kind of global factorization is sometimes extremely useful. Secondly, as we will see later, the theory of quasi-coherent sheaves is much simpler on quasi-projective schemes.

\subsection{Derived fibre products}

Fibre products are an indispensable tool in algebraic geometry, and the same is true for derived algebraic geometry. We begin with a definition.

\begin{defn}\label{DerivedFibreProduct}
Let $\pi_X: X \to S$ and $\pi_Y: Y \to S$ be two derived schemes over a base derived scheme $S$. The \emph{derived fibre product} (sometimes \emph{homotopy fibre product}) $X \times_S^\R Y$ is the $S$-scheme so that given an $S$-scheme $Z$, we have a natural equivalence 
$$\Hom_{-/S}(Z, X \times_S^\R Y) \simeq \Hom_{-/S}(Z, X) \times \Hom_{-/S}(Z, X)$$
of morphism spaces over $S$.

In other words, the space $\Hom(Z, X \times_S^\R Y)$ of morphisms of derived schemes (not over $S$) is naturally equivalent to the space consisting of
\begin{enumerate}
\item morphisms $f_1: Z \to X$ and $f_2: Z \to Y$;
\item a homotopy $\alpha: \pi_X \circ f_1 \simeq \pi_Y \circ f_2$, i.e., a path inside $\Hom(Z, S)$ connecting the two above morphisms.
\end{enumerate}
\end{defn}

As in the classical theory, derived schemes admit all derived fibre products. A commutative square of form
\begin{equation*}
\begin{tikzcd}
X \arrow[->]{r}{\pi'_X} \arrow[->]{d}{\pi'_Y} \times_S^\R Y & Y \arrow[->]{d}{\pi_Y} \\
X \arrow[->]{r}{\pi_X} & S
\end{tikzcd}
\end{equation*}
is called \emph{derived Cartesian} or \emph{homotopy Cartesian} and the map $\pi_X'$ (resp. $\pi_Y'$) is called \emph{derived/homotopy pullback} of $\pi_X$ (resp. $\pi_Y$). Sometimes also the name \emph{derived/homotopy pullback square} is used. Finally, if the two morphisms $\pi_X$ and $\pi_Y$ are closed embeddings, then the derived fibre product is sometimes called \emph{derived/homotopy intersection}.

\begin{rem}
We have chosen to denote the derived fibre product with the symbol $\times_S^\R$, that is, with superscript $\R$. This is done in order to avoid confusing the derived fibre product with the usual one when both make sense, i.e., in the case when all the schemes $X,Y,S$ are classical. The choice of $\R$ stems from the fact that affine locally the fibre product is modeled by the opposite of the derived tensor product of derived rings, which is a non Abelian left derived functor in the sense of Quillen.
\end{rem} 

We note that if $X,Y$ and $S$ are classical schemes, then there exists a natural morphism 
$$X \times_S Y \to X \times_S^\R Y,$$
but this is not an equivalence in general. However, this is true if $X$ and $Y$ are \emph{Tor-independent} over $S$. The following familiar properties remain true for derived fibre products:

\begin{prop}\label{CompositionAndCancellationOfFibreSquares}
Let us have a commutative diagram
\begin{equation*}
\begin{tikzcd}
Z' \arrow[]{r}{} \arrow[]{d}{} & Y' \arrow[]{r}{} \arrow[]{d}{} & X' \arrow[]{d}{} \\
Z \arrow[]{r}{} & Y \arrow[]{r}{} & X.
\end{tikzcd}
\end{equation*}
Now:
\begin{enumerate}
\item if both of the small squares are homotopy Cartesian, then so is the large square;
\item if both the large square and the right small square are homotopy Cartesian, then so is the leftmost small square.
\end{enumerate}
\end{prop}

Derived fibre products commute with truncations:

\begin{prop}\label{FibreProductAndTruncations}
We have a natural equivalence
$$t(X \times^\R_S Y) \to tX \times_{tS} tY.$$
\end{prop}

One consequence is that it is usually rather easy to see when a derived tensor product yields the empty scheme. Moreover, virtually all sensible classes of morphisms are closed under homotopy pullbacks; among them are all the three special classes we have already defined.
\begin{prop}
Let us have a homotopy Cartesian square
\begin{equation*}
\begin{tikzcd}
Y' \arrow[]{r}{} \arrow[]{d}{f'} & X' \arrow[]{d}{f} \\
Y \arrow[]{r}{} & X
\end{tikzcd}
\end{equation*}
Now:
\begin{enumerate}
\item if $f$ is a closed embedding, then so is $f'$;
\item if $f$ is a proper, then so is $f'$;
\item if $f$ is a quasi-projective, then so is $f'$.
\end{enumerate}
\end{prop}
\begin{proof}
The claims $(1)$ and $(2)$ follow immediately from Proposition \fref{FibreProductAndTruncations} and definitions \fref{ClosedEmbeddingDef} and \fref{ProperDef} respectively. Suppose $f$ factors as the composition of $i: X' \hookrightarrow U \times X$ and $p_2: U \times X \to X$, where $U$ is an open subscheme of $\Proj^n$. The square
\begin{equation*}
\begin{tikzcd}
U \times Y \arrow[]{r}{} \arrow[]{d}{p_2} & U \times X \arrow[]{d}{p_2} \\
Y \arrow[]{r}{} & X
\end{tikzcd}
\end{equation*}
is derived Cartesian which easily by applying Proposition\fref{FibreProductAndTruncations} $(2)$ to the diagram

\begin{equation*}
\begin{tikzcd}
U \times Y \arrow[]{r}{} \arrow[]{d}{p_2} & U \times X \arrow[]{r}{} \arrow[]{d}{p_2} & U \arrow[]{d}{} \\
Y \arrow[]{r} & X \arrow[]{r} & pt
\end{tikzcd}
\end{equation*}
As the pullback of $i$ is a closed embedding, Proposition \fref{FibreProductAndTruncations} $(1)$ now expresses $f'$ as a composition of the desired form.
\end{proof}

We close this subsection with an important warning.

\begin{war}
Throughout the introductory section on derived algebraic geometry, in order not to overburden the exposition, we sometimes make simplifications in definitions and statements of theorems. We have tried to warn the reader in each instance in a separate remark or a warning. For the purposes of this article, however, the reader should not worry too much, as all the extra assumptions automatically hold for quasi-projective derived schemes.

There are two main simplifications. First of all, sometimes we need to make some finiteness assumptions on the truncation $tX$ of a derived scheme $X$. Quasi-projective schemes being finite type over $k$, they enjoy very strong finiteness properties that are enough for all the results. Second of all, we will sometimes have to make the assumption that the derived scheme $X$ is quasi-compact and quasi-coherent. As both of these can be checked on the truncation, both the conditions hold whenever $X$ is quasi-projective. 
\end{war} 

\subsection{Quasi-coherent sheaves}

As in classical algebraic geometry, one can learn a lot from a derived scheme $X$ by studying its quasi-coherent sheaves. As we are only going to work with very special kinds of quasi-coherent sheaves (vector bundles), we are not going to give a precise definition here, but refer the interested reader to Chapter 2 of \cite{Lur3}. For the purposes of this paper, it is enough to state that a \emph{quasi-coherent sheaf} on $X$ is 
\begin{enumerate}
\item a sheaf (in the $\infty$-categorical sense) of spectrum objects of simplicial Abelian groups (unbounded chain complexes of Abelian groups) on the underlying topological space $X_\mathrm{top}$;
\item together with an action of the structure sheaf $\OO_X$ (an \emph{$\OO_X$-module});
\item on every affine open $\Spec(A) \subset X$, $\cF \vert_{\Spec(A)}$ is equivalent to the sheaf associated to some $A$ dg-module $M$.
\end{enumerate}
The $\infty$-category of quasi-coherent sheaves on $X$ is denoted by $\QCoh(X)$. Given two quasi-coherent sheaves $\cF$ and $\cE$, the mapping space $\Hom_X(\cF, \cE)$ is the space of morphisms between the underlying $\OO_X$-modules. We note that when $X$ is a classical scheme, then $\QCoh(X)$ is (the $\infty$-category associated to) the unbounded derived dg-category of $X$. 

The category $\QCoh(X)$ has a full subcategory $\Perf(X)$ --- the category of \emph{perfect objects} or \emph{perfect complexes}. It is the full subcategory containing those quasi-coherent sheaves, which locally on $\Spec(A)$ become equivalent to an iterated mapping cone of shifts of projective $A$-modules (locally free modules). Note that whenever $X$ is classical, being a perfect object just means that locally $\cF$ should be quasi-isomorphic to a finite complex of vector bundles, and if $X$ defined over a field of characteristic 0, $\Perf(X)$ is the derived dg-category of perfect complexes on $X$. 

\begin{war}
In bad situations there are multiple different definitions for perfectness which agree in our case of interest. The above definition is correct (in the sense that the theorems we list below remain true) at least when the derived scheme $X$ is quasi-compact and quasi-separated, e.g., when $X$ is quasi-projective.
\end{war}

Both $\QCoh(X)$ and $\Perf(X)$ are \emph{stable}. The interested reader can consult the definition from the first chapter of \cite{Lur2}; we are going to restrict ourselves merely to listing some practical consequences of this. Namely:
\begin{enumerate}
\item both the above categories have zero objects $0$, which respect the inclusion $\Perf(X) \subset \QCoh(X)$;
\item both $\Perf(X)$ and $\QCoh(X)$ admit small limits and colimits;
\item a square
\begin{equation*}
\begin{tikzcd}
A \arrow[]{r} \arrow[]{d} & B \arrow[]{d} \\
C \arrow[]{r} & D
\end{tikzcd}
\end{equation*}
in either of the above categories is a pullback square if and only if it is a pushout square.
\end{enumerate}

A triangle of morphisms $A \to B \to C$ in a stable $\infty$-category is called a \emph{(co)fibre sequence} if it fits into a (co)Cartesian square
\begin{equation*}
\begin{tikzcd}
A \arrow[]{r} \arrow[]{d} & B \arrow[]{d} \\
0 \arrow[]{r} & C
\end{tikzcd}
\end{equation*}
In such a situation $C$ is called the \emph{cofibre} of the morphism $A \to B$ and $A$ is called the \emph{fibre} of the morphism $B \to C$. In our examples, the cofibre is modeled by the \emph{mapping cone construction} familiar from elementary homological algebra. 

Chain complexes have the natural shift functors $[n]$ sending a chain complex $C_{\bullet}$ to the chain complex $C_{\bullet - n}$. The analogue in the setting of stable $\infty$-categories is taking the cofibre of the unique morphism $\cF \to 0$ giving rise to the \emph{suspension functor} $\cF \mapsto \cF[1]$. It has an inverse $\cF \mapsto \cF[-1]$ obtained by taking the fibre of $0 \to \cF$, and one can define using composition shift functors $\cF \mapsto \cF [n]$ for all $n \in \Z$. The axioms of stable $\infty$-categories provide canonical equivalences $\cF[n][m] \simeq \cF[n+m]$. In much of the literature, $\cF[1]$ is denoted by $\Sigma \cF$ owing to the fact that the suspension $\Sigma X$ of a topological space $X$ can be defined as the homotopy pushout of the diagram
\begin{equation*}
\begin{tikzcd}
X \arrow[]{r} \arrow[]{d} & pt \\
pt.
\end{tikzcd}
\end{equation*}
Similarly, $\cF[-1]$ is often denoted by $\Omega \cF$, where $\Omega$ stands for the space of based loops.

\begin{defn}\label{HomotopySheaves}
Any object $\cF$ of $\QCoh(X)$ (and hence $\Perf(X)$) has naturally associated \emph{homotopy sheaves} (sheaves in the classical sense) $\pi_n(\cF)$ for all integers $n \in \Z$. As $k$-linear sheaves, these can be identified with sheaf associated to the presheaf that on an affine open $\Spec(A) \subset X$ takes as the value the $n^{th}$ homology of the underlying chain complex of the dg-$A$ module corresponding to $\cF \vert_{\Spec(A)}$. The sheaves $\pi_n(\cF)$ come naturally equipped with the structure of a quasi-coherent sheaf on the truncation $tX$.

A quasi-coherent sheaf $\cF$ is called \emph{connective} if all the negative homotopy sheaves vanish.
\end{defn}

\begin{war}\label{HomotopyOfGlobalSectVsGlobalSectOfHomotopy}
By definition an element $\cF$ has an underlying sheaf of chain complexes of $\Z$-modules, and we can therefore talk about its global sections $\Gamma(X; \cF)$. As this is a chain complex, it 
also has natural homotopy groups $\pi_n(\Gamma(X; \cF))$, which are defined as its homology groups (this naming convention has its roots in the famous Dold-Kan correspondence). It is \emph{not} true in general that the $k$-vector spaces $\pi_n(\Gamma(X; \cF))$ and $\Gamma_\mathrm{classical}(tX; \pi_n(\cF))$ are isomorphic. Indeed, for a simple example, consider a classical scheme $X$ together with a classical quasi-coherent sheaf $\cF$. In this situation all the nonzero homotopy sheaves vanish, yielding to us the equality $\Gamma_\mathrm{classical}(tX; \pi_n(\cF)) = 0$ for all $n \not = 0$. However, the homotopy groups $\pi_{-n}(\Gamma(X; \cF))$ are naturally identified with the sheaf cohomology groups $H^{n}(X; \cF)$ of $\cF$. As the latter are nontrivial in general, we obtain the desired counterexample.
\end{war}

We are mainly interested in the following special case of an object of $\Perf(X)$.
\begin{defn}\label{VectorBundle}
A \emph{vector bundle of rank $r$} is a quasi-coherent sheaf $E$ that is locally equivalent to the free sheaf $\OO_X^{\oplus r}$ of rank $r$. A vector bundle of rank $1$ is called a \emph{line bundle}. 
\end{defn}
\begin{rem}
When $X$ is a classical scheme, then the above notion of a vector bundle coincides with the usual notion of a vector bundle on $X$. 
\end{rem}

As the homotopy theory of vector bundles is very special,the following definition is sensible:

\begin{defn}
A morphism $f: E \to F$ of vector bundles on $X$ is \emph{surjective} if the truncation $\pi_0(f): \pi_0(E) \to \pi_0(F)$ is a surjective morphism of sheaves on $tX$.
\end{defn} 

\begin{rem}\label{HomologicalAlgebraOfVectorBundles}
If $E$ is a vector bundle of non zero rank, then the shifts $E[n]$ are not vector bundles for any $n \not = 0$. The fibre $K$ of a morphism $f: E \to F$ of vector bundles is a vector bundle if and only if $f$ is a surjection.
\end{rem}

Before going further we will record the following subtlety for the convenience of the non-expert reader.

\begin{rem}[The spave of global sections]
Given a quasi-coherent sheaf $\cF$ on a derived scheme $X$, we have, as in Remark \fref{HomotopyOfGlobalSectVsGlobalSectOfHomotopy} the $k$-linear chain complex $\Gamma(X; \cF)$ of global sections of $\cF$. Moreover, as noted in that remark, the chain complex often has nontrivial homotopy groups (defined as the homology groups) of negative degree.

The famous Dold-Kan correspondence gives a way to naturally associate a space $\abs{C_{\bullet}}$ to a connective (meaning no nonzero homotopy groups) chain complex $C_{\bullet}$. Moreover the homotopy groups $\pi_n \abs{C_\bullet}$ are naturally isomorphic to the homology groups $H_n(C_\bullet)$. Given a quasi-coherent sheaf $\cF$ on $X$ we define the \emph{space of global sections} $\abs{\Gamma(X; \cF)}$ of $\cF$ as the space $\abs{\tau_{\geq 0} \bigl( \Gamma(X; \cF) \bigr)}$, where $\tau_{\geq 0}$ is the \emph{canonical truncation functor}. When $\cF$ is a classical quasi-coherent sheaf on a classical scheme, the space $\abs{\Gamma(X; \cF)}$ is naturally equivalent to the discrete space of global sections $\Gamma_\mathrm{classical}(X; \cF)$.
\end{rem}

\begin{prop}
The space of morphisms $\Hom_X(\OO_X, \cF)$ of quasi-coherent sheaves is naturally equivalent to $\abs{\Gamma(X; \cF)}$.
\end{prop}

We can then make the following definition:

\begin{defn}
A vector bundle $E$ on $X$ is \emph{globally generated} if there exists a surjective morphism $\OO_X^{\oplus n} \to E$ of vector bundles.
\end{defn}

As surjectivity by definition depends only on what happens on the zeroth homotopy groups, a vector bundle $E$ is globally generated if and only if there exists global sections $s_1,...,s_n$ so that their truncations $ts_1, ..., ts_n$ generate $\pi_0(E)$ on $tX$ in the classical sense. We note that this is \emph{not} the same as assuming the truncation to be globally generated as not all global sections can be lifted from the truncation.

The category $\Perf(X)$ admits internal mapping objects (see \cite{Lur3} Section 6.5.3, especially Proposition 6.5.3.6).

\begin{defn}\label{InternalMappingOfPerf}
Let $X$ be a derived scheme. Now the $\infty$-category $\Perf(X)$ of perfect complexes admits internal mapping objects. In other words, there is a bifunctor $\iHom_X(-,-)$ so that for $\cF, \cE, \mathcal G \in \Perf(X)$, we have natural equivalences 
$$\Hom_X(\cF, \iHom_X(\cE, \mathcal{G})) \simeq \Hom_X(\cF \otimes \cE, \mathcal{G}).$$
In particular, when $\cF \simeq \OO_X$, we obtain the equivalence 
$$\abs{\Gamma(X, \iHom_X(\cE, \mathcal{G}))} \simeq \Hom_X(\cE, \mathcal{G}).$$
\end{defn}

We can now define the dual of a perfect complex. 

\begin{defn}
Given an object $\cF$ of $\Perf(X)$, we denote by $\cF^\vee$ its \emph{dual sheaf}, which is defined as the mapping object $\iHom_X(\cF, \OO_X)$.
\end{defn}

The most important properties of the dual construction are summarized in the following proposition.

\begin{prop}\label{PropertiesOfDuals}
The construction $\cF \mapsto \cF^\vee$ satisfies the following properties:
\begin{enumerate}
\item The space of global sections of $\cF^\vee$ is naturally identified with $\Hom_X(\cF, \OO_X)$.
\item The double dual $\cF^{\vee \vee}$ is naturally equivalent to $\cF$.
\item If $E$ is a vector bundle of rank $r$, then so is $E^\vee$.
\item If $\Li$ is a line bundle, then $\Li \otimes \Li^\vee \simeq \OO_X$.
\item The internal mapping object $\underline{\Hom}_X(\cF, \cE)$ of $\cF, \cE \in \Perf(X)$ is naturally identified with $\cF^\vee \otimes \cE$.
\item If $X$ is a classical scheme, and $\cF$ is a perfect complex on $X$ (e.g., a coherent sheaf of finite projective dimension), then $\cF^\vee$ is the \emph{derived dual} of $\cF$. 
\end{enumerate}
\end{prop}

Finally, we need to introduce the \emph{Tor-amplitude} of a perfect complex $\cF$.

\begin{defn}
Let $\cF$ be a perfect complex on a derived scheme $X$. We say that $\cF$ \emph{has Tor-amplitude in $[n,m]$} (here $n \leq m$) if for all \emph{discrete} quasi-coherent sheaves $\cE$ (i.e., $\cE$ can be identified with the homotopy sheaf $\pi_0(\cE)$), the homotopy sheaves $\pi_k(\cF \otimes \cE)$ vanish whenever $k \in \Z$ is not in the closed interval $[n,m]$.
\end{defn}

It is known that a perfect complex always has Tor-amplitude in some finite interval (see \cite{Lur2} Proposition 7.2.4.23 (4)). Moreover, vector bundles can be identified as the perfect complexes with Tor-amplitude in $[0,0]$ (see \cite{An2} Lemma 5.3).

\subsection{Pullbacks and pushforwards of quasi-coherent sheaves}

Like in classical algebraic geometry, given a morphism $f: X \to Y$ of derived schemes, the pushforward functor $f_*$ of sheaves gives naturally rise to a functor
$$f_*: \QCoh(X) \to \QCoh(Y)$$
which admits a left adjoint 
$$f^*: \QCoh(Y) \to \QCoh(X)$$
and both of these functors are naturally functorial in compositions of morphisms of derived schemes.  Whenever $X$ and $Y$ are classical schemes, then these functors coincide with the \emph{derived functors} of pushforward and pullback, which are defined on the level of derived categories. Most of the claims of this section are taken from chapters 2 and 3 of \cite{Lur3}.

\begin{war}
Strictly speaking, we need to assume that the morphism $f$ above is \emph{relatively scalloped} for the above definition to work (see \cite{Lur3} Definition 2.5.4.1 and Proposition 2.5.4.3). However, by Theorem 3.4.2.1 of \emph{loc. cit.} this is true for any morphism between derived schemes that are quasi-compact and quasi-coherent, which is the only case we are interested in the paper.
\end{war}

Moreover, the following proposition follows more or less directly from the definition of the pushforward.

\begin{prop}
Let $X$ be a derived scheme, let $\pi: X \to \Spec(A)$ be morphism to an affine scheme, and let $\cF$ be a quasi-coherent sheaf on $X$. Now the $k$-linear chain complexes $\pi_*(\cF)$ and $\Gamma(X;\cF)$ are equivalent.
\end{prop} 
\begin{proof}
As the target $\Spec(A)$ is affine, the sheaf $\pi_*(\cF)$ is completely and naturally determined by its global sections (\cite{Lur3} Proposition 2.4.1.4 for a more general statement). But by definition
$$\Gamma(\Spec(A); \pi_*(\cF)) := \Gamma(X; \cF),$$
proving the claim.
\end{proof}

We will also need the following two generalizations of classical results.

\begin{prop}[Push-pull formula, \cite{Lur3} Proposition 2.5.4.5]\label{PushPull}
Let us have a derived Cartesian square
\begin{equation*}
\begin{tikzcd}
X' \arrow[]{r}{f'} \arrow[]{d}{g'} & Y' \arrow[]{d}{g} \\
X \arrow[]{r}{f} & Y
\end{tikzcd}
\end{equation*}
of derived schemes. Then there is a natural equivalence
$$g^* f_* \simeq g'^* f'_*$$
of functors $\QCoh(X) \to \QCoh(Y').$
\end{prop}

\begin{prop}[Projection formula, \cite{Lur3} Remark 3.4.2.6]\label{ProjectionFormula}
Given a morphism $f: X \to Y$ of derived schemes and quasi-coherent sheaves $\cF \in \QCoh(X)$ and $\mathcal{G} \in \QCoh(Y)$, then we have a natural equivalence
$$f_*(\cF) \otimes \mathcal{G} \simeq f_* \bigl( \cF \otimes f^*(\mathcal{G}) \bigr)$$
of quasi-coherent sheaves.
\end{prop}

\subsection{Quasi-coherent sheaves on quasi-projective derived schemes}

As already promised, it is much nicer to work with quasi-coherent sheaves on quasi-projective derived schemes. Recall from \cite{An2} Definition 4.2 that a line bundle $\Li$ on $X$ is \emph{ample} if its truncation $\pi_0(\Li)$ is ample in the classical sense on $X$. The following is essentially \emph{loc. cit.} Theorem 4.7 (it is stated for $A = k$ is a field of characteristic 0, but it also works in this generality):

\begin{thm}
Let $A$ be a Noetherian ring and let $X$ be a derived $A$-scheme whose truncation is of finite type over $A$. Then $X$ is quasi-projective over $A$ if and only if it admits an ample line bundle.
\end{thm}

Moreover, every vector bundle can be twisted to be globally generated (this is the Corollary 4.4 of \emph{loc. cit.} as vector bundles are \emph{strong}):

\begin{thm}
Let $X$ be a derived scheme and $\OO_X(1)$ an ample line bundle on $X$. Now, given a vector bundle $E$ on $X$, then for $n \gg 0$ the vector bundle
$$E(n) := E \otimes \OO_X(1)^{\otimes n}$$
is globally generated.
\end{thm}
The above theorem played a fundamental role in the construction of Chern classes in \cite{An}. It is also the reason why we can often prove theorems about vector bundles by reducing first to the globally generated case.

There is also a relative notion of ampleness:

\begin{defn}
Suppose $f: X \to Y$ is a morphism of derived schemes. A line bundle $\Li$ on $X$ is \emph{(relatively) ample over $Y$} (\emph{$f$-ample}) if for every affine open $\Spec(A) \subset Y$ (equivalently every affine open $\Spec(A) \subset Y$ in some affine open cover of $Y$), the restriction $\Li \vert_{f^{-1} \Spec(A)}$ is ample.
\end{defn}

We will need the following theorem (see \cite{An2} Theorem 4.11) in order to show that the derived blow up of a quasi-projective scheme is quasi-projective.

\begin{thm}\label{QProjAndRelativelyAmple}
Suppose $X \to Y$ is a morphism of derived $A$-schemes (with finite type truncations) with $Y$ quasi-projective. Then $X$ is quasi-projective if and only if it admits a relatively ample line bundle over $Y$.
\end{thm}

Finally, perfect complexes on quasi-projective derived schemes admit global resolutions by vector bundles (see the discussion succeeding Lemma 5.3 in \cite{An2} for the precise result).

\begin{prop}\label{ResolutionOfPerfectComplexes}
Let $\cF$ be a perfect complex on a quasi-projective derived scheme $X$ having Tor-amplitude in $[0,d]$. Then $\cF$ admits a resolution $F_d \to F_{d-1} \to \cdots \to F_1 \to F_0$ by vector bundles.
\end{prop}

\subsection{Symmetric algebras of sheaves, geometric vector bundles and derived vanishing loci}

Like in classical algebraic geometry, symmetric powers and symmetric algebras (or rather their global spectra) play an important role in derived algebraic geometry. Let us first recall the definition (cf. \cite{Lur3} Construction 25.2.2.6).

\begin{defn}\label{SymmetricAlgebra}
Suppose $X$ is a derived scheme and $\cF$ is a connective quasi-coherent sheaf on $X$. Now $\cS pec (\Sym^*_X(\cF))$ is an $X$ scheme so that given an $X$-scheme $\pi_S: S \to X$, the space of $X$-maps $S \to \cS pec (\Sym^*_X(\cF))$ is naturally equivalent to the space $\Hom_S(\pi_S^* \cF, \OO_S)$ of morphisms of quasi-coherent sheaves on $S$.
\end{defn}

\begin{rem}
Several remarks are in order to explain the differences of the above construction with \cite{Lur3}. First of all, as we are working over a field $k$ of characteristic 0, the constructions of symmetric powers and their derived versions $\LSym^*$ agree. Hence, in order to lighten the exposition, we are going to denote the symmetric algebra by $\Sym^*$. Secondly, by the adjointness of push-forward and pullback, we have a natural equivalence
\begin{equation*}
\Hom_S(\pi_S^* \cF, \OO_S) \simeq \Hom_X( \cF, \pi_{S*} \OO_S),
\end{equation*}
where the right hand side is what is used in \emph{loc. cit.}
\end{rem}

The above theorem allows us to make precise in derived algebraic geometry the double life of vector bundles: that vector bundles can be regarded at the same time both as quasi-coherent sheaves and as derived schemes over the base scheme of interest.

\begin{defn}
Let $E$ be a vector bundle on a derived scheme $X$. Now the $X$-scheme $\sSpec(\Sym^*_X(E^\vee))$ is called the \emph{geometric vector bundle associated to $E$}. By abuse of notation, we are going to denote it still by $E$
\end{defn}

The following theorem, which also explains why in the above definition we are taking the symmetric algebra of the \emph{dual} $E^\vee$ rather than the original vector bundle $E$, is rather fundamental, although it follows directly from Definition \fref{SymmetricAlgebra}. 

\begin{thm}
Let $X$ be a derived scheme and let $E$ be a vector bundle on $X$. Now the space of $X$-morphisms $X \to \sSpec(\Sym_X^*(E^\vee))$, i.e., the space of sections of the structure morphism $\sSpec(\Sym_X^*(E^\vee)) \to X$, is canonically identified with the underlying space of the global sections of $E$.
\end{thm}
\begin{proof}
Indeed, the space of such morphism is canonically identified with $\Hom_X(E, \OO_X)$, so the claim follows from Proposition \fref{PropertiesOfDuals} $(1)$ and $(2)$.
\end{proof}

From now on, we won't make a distinction between a vector bundle $E$ and its geometric version, nor do we distinguish between a global section of $E$ as a quasi-coherent sheaf and a section of the structure morphism $E \to X$.

Suppose we have a vector bundle $E$ on a derived scheme $X$, and $s$ is a global section of $E$. The \emph{derived vanishing locus} $V(s)$ of the section $s$ is usually defined as the derived fibre product
\begin{equation*}
\begin{tikzcd}
V(s) \arrow[]{r}{} \arrow[]{d}{} & X \arrow[]{d}{s} \\
X \arrow[]{r}{0} & E
\end{tikzcd}
\end{equation*}
where we have denoted by $0$ the zero section of $E$. Because the defining property of the symmetric algebra as defined in \cite{Lur3} Construction 25.2.2.6 is as a left adjoint, it commutes with colimits. Hence, noting that the zero section $X \to E$ is the image of the morphism $E^\vee \to 0$ of quasi-coherent sheaves on $X$, we obtain the following result:

\begin{prop}\label{DerivedVanishingLocus}
Let everything be as above. Then the derived vanishing locus $V(0)$ of the zero section of $E$ is naturally identified with $\sSpec(\Sym_X^*(E^\vee[-1]))$.
\end{prop}

In the special case where we have a section $s$ of a line bundle $\Li$, we call $V(s)$ the \emph{virtual Cartier divisor} associated to $s$. We record here the immediate observation that a line bundle $\Li$ is globally generated if and only if its global sections give a \emph{base point free linear system}, i.e., no closed point of the underlying space $X_\mathrm{top}$ is contained in all the virtual Cartier divisors $V(s)$, $s$ being a section of $\Li$.

\subsection{Projective bundles}

The derived version of the classical construction of a projective bundle will be important later, as the exceptional divisor of a derived blow up takes the form of a projective bundle over the center of the blow up. Let us begin with a definition (cf. \cite{Khan2} Section 3.1).

\begin{defn}
Let $X$ be a derived scheme and $E$ a vector bundle on $X$. Then the \emph{projective bundle} $\Proj(E)$ is the $X$-scheme so that given an $X$-scheme $\pi_S: S \to X$, the space of $X$-morphisms $S \to \Proj(E)$ is naturally equivalent to the space of surjections
$$\pi_S^* E^\vee \to \Li,$$ 
where $\Li$ is a line bundle on $S$. The identity morphism $\Proj(E) \to \Proj(E)$ induces the \emph{canonical surjection} $\pi^* E^\vee \to \OO_{\Proj(E)}(1)$ on $\Proj(E)$, where $\pi$ is the structure morphism $\Proj(E) \to X$. If it should cause no confusion, we sometimes denote $\OO_{\Proj(E)}(1)$ by $\OO(1)$.
\end{defn}

As a special case of the the above definition we obtain $\Proj^{n-1}$ as the projective bundle associated to $k^{\oplus n}$ over $\Spec(k)$. Hence, given a line bundle $\Li$ and global sections $s_0,...,s_{n-1}$ that generate $\Li$, we obtain a morphism $X \to \Proj^{n-1}$. Moreover, multiplying these sections by a common (invertible) scalar does not change the map.

Projective bundles admit the following functorial property.

\begin{defn}\label{ProjectivizedInclusion}
Suppose $E \to F$ is a morphism of vector bundles on $X$ so that the dual morphism $F^\vee \to E^\vee$ is surjective. The induced surjection
$$\pi_{\Proj(E)}^* F^\vee \to \pi_{\Proj(E)}^* E^\vee \to \OO_{\Proj(E)}(1)$$
induces a morphism $\Proj(E) \to \Proj(F)$ over $X$ called the \emph{projectivization} of the inclusion $E \to F$ (such a morphism must necessarily be an injection on $\pi_0$). Its truncation is the classical projectivized inclusion, which is a closed embedding, and therefore the derived version is a closed embedding as well.
\end{defn}

Moreover, as in classical case, tensoring the vector bundle $E$ with a line bundle $\Li$ does not change the $X$-scheme $\Proj(E)$, but it does change the universal line bundle $\OO(1)$. More precisely, we have the following.

\begin{prop}\label{ProjectiveBundlesTensoredWithLineBundles}
Let $X$ be a derived scheme, $E$ a vector bundle on $X$, and $\Li$ a line bundle on $X$. Now the natural surjection $\pi^*_{\Proj(E)}(E^\vee \otimes \Li^\vee) \to \pi_{\Proj(E)}^*(\Li^\vee) \otimes \OO_{\Proj(E)}(1)$ on $\Proj(E)$ induces an equivalence $\Proj(E) \to \Proj(E \otimes \Li)$. 

In other words we have a natural equivalence $\Proj(E) \simeq \Proj(E \otimes \Li)$ of derived schemes over $X$ and a natural equivalence $\OO_{\Proj(E \otimes \Li)}(1) \simeq \pi^*{\Li^\vee} \otimes\OO_{\Proj(E \otimes \Li)}(1)$ of line bundles.
\end{prop}
\begin{proof}
The morphism $\Proj(E) \to \Proj(E \otimes \Li)$ is clearly an equivalence at least when $\Li$ is trivial. On the other hand, every line bundle is locally trivial, proving that $\Proj(E) \to \Proj(E \otimes \Li)$ is an equivalence Zariski locally, and hence globally.
\end{proof}

The final goal of this subsection is to express $\Proj(E) \hookrightarrow \Proj(E \oplus F)$ as a derived vanishing locus of a section of a vector bundle on $\Proj(E \oplus F)$. Before that, however, we need to record the following results concerning pushforwards along the structure morphism.

\begin{prop}\label{PushforwardsInProjectiveBundle}
Let $X$ be a derived scheme, and let $\pi: \Proj(E) \to X$ be a projective bundle over $X$. Then
\begin{enumerate}
\item given a vector bundle $F$ on $X$, the natural unit morphism
$$F \to \pi_* \pi^* F$$
is an equivalence;

\item applying the pushforward functor $\pi_*$ to the natural surjection $\pi^* E^\vee \to \OO(1)$, we obtain a natural equivalence
$$\pi_* \pi^* E^\vee \to \pi_* \OO(1).$$
\end{enumerate}
Combining parts $(1)$ and $(2)$, we obtain a natural equivalence $E^\vee \to \pi_* \OO(1)$.
\end{prop}
\begin{proof}
In both of the cases, it is enough to check that the globally given morphism is locally an equivalence, and therefore we can reduce to the case where both $E$ and $F$ are trivial vector bundles. In order to prove the first claim, we note that it follows from the push-pull formula (Proposition \fref{PushPull}) applied to the derived Cartesian square
\begin{equation*}
\begin{tikzcd}
\Proj^n \times X \arrow[]{r} \arrow[]{d}{} & \Proj^n \arrow[]{d} \\
X  \arrow[]{r} & \Spec(k)
\end{tikzcd}
\end{equation*}
with $\OO_{\Proj^n}^{\oplus r}$ in the top right corner, that the sheaf $\pi_* \pi^* \OO_X^{\oplus r}$ is a (trivial) vector bundle. Hence, $F \to \pi_* \pi^* F$ is a map of vector bundles, and it is an equivalence if and only if its truncation is (vector bundles are strong, see \cite{TV} Definition 2.2.2.1). But the latter claim follows from well known classical results.

To prove the second claim, consider again the homotopy Cartesian diagram
\begin{equation*}
\begin{tikzcd}
\Proj^n \times X \arrow[]{r}{p} \arrow[]{d}{\pi} & \Proj^n \arrow[]{d}{\pi} \\
X \arrow[]{r}{} & \Spec(k).
\end{tikzcd}
\end{equation*}
As the natural surjection $\pi^* \OO_X^{\oplus n+1} \to \OO_{\Proj^n \times X}(1)$ is identified as the pullback of the natural surjection $\OO_{\Proj^n}^{\oplus n+1} \to \OO_{\Proj^n}(1)$ via $p$, this claim too follows from combining well known classical results with the push-pull formula \fref{PushPull}.
\end{proof}

\begin{rem}
The first part of the above claim holds in greater generality. Indeed, knowing it for vector bundles, it follows immediately for all perfect complexes as both $\pi_*$ and $\pi^*$ preserve finite colimits.
\end{rem}

Consider now the projective bundle $\pi: \Proj(E \oplus F) \to X$ and the vector bundle $\pi^*(F) \otimes \OO(1)$ on it. Then, using functoriality of pushforwards and the projection formula \fref{ProjectionFormula}, we have a natural identification of global sections of $\pi^*(F) \otimes \OO(1)$ with the global sections of its pushforward $F \otimes (E^\vee \oplus F^\vee)$. On the other hand, as 
$$ \abs{\Gamma(X; F \otimes (E^\vee \oplus F^\vee))} \simeq \Hom_X(F^\vee, E^\vee \oplus F^\vee)$$
by Proposition \fref{PropertiesOfDuals} part $(5)$, we obtain a global section $s$ of $\pi^*(F) \otimes \OO(1)$ corresponding to the canonical inclusion $F^\vee \hookrightarrow E^\vee \oplus F^\vee$.

\begin{prop}\label{ProjectivizedEmbeddingAsCompleteIntersection}
Let everything be as above. Then the inclusion of the derived vanishing locus $i: V(s) \hookrightarrow \Proj(E \oplus F)$ can be identified with the projectivized embedding $j: \Proj(E) \hookrightarrow \Proj(E \oplus F)$.
\end{prop}
\begin{proof}
The proof consists of two parts: we have to first find a morphism $\Proj(E) \to V(s)$ over $\Proj(E \oplus F)$ and then show that it is an equivalence of derived schemes.
\begin{enumerate}
\item By a global version of \cite{Khan} Lemma 2.3.5, the space of $\Proj(E \oplus F)$-morphisms $\Proj(E) \to V(s)$ is equivalent to the space of paths $\alpha: j^* (s) \sim 0$ inside the space of global sections of $j^*(\pi_{\Proj(E \oplus F)}^*F \otimes \OO(1)) \simeq \pi_{\Proj(E)}^* (F) \otimes \OO_{\Proj(E)}(1)$, which is equivalent to $\Hom_X(F^\vee, E^\vee)$ via a calculation identical to the one preceding the proposition. 

We note first that the map
$$\eta : \pi_{\Proj(E \oplus F)}^*(F) \otimes \pi_{\Proj(E \oplus F)}^*(E^\vee \oplus F^\vee) \to \pi_{\Proj(E \oplus F)}^*(F) \otimes \OO_{\Proj(E \oplus F)}(1)$$
induced by the canonical surjection induces an equivalence on global sections, as does the pullback morphism
$$j^*: \Gamma(\Proj(E \oplus F); \pi_{\Proj(E \oplus F)}^*(F \otimes (E^\vee \oplus F^\vee))) \to \Gamma(\Proj(E); \pi_{\Proj(E)}^* (F \otimes (E^\vee \oplus F^\vee))).$$
By the defining property of the linearized embedding $j$, there is a natural factorization of the pullback $j^* \eta$ of $\eta$ expressing it as the composition
\begin{align*}
\pi_{\Proj(E)}^*(F) \otimes \pi_{\Proj(E)}^*(E^\vee \oplus F^\vee) &\to \pi_{\Proj(E)}^*(F) \otimes \pi_{\Proj(E)}^*(E^\vee) \\
&\stackrel{\eta'}{\to} \pi_{\Proj(E)}^*(F) \otimes \OO_{\Proj(E)}(1)
\end{align*} 
and the induced map on global sections is canonically identified with the map
$$j^*: \Gamma(\Proj(E \oplus F); \pi^*_{\Proj(E \oplus F)}(F) \otimes \OO_{\Proj(E \oplus F)}(1)) \to \Gamma(\Proj(E); \pi^*_{\Proj(E)}(F) \otimes \OO_{\Proj(E)}(1)).$$
Applying $\pi_{\Proj(E)*}$ to $j^* \eta$, using projection formula, and recalling that $\eta'$ induces an equivalence on global sections, we see that the morphism
$$\pi_{\Proj(E)*}j^* \eta: F \otimes (E^\vee \oplus F^\vee) \to F \otimes E^\vee$$
is the one induced by the canonical projection $E^\vee \oplus F^\vee \to E^\vee$. Hence, using the identification given before the proposition, the map
$$j^*: \Hom_X(F^\vee, E^\vee \oplus F^\vee) \to \Hom_X(F^\vee, E^\vee)$$
is given by composing with the canonical projection $E^\vee \oplus F^\vee \to E^\vee$, and the universal property of direct sums provides a canonical path $\alpha: j^*(s) \simeq 0$ (the constant path), therefore providing also a $\Proj(E \oplus F)$-morphism $\Proj(E) \to V(s)$.

\item As usual we can check being an equivalence locally, and we can therefore assume that both $E$ and $F$ are trivial vector bundles $\OO_X^{\oplus r}$ and $\OO_X^{\oplus s}$ respectively. Denoting by $e_1,...,e_{r+s}$ the usual basis for $\OO_X^{\oplus r} \oplus \OO_X^{\oplus s}$, and by $x_1,...,x_{r+s}$ its dual basis, we see that the section $s$ is given by 
$$e_{r+1} \otimes x_{r+1} + \cdots + e_{r+1} \otimes x_{r+1} \in \Gamma\bigl(\Proj(\OO_X^{\otimes r} \oplus \OO_X^{\otimes r}); \OO_X^{\oplus r} \otimes \OO(1)\bigr).$$
Hence the derived vanishing locus $V(s)$ is just the derived vanishing locus of the sections $x_{r+1},...,x_{r+s}$ of $\OO(1)$, so it is at least abstractly equivalent to $\Proj(\OO_X^{\oplus r})$.

In classical algebraic geometry, closed embeddings have no nontrivial endomorphism, and therefore the induced map $t\Proj(E) \to tV(s)$ on truncations must be an isomorphism of schemes. Moreover, as both of the derived schemes are smooth over $X$, and of the same relative virtual dimension, we conclude that $\Proj(E) \to V(s)$ is a quasi-smooth embedding of virtual codimension $0$ (see Theorem \fref{QSProperties} in the following subsection). On the other hand, such a morphism is necessarily an equivalence, for example because such an embedding is locally cut out by 0 equations, see \cite{Khan} Section 2.
\end{enumerate}
\end{proof}

\subsection{Cotangent complex and quasi-smooth morphisms}\label{CotangentComplexSection}
The cotangent complex is perhaps the single most important object in derived algebraic geometry. Via derived deformation theory, it allows great control over a derived scheme whose truncation and cotangent complex are understood (for example various derived moduli spaces). Even though we are not going to make a serious use of the more advanced properties of cotangent complex, even the simplest ones help us a long way in our eventual construction of the bivariant theories of interest. Let us begin with a definition.

\begin{defn}
Given a morphism $X \to Y$ of derived schemes, there is a naturally associated (derived) quasi-coherent $\OO_X$-module $L_{X/Y}$ on $X$ called the \emph{relative cotangent complex}. The relative cotangent complex (of an $A$-scheme) associated to the structure morphism $X \to \Spec(A)$ is called the \emph{absolute cotangent complex}, and is denoted by $L_X$.
\end{defn}
\begin{thm} The cotangent complex satisfies the following basic properties:
\begin{enumerate}
\item given a derived Cartesian square
\begin{equation*}
\begin{tikzcd}
X' \arrow[]{r}{} \arrow[]{d}{f} & Y' \arrow[]{d}{} \\
X \arrow[]{r}{} & Y
\end{tikzcd}
\end{equation*}
we have a natural equivalence $f^* L_{X/Y} \simeq L_{X'/Y'}$;

\item given a triangle $X \stackrel{g}{\to} Y \to Z$ of derived schemes, we have a natural cofibre sequence (a distinguished triangle)
$$
g^* L_{Y/Z} \to L_{X/Z} \to L_{X/Y} 
$$
of quasi-coherent sheaves on $X$;

\item Given a commutative square (in the $\infty$-categorical sense, i.e., equipped with a homotopy realizing the equivalence of the two sides of the square)
\begin{equation*}
\begin{tikzcd}
X' \arrow[]{r}{} \arrow[]{d}{f} & Y' \arrow[]{d}{} \\
X \arrow[]{r}{} & Y,
\end{tikzcd}
\end{equation*}
we have a naturally induced morphism $f^* L_{X/Y} \to L_{X'/Y'}$, which is functorial in the vertical composition of squares;

\item given a derived Cartesian square
\begin{equation*}
\begin{tikzcd}
X' \arrow[]{r}{g} \arrow[]{d}{f} & Y' \arrow[]{d}{} \\
X \arrow[]{r}{} & Y,
\end{tikzcd}
\end{equation*}
the cofibre sequences
$$f^* L_{X/Y} \to L_{X'/Z} \to L_{X'/X}$$
and
$$g^* L_{Y'/Y} \to L_{X'/Z} \to L_{X'/Y'}$$
together with the natural identifications of $(1)$ naturally identify $L_{X'/Y}$ as the direct sum $f^* L_{X/Y} \oplus g^* L_{Y'/Y}$, and the maps as the canonical injections and surjections.
\end{enumerate}
\end{thm}

\begin{rem}
Note that the classical analogue of (1) is \emph{not} true. Indeed, in classical algebraic geometry, the relative cotangent complex is known to be stable under pullbacks only for \emph{Tor-independent squares}. In turn, these are the squares for which the derived fibre product agrees with classical one. This stability of the cotangent complex under derived pullbacks was the crucial technical ingredient that allowed an easy construction of bivariant derived algebraic cobordism in \cite{An}, and the of theory we are going to consider in this article. 
\end{rem}

\begin{defn}
A morphism $X \to Y$ is called \emph{quasi-smooth} if the relative cotangent complex $L_{X/Y}$ is perfect and has (homological) Tor-amplitude in $[1,0]$. A derived scheme $X$ is called \emph{quasi-smooth} if its absolute cotangent complex $L_X$ has amplitude in $[1,0]$. Similarly, a derived scheme (or a morphism of derived schemes) is \emph{smooth} if and only if the cotangent complex is a vector bundle.

A quasi-smooth morphism that is also a closed embedding is called a \emph{quasi-smooth embedding} or a \emph{quasi-smooth immersion}. Often, also the term \emph{derived regular embedding} is used.
\end{defn}

\begin{rem}
If $X$ is quasi-projective, then a perfect complex $\cF$ having Tor-amplitude $[0,1]$ coincides with requiring $\cF$ to be equivalent to the cofibre of a morphism $E_1 \to E_0$ of vector bundles on $X$ (this is a special case of Proposition \fref{ResolutionOfPerfectComplexes}).
\end{rem}

\begin{war}
One should be careful with the definition of (quasi-)smoothness given above. Usually one makes the additional assumption that $f: X \to Y$ should be of finite presentation in homotopical sense. However, it is known (see \cite{Lur2} Theorem 7.4.3.18), that whenever the truncation $tf: tX \to tY$ is of finite presentation in the classical sense, then the cotangent complex of $f$ being an object of $\Perf(X)$ implies the homotopical finite presentation of $f$. As we are mainly interested in derived schemes quasi-projective over a Noetherian ring $A$, this extra finiteness condition on the truncation is always fulfilled, allowing us to use a simpler definition.
\end{war}

\begin{ex}
A morphism $X \to Y$ of classical schemes is quasi-smooth if and only if it is local complete intersection (abbr., l.c.i.).
\end{ex}

\begin{ex}\label{CotangentComplexOfVanishingLoci}
The structure morphisms of vector bundles and projective bundles over $X$ are smooth. The first claim has the following more precise version (as well as a generalization, see \cite{Lur3} Example 25.3.2.2): if $\cF$ is a connective quasi-coherent sheaf over a derived scheme $X$, then the relative cotangent complex of 
$$\pi : \sSpec(\Sym_X^*(\cF)) \to X$$
is naturally equivalent to $\pi^* \cF$.

From the above result, one can conclude that the cotangent complex of any section $X \to E$ of any vector bundle is naturally equivalent to the shifted dual bundle $E^\vee [-1]$. By invariance of the cotangent complex in homotopy pullbacks, we can also conclude that 
$$L_{V(s)/X} \simeq E^\vee [-1] \vert_{V(s)}.$$
\end{ex}

The following theorem follows easily from the basic properties of cotangent complex and from the fact that the cotangent complex of a smooth morphism is just a vector bundle.

\begin{thm}\label{QSProperties}
Quasi-smooth morphisms have the following properties.
\begin{enumerate}
\item Quasi-smooth morphisms are stable under homotopy pullbacks.
\item Quasi-smooth morphisms are stable under composition.
\item Given maps $X \to Y \to Z$ with the latter map smooth, we have that the map $X \to Y$ is quasi-smooth if and only if the composition $X \to Z$ is. 
\end{enumerate}
\end{thm}

\begin{ex}\label{ProjectivizedInclusionIsQs}
As the structure morphisms of projective bundles are smooth, the projectivized inclusion associated to a map of vector bundles as in \fref{ProjectivizedInclusion} is quasi-smooth by the third part of the above theorem.

Moreover, as $\Proj(E) \to \Proj(E \oplus F)$ is by the Proposition \fref{ProjectivizedEmbeddingAsCompleteIntersection} the derived vanishing locus of a section of the vector bundle $\pi_{\Proj(E \oplus F)}^*(F) \otimes \OO_{\Proj(E \oplus F)}(1)$, we can use Example \fref{CotangentComplexOfVanishingLoci} to conclude that
$$\cN^\vee_{\Proj(E) / \Proj(E \oplus F)} \simeq \pi^*_{\Proj(E)} F^\vee \otimes \OO_{\Proj(E)}(-1)$$
\end{ex}

For the convenience of the reader, we also record the following proposition (this is Proposition 2.3.8 of \cite{Khan}).

\begin{prop}\label{ConormalBundle}
A closed embedding $Z \hookrightarrow X$ is quasi-smooth if and only if the shifted cotangent complex $L_{Z/X}[1]$ is a vector bundle. We define $L_{Z/X}[1]$ to be the \emph{conormal bundle} of $Z \hookrightarrow X$, and denote it by $\cN^\vee_{Z/X}$.
\end{prop}

There is also a well defined notion of relative dimension for quasi-smooth maps.

\begin{defn}
The \emph{relative virtual dimension} of a quasi-smooth morphism $X \to Y$ is defined to be $\mathrm{dim}(E_0) - \mathrm{dim}(E_1)$ where $E_0 \to E_1$ is any resolution of the cotangent complex $L_{X/Y}$ (recall that such a resolution always exists if $X$ is quasi-projective). Similarly, using the absolute cotangent complex, one may define the \emph{virtual dimension} of a quasi-smooth derived scheme $X$.
\end{defn}

Of course, the relative virtual dimension can be more generally using either local resolutions or trace methods. It is clear from the above definition that the relative virtual dimension of an equivalence is 0 and that relative virtual dimension is stable under derived pullbacks. 

\begin{ex}
The inclusion of a virtual Cartier divisor $D$ associated to a section of a line bundle $\Li$ on $X$ has relative virtual dimension $-1$. Indeed, analogously to the classical case, conormal bundle of $D \hookrightarrow X$ is equivalent to $\Li^\vee \vert_D$.

Suppose now $X = \Spec(k)$ is the point, $\Li = \OO_X$ is the trivial line bundle, and $s = 0$ is the zero section. The virtual Cartier divisor in this case is a nontrivial quasi-smooth derived scheme of virtual dimension $-1$. Its underlying classical scheme is clearly just $\Spec(k)$. This is one of the technical disadvantages of working with derived schemes: the existence of nontrivial derived schemes of negative dimension often causes complications to proofs that proceed by induction on dimension. 
\end{ex}

\subsection{Derived blow ups}

As derived version of blowing up will be an indispensable tool for us, we recall some of the details from \cite{Khan}. First of all, we recall the definition from the Section 4 of \emph{loc. cit.}:

\begin{defn}
Let $X$ be a derived scheme, and let $Z \hookrightarrow X$ be a quasi-smooth embedding. Now the \emph{derived blow up} $\bl_Z(X)$ of $X$ at $Z$ is the derived $X$-scheme so that the space of $X$-maps $S \to \bl_Z(X)$ is naturally equivalent to the space of virtual Cartier divisor on $S$ lying over $Z$, i.e., the following data:
\begin{enumerate}
\item a commutative diagram
\begin{center}
\begin{tikzcd}
D \arrow[hookrightarrow]{r}{i_D} \arrow[->]{d}{g} & S \arrow[->]{d} \\ 
Z \arrow[hookrightarrow]{r} & X
\end{tikzcd}
\end{center}
so that the top arrow is an inclusion of a virtual Cartier divisor;
\item the above square truncates to a Cartesian (not necessarily homotopy Cartesian) square of classical schemes;
\item the canonical map $g^* \cN^\vee_{Z/X} \to \cN^\vee_{D/S}$ between the conormal bundles (see Proposition \fref{ConormalBundle} for the definition) is a surjection of vector bundles, i.e., $\pi_0 \bigl( g^* \cN^\vee_{Z/X}\bigr) \to \pi_0 \bigl(\cN^\vee_{D/S} \bigr)$ is surjective morphism of sheaves in the classical sense.
\end{enumerate} 
We denote the \emph{structure morphism} $\bl_Z(X) \to \bl_Z(X)$ usually by $\pi$. Note that the identity map  $\bl_Z(X) \to \bl_Z(X)$ gives rise to the \emph{universal virtual Cartier divisor lying over $Z$}, which we are going to call the \emph{exceptional divisor}, and denote by $E \hookrightarrow \bl_Z(X)$.
\end{defn}

The following basic properties come from the Theorem 4.1.5. of \cite{Khan}, and from its proof.

\begin{thm}\label{BasicPropertiesOfBlowUp}
The derived blow up satisfies the following basic properties:
\begin{enumerate}
\item Given a quasi-smooth embedding $Z \hookrightarrow X$, the derived blow up $\bl_Z(X)$ always exists, and is unique up to a contractible choice of isomorphisms.

\item The structure morphism $\pi: \bl_Z(X) \to X$ is proper and quasi-smooth, and induces an equivalence $\bl_Z(X) - E \simeq X - Z$.

\item The canonical surjection $g^* \cN^\vee_{Z/X} \to \cN^\vee_{E/\bl_Z(X)}$ induced by the universal property of the derived blow up, identifies the exceptional divisor as a $Z$-scheme with the projective bundle $\Proj(\cN_{Z/X})$.  (This statement is contained in the discussions of \cite{Khan} in 4.3.4 and 4.3.5.) 

In other words, there is a natural identification $E \simeq \Proj(\cN_{Z/X})$ of derived schemes over $Z$, and a natural identification of surjections $g^* \cN^\vee_{Z/X} \to \cN^\vee_{E/\bl_Z(X)}, g^* \cN^\vee_{Z/X} \to \OO(1)$ of vector bundles on $X$.

\item Given a chain $Z \hookrightarrow Y \hookrightarrow X$ of quasi-smooth immersions, we can form the following homotopy commutative square
\begin{center}
\begin{tikzcd}
E \arrow[hookrightarrow]{r} \arrow[->]{d}{g} & \bl_Z(Y) \arrow[->]{d}{h} \\ 
Z \arrow[hookrightarrow]{r} & X
\end{tikzcd}
\end{center}
where $h$ is the composition $\bl_Z(Y) \to Y \hookrightarrow X$ and $g$ is the projection $\Proj(\cN_{Z/Y}) \to Z$. This square satisfies the conditions necessary for it to define a morphism $\bl_Z(Y) \to \bl_Z(X)$ called the \emph{strict transform} of $Y \hookrightarrow X$. Moreover, the strict transform is a quasi-smooth embedding.

\item If $Z \hookrightarrow X$ is a virtual Cartier divisor, the natural map $\bl_Z(X) \to X$ is an equivalence.

\item If both $X$ and $Z$ are classical schemes, then there is a canonical equivalence $\bl_Z(X) \cong \bl_Z^{\mathrm{cl}}(X)$, where $\bl_Z^{\mathrm{cl}}(X)$ is the classical blow up of $X$ at $Z$.
\end{enumerate}
\end{thm}

As we need to restrict our attention to quasi-projective derived schemes, we need the following stronger version of (2) above.

\begin{prop}
Let $X$ be a derived scheme quasi-projective over a Noetherian ring $A$, and let $Z \hookrightarrow X$ be a regular embedding of virtual codimension $r$. Then the derived blow up $\bl_Z(X)$ is quasi-projective, and therefore the structure morphism $\pi: \bl_Z(X) \to X$ is projective.
\end{prop}
\begin{proof}
Recall from \fref{QProjAndRelativelyAmple} that it is enough to find a line bundle $\Li$ on $\bl_Z(X)$ which is relatively ample over $X$, i.e., it becomes ample when restricted to any preimage of an affine open set of $X$ (or, equivalently, all in any affine open cover of $X$). We claim that we can choose $\Li = \OO(-E)$. Indeed, let $(U_i)_{i \in I}$ be an affine open cover of $X$ with $U_i \simeq \Spec(A_i)$ so that each $U_i$ either does not meet $Z$ at all, or that $Z \vert_{U_i}$ is cut out by $r$ elements $a_1,...a_r \in A_i$. In the former case $\OO(-E) \vert_{\pi^{-1} U_i} \simeq \OO_{U_i}$ is ample as the structure sheaf of an affine scheme. On the other hand, in the latter case, we have the following homotopy Cartesian square
\begin{equation*}
\begin{tikzcd}
\bl_{Z \cap U_i}(U_i) \arrow[]{r} \arrow[]{d} & U_i \arrow[]{d} \\
\bl_{\{0\}}(\A^r) \arrow[]{r} & \A^r.
\end{tikzcd} 
\end{equation*}
Recalling that relative ampleness is stable under pullbacks, and that the pullback of $\OO(1)$ on $\bl_{\{0\}}(\A^r)$ is $\OO(-E) \vert_{\pi^{-1} U_i}$, we see that the latter line bundle is relatively ample over the affine scheme $U_i$ and hence is ample itself.
\end{proof}

We will also need the following lemma later.

\begin{lem}\label{LinearizationInBlowUp}
Let us have a triangle of quasi-smooth immersions $Z \hookrightarrow Y \hookrightarrow X$. Then the derived pullback of the strict transform $\bl_Z(Y) \hookrightarrow \bl_Z(X)$ to the exceptional divisor $E \simeq \Proj(\mathcal{N}_{Z/X})$ is naturally identified as the projectivized inclusion $\Proj(\mathcal{N}_{Z/Y}) \hookrightarrow \Proj(\mathcal{N}_{Z/X})$ associated to the inclusion $\cN_{Z/Y} \hookrightarrow \cN_{Z/X}$ of normal sheaves (which is dual of the surjection of conormal sheaves).
\end{lem}
\begin{proof}
Recall that the strict transform is induced by the outer square in
\begin{center}
\begin{tikzcd}
\Proj(\cN_{Z/Y}) \arrow[hookrightarrow]{r} \arrow[->]{d}{g} & \bl_Z(Y) \arrow[->]{d} \\ 
Z \arrow[hookrightarrow]{r} \arrow[->]{d}{\mathrm{Id}_Z} & Y \arrow[->]{d} \\
Z \arrow[hookrightarrow]{r} & X
\end{tikzcd}
\end{center}
where the upper square is the universal square associated to $\bl_Z(Y)$ and the lower square is induced by the triangle $Z \hookrightarrow Y \hookrightarrow X$. By the basic functoriality properties of the cotangent complex, the morphism $g^* \cN_{Z/X}^\vee \to \cN_{\Proj(\cN_{Z/Y})/\bl_Z(Y)}^\vee$ is naturally identified with the composition
$$g^* \cN_{Z/X}^\vee \to g^* \cN_{Z/Y}^\vee \to \cN_{\Proj(\cN_{Z/Y})/\bl_Z(Y)}^\vee$$
which is exactly the defining property of the projectivized inclusion given in \fref{ProjectivizedInclusion}.
\end{proof}

This allows us to prove the following useful proposition about blowing up at a homotopy intersection.

\begin{prop}\label{BlowUpOfIntersection}
Let $X$ be a derived scheme, and let $Z_1 \hookrightarrow X$ and $Z_2 \hookrightarrow X$ be two derived regular embeddings. Let us denote by $Z_{12} \hookrightarrow X$ the inclusion of the homotopy intersection of $Z_1$ and $Z_2$ inside $X$. Then
\begin{enumerate}
\item the strict transforms $\bl_{Z_{12}}(Z_1) \hookrightarrow \bl_{Z_{12}}(X)$ and $\bl_{Z_{12}}(Z_2) \hookrightarrow \bl_{Z_{12}}(X)$ do not meet inside the blow up $\bl_{Z_{12}}(X)$;

\item the exceptional divisor $\Proj(\cN_{Z_1/X} \vert_{Z_{12}} \oplus \cN_{Z_1/X}\vert_{Z_{12}}) \hookrightarrow \bl_{Z_{12}}(X)$ meets the strict transform $\bl_{Z_{12}}(Z_i) \hookrightarrow \bl_{Z_{12}}(X)$ in $\Proj(\cN_{Z_{3-i}/X} \vert_{Z_{12}})$;

\item the conormal bundle of the above inclusion $\Proj(\cN_{Z_i / X} \vert_{Z_{12}}) \hookrightarrow \Proj(\cN_{Z_1 / X} \vert_{Z_{12}} \oplus \cN_{Z_2 / X} \vert_{Z_{12}})$ is naturally equivalent to $\cN_{Z_{3-i}/X}^\vee \vert_{Z_{12}}$.
\end{enumerate}
\end{prop}
\begin{proof}
Everything follows more or less directly from things we already know. First of all, the natural maps $\cN^\vee_{Z_{12}/X} \to \cN^\vee_{Z_{12}/Z_i}$ are naturally identified as the projections expressing $\cN^\vee_{Z_{12}/X}$ as the direct sum $\cN^\vee_{Z_1/X} \vert_{Z_{12}}  \oplus \cN^\vee_{Z_2/X} \vert_{Z_{12}}$. This, together with Lemma \fref{LinearizationInBlowUp} proves $(2)$. Moreover, as the strict transforms would have to meet inside the exceptional divisor, and as the intersection of the inclusions $\Proj(\cN_{Z_i / X}) \hookrightarrow \Proj(\cN_{Z_1 / X} \oplus \cN_{Z_2 / X})$ is clearly empty, we also obtain $(1)$. The final claim follows from Example \fref{ProjectivizedInclusionIsQs}.
\end{proof}

\section {Fulton--MacPherson's bivariant theory and a universal bivariant theory}\label{FM-BT}

We make a quick review of Fulton--MacPherson's bivariant theory \cite {FM} (also see \cite{Fulton-book}) and a universal bivariant theory \cite{Yo1}. 

Let $\mathcal  V$ be a category which has a final object $pt$ and on which the fiber product or fiber square is well-defined. Also we consider a class of maps, called ``confined maps" (e.g., proper maps, projective maps, in algebraic geometry), which are \emph{closed under composition and base change and contain all the identity maps}, and a class of fiber squares, called ``independent squares" (or ``confined squares``
'', e.g., ``Tor-independent" in algebraic geometry, a fiber square with some extra conditions required on morphisms of the square), which satisfy the following:

(i) if the two inside squares in  
$$\CD
X''@> {h'} >> X' @> {g'} >> X \\
@VV {f''}V @VV {f'}V @VV {f}V\\
Y''@>> {h} > Y' @>> {g} > Y \endCD
\quad \quad \qquad \text{or} \qquad \quad \quad 
\CD
X' @>> {h''} > X \\
@V {f'}VV @VV {f}V\\
Y' @>> {h'} > Y \\
@V {g'}VV @VV {g}V \\
Z'  @>> {h} > Z \endCD
$$
are independent, then the outside square is also independent,

(ii) any square of the following forms are independent:
$$
\xymatrix{X \ar[d]_{f} \ar[r]^{\op {id}_X}&  X \ar[d]^f & & X \ar[d]_{\op {id}_X} \ar[r]^f & Y \ar[d]^{\op {id}_Y} \\
Y \ar[r]_{\op {id}_X}  & Y && X \ar[r]_f & Y}
$$
where $f:X \to Y$ is \emph{any} morphism. 
\begin{defn}\label{BivariantTheory}
A \emph{bivariant theory} $\mathbb B$ on a category $\mathcal  V$ with values in the category of graded abelian groups is an assignment to each morphism
$ X  \xrightarrow{f} Y$
in the category $\Cal V$ a graded abelian group (in most cases we ignore the grading )
$\bB(X  \xrightarrow{f} Y)$
which is equipped with the following three basic operations. The $i$-th component of $\bB(X  \xrightarrow{f} Y)$, $i \in \bZ$, is denoted by $\bB^i(X  \xrightarrow{f} Y)$.
\begin{enumerate}
\item {\bf Product}: For morphisms $f: X \to Y$ and $g: Y
\to Z$, the product operation
$$\bullet: \bB^i( X  \xrightarrow{f}  Y) \otimes \bB^j( Y  \xrightarrow{g}  Z) \to
\bB^{i+j}( X  \xrightarrow{gf}  Z)$$
is  defined.

\item {\bf Pushforward}: For morphisms $f: X \to Y$
and $g: Y \to Z$ with $f$ \emph {confined}, the pushforward operation
$$f_*: \bB^i( X  \xrightarrow{gf} Z) \to \bB^i( Y  \xrightarrow{g}  Z) $$
is  defined.

\item {\bf Pullback} : For an \emph{independent} square \qquad $\CD
X' @> g' >> X \\
@V f' VV @VV f V\\
Y' @>> g > Y, \endCD
$

the pullback operation
$$g^* : \bB^i( X  \xrightarrow{f} Y) \to \bB^i( X'  \xrightarrow{f'} Y') $$
is  defined.
\end{enumerate}
These three operations are required to satisfy the following seven compatibility axioms (\cite [Part I, \S 2.2]{FM}):

\begin{enumerate}
\item[($A_1$)] {\bf Product is associative}: given a diagram $X \xrightarrow f Y  \xrightarrow g Z \xrightarrow h  W$ with $\alp \in \bB(X \xrightarrow f Y),  \be \in \bB(Y \xrightarrow g Z), \ga \in \bB(Z \xrightarrow h W)$,
$$(\alp \bullet\be) \bullet \ga = \alp \bullet (\be \bullet \ga).$$
\item[($A_2$)] {\bf Pushforward is functorial} : given a diagram $X \xrightarrow f Y  \xrightarrow g Z \xrightarrow h  W$ with $f$ and $g$ confined and $\alp \in \bB(X \xrightarrow {h\circ g\circ f} W)$
$$(g\circ f)_* (\alp) = g_*(f_*(\alp)).$$
\item[($A_3$)] {\bf Pullback is functorial}: given independent squares
$$\CD
X''@> {h'} >> X' @> {g'} >> X \\
@VV {f''}V @VV {f'}V @VV {f}V\\
Y''@>> {h} > Y' @>> {g} > Y \endCD
$$
and $\alp \in \bB(X \xrightarrow f Y)$,
$$(g \circ h)^*(\alp) = h^*(g^*(\alp)).$$
\item[($A_{12}$)] {\bf Product and pushforward commute}: given a diagram $X \xrightarrow f Y  \xrightarrow g Z \xrightarrow h  W$ with $f$ confined and $\alp \in \bB(X \xrightarrow {g \circ f} Z),  \be \in \bB(Z \xrightarrow h W)$,
$$f_*(\alp \bullet\be)  = f_*(\alp) \bullet \be.$$
\item[($A_{13}$)] {\bf Product and pullback commute}: given independent squares
$$\CD
X' @> {h''} >> X \\
@V {f'}VV @VV {f}V\\
Y' @> {h'} >> Y \\
@V {g'}VV @VV {g}V \\
Z'  @>> {h} > Z \endCD
$$
with $\alp \in \bB(X \xrightarrow {f} Y),  \be \in \bB(Y \xrightarrow g Z)$,
$$h^*(\alp \bullet\be)  = {h'}^*(\alp) \bullet h^*(\be).$$
\item[($A_{23}$)] \label{push-pull}{\bf Pushforward and pullback commute}: given independent squares
$$\CD
X' @> {h''} >> X \\
@V {f'}VV @VV {f}V\\
Y' @> {h'} >> Y \\
@V {g'}VV @VV {g}V \\
Z'  @>> {h} > Z \endCD
$$
with $f$ confined and $\alp \in \bB(X \xrightarrow {g\circ f} Z)$,
$$f'_*(h^*(\alp))  = h^*(f_*(\alp)).$$
\item[($A_{123}$)] {\bf Projection formula}: given an independent square with $g$ confined and $\alp \in \bB(X \xrightarrow {f} Y),  \be \in \bB(Y' \xrightarrow {h \circ g} Z)$
$$\CD
X' @> {g'} >> X \\
@V {f'}VV @VV {f}V\\
Y' @>> {g} > Y @>> h >Z \\
\endCD
$$
and $\alp \in \bB(X \xrightarrow {f} Y),  \be \in \bB(Y' \xrightarrow {h \circ g} Z)$,
$$g'_*(g^*(\alp) \bullet \be)  = \alp \bullet g_*(\be).$$
\end{enumerate}
Finally, we also require the theory $\bB$ to have multiplicative units:
\begin{enumerate}
\item[(U)] For all $X \in \mathcal{V}$, there is an element $1_X \in \bB^0( X  \xrightarrow{\op {id}_X} X)$ such that $\alp \bullet 1_X = \alp$ for all morphisms $W \to X$ and all $\alp \in \bB(W \to X)$, and such that $1_X \bullet \beta = \beta $ for all morphisms $X \to Y$ and all $\beta \in \bB(X \to Y)$, and such that $g^*1_X = 1_{X'}$ for all $g: X' \to X$.
\end{enumerate}

\end{defn}

The theories we are going to encounter in this paper satisfy the following extra condition.

\begin{defn}\label{Commutativity}
A bivariant theory $\bB$ is called \emph{commutative} if whenever both
$$\CD
W @> {g'} >> X \\
@V {f'}VV @VV {f}V\\
Y @>> {g} > Z  \\
\endCD  
\quad \quad \text{and} \quad \quad 
\CD
W @> {f'} >> Y \\
@V {g'}VV @VV {g}V\\
X @>> {g} > Z \\
\endCD  
$$
are independent squares with $\alp \in \bB(X \xrightarrow f Z)$ and $\be \in \bB(Y \xrightarrow g Z)$,
$$g^*(\alp) \bullet \be = f^*(\be) \bullet \alp .$$
\end{defn}

We will also recall the correct notion of a morphism between two bivariant theories:

\begin{defn}\label{groth}
Let $\bB, \bB'$ be two bivariant theories on a category $\mathcal V$. A {\it Grothendieck transformation} from $\bB$ to $\bB'$, $\ga : \bB \to \bB'$
is a collection of homomorphisms
$\bB(X \to Y) \to \bB'(X \to Y)$
for a morphism $X \to Y$ in the category $\mathcal V$, which preserves the above three basic operations: 
\begin{enumerate}
\item $\ga (\alp \bullet_{\bB} \be) = \ga (\alp) \bullet _{\bB'} \ga (\be)$, 
\item $\ga(f_{*}\alp) = f_*\ga (\alp)$, and 
\item $\ga (g^* \alp) = g^* \ga (\alp)$. 
\end{enumerate}
\end{defn}

\begin{defn}\label{BI}
Let $\bB$ be a bivariant theory. A \emph{bivariant ideal} $\bI \subset \bB$ consists of (graded) subgroups $\bI(X \xrightarrow {f} Y) \subset \bB(X \xrightarrow {f} Y)$ for each $f:X \to Y$ so that
\begin{enumerate}
\item if $\alpha \in \bI(X \xrightarrow {g\circ f} Z)$, then $f_* \alpha \in \bI(Y \xrightarrow {g}  Z)$ for $f: X \to Y$ confined;

\item if $\alpha \in \bI(X \xrightarrow {f}  Y)$, then $g^* \alpha \in \bI(X' \xrightarrow {f'}  Y')$ for all $g: Y' \to Y$ so that the Cartesian square
$$\CD
X' @> {g'} >> X \\
@V {f'} VV @VV {f} V\\
Y' @>> {g} > Y \endCD
$$
is independent;

\item if $\alpha \in \bI(X \to Y)$, then $\beta \bullet \alpha \in \bI(X' \to Y)$ for any $\beta \in \bB(X' \to X)$ and $\alpha \bullet \gamma \in \bB(X \to Y')$ for any $\gamma \in \bB(Y \to Y')$.
\end{enumerate}
\end{defn}

\noindent Bivariant ideals are clearly to bivariant theories what ideals are to rings. Namely: 

\begin{prop}
\noindent
\begin{enumerate}
\item The (object-wise) kernel of a Grothendieck transformation $\gamma: \bB \to \bB'$ is a bivariant ideal.

\item Given a bivariant ideal $\bI \subset \bB$, one may form the \emph{quotient bivariant theory} $\bB / \bI$ by setting $(\bB / \bI)(X \to Y) := \bB(X \to Y) / \bI(X \to Y)$ and by taking the bivariant operations to be the ones induced by $\bB$. Namely, they are defined as follows: 
\begin{enumerate}
\item {\bf Product}: For morphisms $f: X \to Y$ and $g: Y
\to Z$, the product operation
$$\bullet: (\bB / \bI)^i( X  \xrightarrow{f}  Y) \otimes (\bB / \bI)^j( Y  \xrightarrow{g}  Z) \to
(\bB / \bI)^{i+j}( X  \xrightarrow{gf}  Z)$$
is  defined by $[\alp] \bullet [\be]:=[\alp \bullet \be].$

\item {\bf Pushforward}: For morphisms $f: X \to Y$
and $g: Y \to Z$ with $f$ \emph {confined}, the pushforward operation
$$f_*: (\bB / \bI)^i( X  \xrightarrow{gf} Z) \to (\bB / \bI)^i( Y  \xrightarrow{g}  Z) $$
is  defined by $f_*([\alp]):=[f_*\alp]$.

\item {\bf Pullback} : For an \emph{independent} square 
$$\CD
X' @> g' >> X \\
@V f' VV @VV f V\\
Y' @>> g > Y, \endCD
$$
\newline the pullback operation
$$g^* : (\bB / \bI)^i( X  \xrightarrow{f} Y) \to (\bB / \bI)^i( X'  \xrightarrow{f'} Y') $$
is  defined by $g^*([\alp]):=[g^*\alp]$.
\end{enumerate}
\end{enumerate}
\end{prop}

\begin{proof}
The proofs are easy, but 
we give a proof.
\begin{enumerate}
\item The conditions (1), (2) and (3) of Definition \ref{BI} of a bivariant ideal are satisfied respectively by the requirements (2), (3) 
and (1) of Definition \ref{groth} of a Grothendieck transformation. \newline
\item It suffices to show that the above associated three bivariant operations are well-defined, i.e., do not depend on representatives.

\begin{enumerate}
\item Suppose that $[\alp]=[\alp'] \in (\bB / \bI)^i( X  \xrightarrow{f}  Y)$ and $[\be]=[\be'] \in (\bB / \bI)^j( Y  \xrightarrow{g}  Z) $, i.e., $\alp = \alp' + a$ with $a\in \bI(X  \xrightarrow{f} Y)$
and $\be = \be' + b$ with $b\in \bI(Y  \xrightarrow{g} Z)$. Then 
$$\alp \bullet \be = (\alp' + a) \bullet (\be' + b) =\alp' \bullet \be' + \alp' \bullet b + a \bullet \be' + a \bullet b$$
It follows from the condition $(3)$ of Definition \ref{BI} that the last three terms $\alp' \bullet b + a \bullet \be' + a \bullet b$ belong to $\bI(X \xrightarrow{g \circ f } Z)$. Hence $[\alp \bullet \be] = [\alp' \bullet \be'].$
\item Suppose that $[\alp]=[\alp'] \in (\bB / \bI)^i( X  \xrightarrow{gf} Z) $, i.e., $\alp = \alp' + a$ with $a\in \bI(X  \xrightarrow{gf}  Z)$.
Then $f_*\alp = f_*\alp'  + f_*a$ and it follows from the condition $(1)$ of Definition \ref{BI} that $f_*a \in \bI ( Y \xrightarrow{g} Z) $. Therefore we have $[f_*\alp] = [f_*\alp']$.
\item Suppose that $[\alp]=[\alp'] \in (\bB / \bI)^i( X  \xrightarrow{f} Y) $, i.e., $\alp = \alp' + a$ with $a\in \bI(X  \xrightarrow{f}  Y)$.
Then $g^*\alp = g^*\alp'  + g^*a$ and it follows from the condition $(2)$ of Definition \ref{BI} that $g^*a \in \bI ( X' \xrightarrow{f'} Y') $. Therefore we have $[g^*\alp] = [g^*\alp']$
\end{enumerate}

Then the seven axioms automatically follow because the seven axioms hold for representatives of equivalence classes. For example, since $([\alp]\bullet [\be]) \bullet [\gamma] = [(\alp \bullet \be) \bullet \gamma]$ and $[\alp]\bullet ([\be] \bullet [\gamma]) = [\alp \bullet (\be \bullet \gamma)]$ by the definition of the bivariant product $\bullet$ for $\mathbb B/\mathbb I$, the associativity 
$$([\alp]\bullet [\be]) \bullet [\gamma] = [\alp]\bullet ([\be] \bullet [\gamma])$$
follows from $(\alp \bullet \be) \bullet \gamma = \alp \bullet (\be \bullet \gamma)$
on the level of representatives.
\end{enumerate}
\end{proof}

\begin{rem}
\begin{enumerate}
\item
The definitions of the above three bivariant operations for $\mathbb B/ \mathbb I$ given in (2) of Proposition \ref{BI} should be denoted differently to avoid some possible confusion with those on the original one $\mathbb B$, e.g., the product $\bullet_{\mathbb I}$, the pushforward $[f_*]$ and the pullback $[g^*]$, but we use the same symbols.

\item These defintions, i.e., $[\alp]\bullet [\beta] =[\alp \bullet \beta] $, $f_*([\alp]) =[f_*(\alp)]$ and $g^*([\alp]) =[g^*(\alp)]$, also mean in other words that
the quotient map $\Theta: \mathbb B \to \mathbb B/\mathbb I$ defined by $\Theta (\alp):=[\alp]$ is a Grothendieck transformation, i.e.,
$\Theta (\alp \bullet \beta) =\Theta (\alp) \bullet \Theta (\beta) $, $\Theta(f_*(\alp))= f_*\Theta (\alp)$ and $\Theta(g^*(\alp)) =g^*\Theta(\alp)$
\end{enumerate}
\end{rem}

As in the case of rings, we do have a simple description of the bivariant ideal generated by a subset, at least in good situations. By a \emph{(bivariant) subset} $S$ of $\bB$ (denoted $S \subset \bB$), we mean a collection of subsets $S(X \to Y) \subset \bB(X \to Y)$ --- one for each map $X \to Y$ in $\mathcal V$. Given a subset $S \subset \bB$, we denote by $\langle S \rangle$ the \emph{bivariant ideal generated by $S$}, i.e., the smallest bivariant ideal of $\bB$ containing $S$. When we need to make it clear in which bivariant theory $\mathbb B$ you consider such a bivariant ideal $\langle S \rangle$, we denote it by $\langle S \rangle_{\mathbb B}$.

\begin{prop}\label{GeneratedBivariantIdeal}
Let $S$ be a bivariant subset of $\bB$. Moreover, assume that \emph{all Cartesian squares are independent}. Now $\langle S \rangle(X \stackrel h \to Y)$ consists of elements of the form
\begin{equation*}
f_*(\alpha \bullet g^*(s) \bullet \beta)
\end{equation*}
where $f, g, \alpha, \beta$ and $s$ are as in the following diagram and $s \in S(A \to B)$.
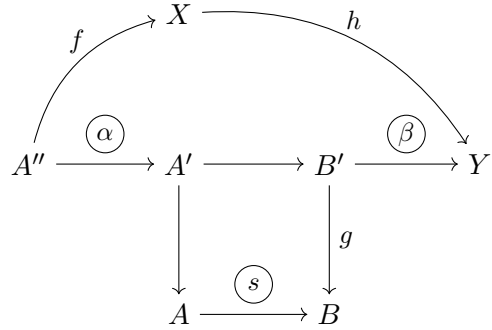
\begin{center}
\begin{tikzpicture}[scale=2]
\node (A1) at (1,0) {$A$};
\node (A2) at (2,0) {$B$};
\node (B0) at (0,1) {$A''$};
\node (B1) at (1,1) {$A'$};
\node (B2) at (2,1) {$B'$};
\node (B3) at (3,1) {$Y$};
\node (X) at (1,2) {$X$};
\path[every node/.style={font=\sffamily\small}]
(A1) edge[->] node[yshift=0.4cm,circle,draw,inner sep=0pt, minimum size=0.5cm]{$s$} (A2)
(B0) edge[->] node[yshift=0.4cm,circle,draw,inner sep=0pt, minimum size=0.5cm]{$\alpha$} (B1)
(B1) edge[->] (B2)
(B2) edge[->] node[yshift=0.4cm,circle,draw,inner sep=0pt, minimum size=0.5cm]{$\beta$} (B3)
(B1) edge[->] (A1)
(B2) edge[->] node[right]{$g$} (A2)
(B0) edge[bend left,->] node[above]{$f$} (X)
(X) edge[bend left,->] node[above]{$h$} (B3)
;
\end{tikzpicture}
\end{center}
In the above diagram, the bottom square is assumed to be Cartesian (hence independent), and $f$ to be confined.
\end{prop}
\begin{rem}
Following \cite{FM}, in the above diagram 
\begin{center}
\begin{tikzpicture}[scale=2]
\node (B0) at (0,1) {$A''$};
\node (B1) at (1,1) {$A'$};
\path[every node/.style={font=\sffamily\small}]
(B0) edge[->] node[yshift=0.4cm,circle,draw,inner sep=0pt, minimum size=0.5cm]{$\alpha$} (B1)
;
\end{tikzpicture}
\end{center}
means $\alpha \in \B^*(A'' \to A')$ is a bivariant element.
\end{rem}

\begin{proof}
Clearly elements of this form all lie in $\langle S \rangle$, so we only need to show that the description above gives a bivariant ideal.
\begin{enumerate}
\item Suppose $f'$ is confined. It is now enough to show that $f'_* (f_*(\alpha \bullet g^*(s) \bullet \beta))$ can be expressed in the above form. But this is trivial, as by functoriality of bivariant pushforward the above element is just $(f' \circ f)_*(\alpha \bullet g^*(s) \bullet \beta)$.

\item Suppose $i: Y' \to Y$ is a map and consider the the Cartesian diagram
$$\CD
C'' @>> {f'} > X' @> {h'} >> Y' \\
@V {i_3} VV @V {i_2} VV @VV i V\\
A'' @>> {f} > X @>> {h} > Y \endCD.
$$
As all the squares are independent, we can use the bivariant axiom $(A_{23})$ to conclude that 
$$i^*(f_*(\alpha \bullet g^*(s) \bullet \beta)) = f'_*(i^*(\alpha \bullet g^*(s) \bullet \beta)).$$ We can also consider the Cartesian diagram
$$\CD
C'' @>> {} > C' @> {} >> D' @>> {} > Y' \\
@V {i_3} VV @V {i_5} VV @V {i_4} VV @VV i V\\
A'' @>> {} > A' @>> {} > B' @>> {} > Y \endCD
$$
and again, as every square is independent, we can conclude using the bivariant axiom $(A_{13})$ that
$$i^*(\alpha \bullet g^*(s) \bullet \beta) = i_5^*(\alpha) \bullet i_4^*(g^*(s)) \bullet i^*(\beta).$$

Therefore we have
$$i^*(f_*(\alpha \bullet g^*(s) \bullet \beta)) = f'_* \Bigl (i_5^*(\alpha) \bullet i_4^*(g^*(s)) \bullet i^*(\beta) \Bigr),$$
which  is of the form $f'_*(\alpha' \bullet (g \circ i_4)^*(s) \bullet \beta')$ where the situation is described in the following diagram
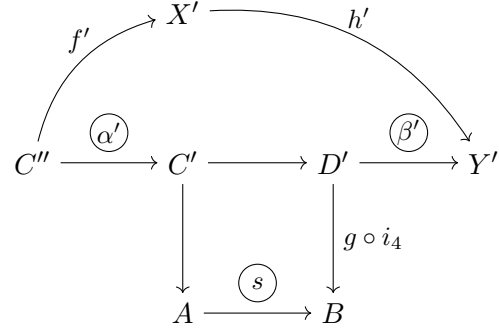
\begin{center}
\begin{tikzpicture}[scale=2]
\node (A1) at (1,0) {$A$};
\node (A2) at (2,0) {$B$};
\node (B0) at (0,1) {$C''$};
\node (B1) at (1,1) {$C'$};
\node (B2) at (2,1) {$D'$};
\node (B3) at (3,1) {$Y'$};
\node (X) at (1,2) {$X'$};
\path[every node/.style={font=\sffamily\small}]
(A1) edge[->] node[yshift=0.4cm,circle,draw,inner sep=0pt, minimum size=0.5cm]{$s$} (A2)
(B0) edge[->] node[yshift=0.4cm,circle,draw,inner sep=0pt, minimum size=0.5cm]{$\alpha'$} (B1)
(B1) edge[->] (B2)
(B2) edge[->] node[yshift=0.4cm,circle,draw,inner sep=0pt, minimum size=0.5cm]{$\beta'$} (B3)
(B1) edge[->] (A1)
(B2) edge[->] node[right]{$g \circ i_4$} (A2)
(B0) edge[bend left,->] node[above]{$f'$} (X)
(X) edge[bend left,->] node[above]{$h'$} (B3)
;
\end{tikzpicture}
\end{center}
This concludes the proof.

\item First of all, given an element $\gamma \in \bB(Y \to Y')$, we conclude immediately using $(A_{12})$ and  $(A_1)$ that 
\begin{align*}
f_* \Bigl (\alpha \bullet g^*(s) \bullet \beta \Bigr ) \bullet \gamma & = f_* \Bigl (\bigl (\alpha \bullet g^*(s) \bullet \beta \bigr ) \bullet \gamma  \Bigr)  \quad \text{(by $(A_{12})$)} \\
& = f_* \Bigl (\alpha \bullet g^*(s) \bullet (\beta \bullet \gamma) \Bigr ) \quad \text{(by $(A_1)$).} \\
\end{align*}
On the other hand, given $\gamma \in \bB(X' \xrightarrow{g} X)$ we may form the Cartesian diagram
$$\CD
E'' @>> {f'} > X'  \\
@V {g'} VV @V {g} VV \\
A'' @>> {f} > X @>> {h} > Y \endCD
$$
and use the bivariant projection formula $(A_{123})$ and  $(A_1)$ to conclude that 
\begin{align*}
\gamma \bullet f_*(\alpha \bullet g^*(s) \bullet \beta)  & = f'_* \Bigl (f^*(\gamma) \bullet \bigl (\alpha \bullet g^*(s) \bullet \beta \bigr )\Bigr )  \quad \text{(by $(A_{123})$)} \\
& = f'_* \Bigl ( \bigl (f^*(\gamma) \bullet \alpha \bigr )  \bullet g^*(s) \bullet \beta  \Bigr )\quad \text{(by $(A_1)$).} \\
&  =: f'_* \Bigl ( \alpha'  \bullet g^*(s) \bullet \beta  \Bigr )
\end{align*}
where, again, the situation can be described by the diagram
\begin{center}
\begin{tikzpicture}[scale=2]
\node (A1) at (1,0) {$A$};
\node (A2) at (2,0) {$B$};
\node (B0) at (0,1) {$E''$};
\node (B1) at (1,1) {$A'$};
\node (B2) at (2,1) {$B'$};
\node (B3) at (3,1) {$Y$};
\node (X) at (1,2) {$X'$};
\path[every node/.style={font=\sffamily\small}]
(A1) edge[->] node[yshift=0.4cm,circle,draw,inner sep=0pt, minimum size=0.5cm]{$s$} (A2)
(B0) edge[->] node[yshift=0.4cm,circle,draw,inner sep=0pt, minimum size=0.5cm]{$\alpha'$} node[below]{$f \circ g'$} (B1)
(B1) edge[->] (B2)
(B2) edge[->] node[yshift=0.4cm,circle,draw,inner sep=0pt, minimum size=0.5cm]{$\beta$} (B3)
(B1) edge[->] (A1)
(B2) edge[->] node[right]{$g$} (A2)
(B0) edge[bend left,->] node[above]{$f'$} (X)
(X) edge[bend left,->] node[above]{$h \circ g$} (B3)
;
\end{tikzpicture}
\end{center}
Thus we are done.
\end{enumerate}
\end{proof}

A bivariant theory unifies both a covariant theory and a contravariant theory in the following sense:
\begin{itemize}
\item We have the associated \emph{homology groups} $\bB_*(X):= \bB(X \to pt)$ which are covariant for confined morphisms, and where the grading is given by $\bB_i(X):= \bB^{-i}(X  \xrightarrow{id}  pt)$.
\item We have the associated \emph{cohomology groups} $\bB^*(X) := \bB(X  \xrightarrow{id}  X)$ which are contravariant for all morphisms, and whose grading is given by $\bB^j(X):= \bB^j(X  \xrightarrow{id}  X)$.
\end{itemize}
\noindent
A Grothendieck transformation $\ga: \bB \to \bB'$ induces natural transformations $\ga_*: \bB_* \to \bB_*'$ and $\ga^*: \bB^* \to {\bB'}^*$.
\begin{rem}
The cohomology groups $\bB^*(X)$ are closed under the bivariant product, and actually form a ring (the associated \emph{cohomology rings}). The axiom $(U)$ in Definition \fref{BivariantTheory} of bivariant theory makes these rings unital, and if the extra assumption of commutativity (as in Definition \fref{Commutativity}) holds, then the rings are commutative. Moreover, it follows from the bivariant axiom $(A_{13})$ that the contravariant pullback maps of the associated cohomology theory $\bB^*$ respect multiplication, and hence $\bB^*$ is a multiplicative cohomology theory.
\end{rem}
 
\begin{defn}\label{canonical}(\cite[Part I, \S 2.6.2 Definition]{FM}) Let $\Cal S$ be a class of maps in $\Cal V$, which is closed under compositions and containing all identity maps. Suppose that to each $f: X \to Y$ in $\Cal S$ there is assigned an element
$\theta(f) \in \bB(X  \xrightarrow {f} Y)$ satisfying that
\begin{enumerate}
\item [(i)] $\theta (g \circ f) = \theta(f) \bullet \theta(g)$ for all $f:X \to Y$, $g: Y \to Z \in \Cal S$ and
\item [(ii)] $\theta(\op {id}_X) = 1_X $ for all $X$ with $1_X \in \bB^*(X):= B(X  \xrightarrow{\op {id}_X} X)$ the unit element.
\end{enumerate}
Then $\theta(f)$ is called an {\it orientation} of $f$. (In \cite[Part I, \S 2.6.2 Definition]{FM} it is called a {\it canonical orientation} of $f$, but in this paper it shall be simply called an orientation.)
\end{defn} 

\begin{defn}\label{strong-orient}
Let $\bB$ be a bivariant theory with an orientation $\theta$. We say that the orientation $\theta(g)$ of $g:Y \to Z$ is \emph{strong} if the maps
\begin{equation*}
- \bullet \theta(f) :\bB^*(X  \xrightarrow f Y) \to \bB^*(X \xrightarrow {g \circ f} Z)
\end{equation*}
are isomorphisms for all $f:X \to Y$.
\end{defn}

\noindent
Gysin homomorphisms \label{gysin}:
Note that such an orientation makes the covariant functor $\bB_*(X)$ a contravariant functor for morphisms in $\Cal S$, and also makes the contravariant functor $\bB^*$ a covariant functor for morphisms in $\Cal C \cap \Cal S$: Indeed, 
\begin{enumerate}
\item As to the covariant functor $\bB_*(X)$: For a morphism $f: X \to Y \in \Cal S$ and the orientation $\theta$ on $\Cal S$ the following {\it Gysin (pullback) homomorphism}
$$f^!: \bB_*(Y) \to \bB_*(X) \quad \text {defined by} \quad  f^!(\alp) :=\theta(f) \bullet \alp$$
 is {\it contravariantly functorial}. 
\item As to contravariant functor $\bB^*$: For a fiber square (which is an independent square by hypothesis)
$$\CD
X @> f >> Y \\
@V {\op {id}_X} VV @VV {\op {id}_Y}V\\
X @>> f > Y, \endCD
$$
where $f \in \Cal C \cap  \Cal S$, the following {\it Gysin (pushforward) homomorphism}
$$f_!: \bB^*(X) \to \bB^*(Y) \quad \text {defined by} \quad
f_!(\alp) := f_*(\alp \bullet \theta (f))$$
is {\it covariantly functorial}.
\end{enumerate}

The above notation $f^!$ and $f_!$ should carry the information of $\Cal S$ and the orientation $\theta$, but it will be usually omitted if it is not necessary to be mentioned. Note that the above conditions (i) and (ii) of Definition (\ref{canonical}) are certainly \emph{necessary} for the above Gysin homomorphisms to be functorial. 

\begin{defn} (i) Let $\Cal S$ be another class of maps called ``specialized maps" (e.g., smooth maps in algebraic geometry) in $\Cal V$ , which is closed under composition, closed under base change and containing all identity maps. Let $\bB$ be a 
bivariant theory. If $\Cal S$ has  orientations in $\bB$, then we say that $\Cal S$ is $\bB$-oriented and an element of $\Cal S$ is called a $\bB$-oriented morphism. (Of course $\Cal S$ is also a class of confined maps, but since we consider the above extra condition of $\bB$-orientation  on $\Cal S$, we give a different name to $\Cal S$.)

(ii) Let $\Cal S$ be as in (i). Let $\bB$ be a bivariant theory and $\Cal S$ be $\bB$-oriented. Furthermore, if the orientation $\theta$ on $\Cal S$ satisfies that for an independent square with $f \in \Cal S$
$$
\CD
X' @> g' >> X\\
@Vf'VV   @VV f V \\
Y' @>> g > Y
\endCD
$$
the following condition holds: 
$\theta (f') = g^* \theta (f)$, 
(which means that the orientation $\theta$ preserves the pullback operation), then we call $\theta$ a {\it nice canonical orientation} and say that $\Cal S$ is {\it nice canonically $\bB$-oriented} and an element of $\Cal S$ is called {\it a nice canonically  $\bB$-oriented morphism} . (Note that in \cite{SY1} they are respectively called {\it a stable orientation}, {\it stably $\bB$-oriented} and {\it a stably $\bB$-oriented morphism}.)

\end{defn}

The following theorem is about {\it the existence of a universal one} of the bivariant theories for a given category $\Cal V$ with a class $\Cal C$ of confined morphisms, a class of independent squares and a class $\Cal S$ of specialized morphisms. 

\begin{thm}(\cite[Theorem 3.1]{Yo1})(A universal bivariant theory) \label{UBT} Let  $\Cal V$ be a category with a class $\Cal C$ of confined morphisms, a class of independent squares and a class $\Cal S$ of specialized maps.  We define 
$$\bM^{\Cal C} _{\Cal S}(X  \xrightarrow{f}  Y)$$
to be the free abelian group generated by the set of isomorphism classes of confined morphisms $h: W \to X$  such that the composite of  $h$ and $f$ is a specialized map:
$$h \in \Cal C \quad \text {and} \quad f \circ h: W \to Y \in \Cal S.$$
\begin{enumerate}
\item  The association $\bM^{\Cal C} _{\Cal S}$ is a bivariant theory if the three bivariant operations are defined as follows:
\begin{enumerate}
\item  {\bf Product}: For morphisms $f: X \to Y$ and $g: Y
\to Z$, the product operation
$$\bullet: \bM^{\Cal C} _{\Cal S} ( X  \xrightarrow{f}  Y) \otimes \bM^{\Cal C} _{\Cal S} ( Y  \xrightarrow{g}  Z) \to
\bM^{\Cal C} _{\Cal S} ( X  \xrightarrow{gf}  Z)$$
is  defined by
$$[V \xrightarrow{h}  X] \bullet [W  \xrightarrow{ k}  Y]:= [V'  \xrightarrow{ h \circ {k}''}  X]$$
and extended linearly, where we consider the following fiber squares
$$\CD
V' @> {h'} >> X' @> {f'} >> W \\
@V {{k}''}VV @V {{k}'}VV @V {k}VV\\
V@>> {h} > X @>> {f} > Y @>> {g} > Z .\endCD
$$
\item {\bf Pushforward}: For morphisms $f: X \to Y$
and $g: Y \to Z$ with $f$ 
confined, the pushforward operation
$$f_*: \bM^{\Cal C} _{\Cal S} ( X  \xrightarrow{gf} Z) \to \bM^{\Cal C} _{\Cal S} ( Y  \xrightarrow{g}  Z) $$
is  defined by
$$f_*\left ([V \xrightarrow{h}  X] \right) := [V  \xrightarrow{f \circ h}  Y]$$
and extended linearly.

\item {\bf Pullback}: For an independent square
$$\CD
X' @> g' >> X \\
@V f' VV @VV f V\\
Y' @>> g > Y, \endCD
$$
the pullback operation
$$g^* : \bM^{\Cal C} _{\Cal S} ( X  \xrightarrow{f} Y) \to \bM^{\Cal C} _{\Cal S}( X'  \xrightarrow{f'} Y') $$
is  defined by
$$g^*\left ([V  \xrightarrow{h}  X] \right):=  [V'  \xrightarrow{{h}'}  X']$$
and extended linearly, where we consider the following fiber squares:
$$\CD
V' @> g'' >> V \\
@V {{h}'} VV @VV {h}V\\
X' @> g' >> X \\
@V f' VV @VV f V\\
Y' @>> g > Y. \endCD
$$
\end{enumerate}
\item Let $\Cal {BT}$ be a class of bivariant theories $\bB$ on the same category $\Cal V$ with a class $\Cal C$ of confined morphisms, a class of independent squares and a class $\Cal S$ of specialized maps. 
Let $\Cal S$ be nice canonically $\bB$-oriented for any bivariant theory $\bB \in \Cal {BT}$. Then, for each bivariant theory $\bB \in \Cal {BT}$ there exists a unique Grothendieck transformation
$$\ga_{\bB} : \bM^{\Cal C} _{\Cal S} \to \bB$$
such that for a specialized morphism $f: X \to Y \in \Cal S$ the homomorphism
$\ga_{\bB} : \bM^{\Cal C} _{\Cal S}(X  \xrightarrow{f}  Y) \to \bB(X  \xrightarrow{f}  Y)$
satisfies the normalization condition that $$\ga_{\bB}([X  \xrightarrow{\op {id}_X}  X]) = \theta_{\bB}(f).$$
\end{enumerate}
\end{thm}

\section{Oriented bivariant theory and a universal oriented bivariant theory }

Levine--Morel's algebraic cobordism is the universal one among the so-called  \emph {oriented} Borel--Moore functors with products for algebraic schemes. Here ``oriented" means that the given Borel--Moore functor $H_*$ is equipped with the endomorphsim $\tilde c_1(L): H_*(X) \to H_*(X)$ for a line bundle $L$ over the scheme $X$. Motivated by this ``orientation" (which is different from the one given in Definition \ref{canonical}, but we still call this ``orientation" using a different symbol so that the reader will not be confused with terminologies), in \cite[\S4]{Yo1} we introduce an orientation to bivariant theories for any category, using the notion of \emph {fibered categories} in abstract category theory (e.g, see \cite{Vistoli}) and such a bivariant theory equipped with such an orientation (Chern class operator) is called \emph{an oriented bivariant theory}. 

\begin{defn}
Let $\Cal L$ be a fibered category over $\Cal V$. An object in the fiber $\Cal L(X)$ over an object $X \in \Cal V$ is called an \emph {``fiber-object over $X$"}, abusing words, and denoted by $L$, $M$, etc.
\end{defn}



\begin{defn}(\cite[Definition 4.2]{Yo1}) (an oriented bivariant theory) \label{orientation} Let $\bB$ be a bivariant theory on a category $\Cal V$.

\begin{enumerate}
\item For a fiber-object $L$ over $X$, the \emph{``operator" on $\bB$ associated to $L$}, denoted by $\phi(L)$, is defined to be an endomorphism
$$\phi(L): \bB(X  \xrightarrow{f}  Y) \to \bB(X  \xrightarrow{f}  Y) $$
which satisfies the following properties:

(O-1) {\bf identity}: If $L$ and $L'$ are line bundles over $X$ and isomorphic (i.e., if $f:L\to X$ and $f':L' \to X$, then there exists an isomorphism $i: L\to L'$ such that $f = f' \circ i$) , then we have
$$\phi(L) = \phi(L'): \bB(X  \xrightarrow{f}  Y) \to \bB(X  \xrightarrow{f}  Y).$$

(O-2)  {\bf commutativity}: Let $L$ and $L'$ be two fiber-objects  over $X$, then we have
$$\phi(L) \circ \phi(L') = \phi(L') \circ \phi(L) :\bB(X  \xrightarrow{f}  Y) \to \bB(X  \xrightarrow{f}  Y). $$

(O-3)   {\bf compatibility with product}: For morphisms $f:X \to Y$ and $g:Y \to Z$,  $\alp \in \bB(X  \xrightarrow{f} Y)$ and $ \be \in \bB(Y  \xrightarrow{g} Z)$,  a fiber-object  $L$ over $X$ and a fiber-object  $M$ over $Y$, we have
 $$ \phi(L) (\alp \bullet \be) = \phi(L)(\alp) \bullet \be, \quad  \phi(f^*M) (\alp \bullet \be) = \alp \bullet \phi(M)(\be).$$

(O-4)  {\bf compatibility with pushforward}: For a confined morphism $f:X \to Y$ and a fiber-object $M$ over $Y$ we have 
$$ f_*\left (\phi(f^*M)(\alp) \right ) = \phi(M)(f_*\alp).$$

(O-5)   {\bf compatibility with pullback}: For an independent square and a fiber-object  $L$ over $X$
$$
\CD
X' @> g' >> X \\
@V f' VV @VV f V\\
Y' @>> g > Y 
 \endCD
$$
we have 
$$g^*\left (\phi(L)(\alp) \right ) = \phi({g'}^*L)(g^*\alp).$$

The above operator is called an ``{\it orientation}" and a bivariant theory equipped with such an orientation is called an {\it oriented bivariant theory}, denoted by $\bOB$. 
\item An {\it oriented Grothendieck transformation} between two oriented bivariant theories is a Grothendieck transformation which preserves or is compatible with the operator, i.e., for two oriented bivariant theories $\bOB$ with an orientation $\phi$ and $\bOB'$ with an orientation $\phi'$ the following diagram commutes
$$
\CD
\bOB (X  \xrightarrow{f}  Y)  @> {\phi(L)}>> \bOB (X  \xrightarrow{f}  Y) \\
@V \ga VV @VV \ga V\\
\bOB' (X  \xrightarrow{f}  Y) @>>{\phi'(L)} > \bOB' (X  \xrightarrow{f}  Y).
 \endCD
$$
\end{enumerate}
\end{defn} 

\begin{thm} (\cite[Theorem 4.6]{Yo1}) \label{obt}(A universal oriented bivariant theory)
 Let  $\Cal V$ be a category with a class $\Cal C$ of confined morphisms, a class of independent squares, a class  $\Cal S$ of specialized morphisms and $\Cal L$ a fibered category over $\Cal V$.  We define 
$$\bOM^{\Cal C} _{\Cal S}(X  \xrightarrow{f}  Y)$$
to be the free abelian group generated by the set of isomorphism classes of cobordism cycles over $X$
$$[V  \xrightarrow{h} X; L_1, L_2, \cdots, L_r]$$
such that $h \in \Cal  C$, $f \circ h: W \to Y \in \Cal S$ and $L_i$ a fiber-object over $V$.
\begin{enumerate}
\item The association $\bOM^{\Cal C} _{\Cal S}$ becomes an oriented bivariant theory if the four operations are defined as follows:
\begin{enumerate}

\item  {\bf Orientation $\Phi$}: For a morphism $f:X \to Y$ and a fiber-object $L$ over $X$, the operator
$$\phi (L):\bOM^{\Cal C} _{\Cal S} ( X  \xrightarrow{f}  Y) \to \bOM^{\Cal C} _{\Cal S} ( X  \xrightarrow{f}  Y) $$
is defined by
$$\phi(L)([V  \xrightarrow{h} X; L_1, L_2, \cdots, L_r]):=[V  \xrightarrow{h} X; L_1, L_2, \cdots, L_r, h^*L].$$
and extended linearly.

\item {\bf Product}: For morphisms $f: X \to Y$ and $g: Y
\to Z$, the product operation
$$\bullet: \bOM^{\Cal C} _{\Cal S} ( X  \xrightarrow{f}  Y) \otimes \bOM^{\Cal C} _{\Cal S} ( Y  \xrightarrow{g}  Z) \to
\bOM^{\Cal C} _{\Cal S} ( X  \xrightarrow{gf}  Z)$$
is  defined as follows: The product is defined by
\begin{align*}
& [V  \xrightarrow{h} X; L_1, \cdots, L_r]  \bullet [W  \xrightarrow{k} Y; M_1, \cdots, M_s] \\
& :=  [V'  \xrightarrow{h \circ {k}''}  X; {{k}''}^*L_1, \cdots,{{k}''}^*L_r, (f' \circ {h}')^*M_1, \cdots, (f' \circ {h}')^*M_s ]\
\end{align*}
and extended bilinearly. Here we consider the following fiber squares
$$\CD
V' @> {h'} >> X' @> {f'} >> W \\
@V {{k}''}VV @V {{k}'}VV @V {k}VV\\
V@>> {h} > X @>> {f} > Y @>> {g} > Z .\endCD
$$
\item {\bf Pushforward}: For morphisms $f: X \to Y$
and $g: Y \to Z$ with $f$ confined, the pushforward operation
$$f_*: \bOM^{\Cal C} _{\Cal S} ( X  \xrightarrow{gf} Z) \to \bOM^{\Cal C} _{\Cal S} ( Y  \xrightarrow{g}  Z) $$
is  defined by
$$f_*\left ([V  \xrightarrow{h} X; L_1, \cdots, L_r]  \right) := [V  \xrightarrow{f \circ h} Y; L_1, \cdots, L_r]$$
and extended linearly.

\item {\bf Pullback}: For an independent square
$$\CD
X' @> g' >> X \\
@V f' VV @VV f V\\
Y' @>> g > Y, \endCD
$$
the pullback operation
$$g^*: \bOM^{\Cal C} _{\Cal S} ( X  \xrightarrow{f} Y) \to \bOM^{\Cal C} _{\Cal S}( X'  \xrightarrow{f'} Y') $$
is  defined by
$$g^*\left ([V  \xrightarrow{h} X; L_1, \cdots, L_r] \right):=  [V'  \xrightarrow{{h}'}  X'; {g''}^*L_1, \cdots, {g''}^*L_r]$$
and extended linearly, where we consider the following fiber squares:
$$\CD
V' @> g'' >> V \\
@V {h'} VV @VV {h}V\\
X' @> g' >> X \\
@V f' VV @VV f V\\
Y' @>> g > Y. \endCD
$$
\end{enumerate}

\item Let $\Cal {OBT}$ be a class of oriented bivariant theories $\bOB$ on the same category $\Cal V$ with a class $\Cal C$ of confined morphisms, a class of independent squares, a class $\Cal S$ of specialized morphisms and a fibered category $\Cal L$ over $\Cal V$. Let $\Cal S$ be nice canonically $\bOB$-oriented for any oriented bivariant theory $\bOB \in \Cal {OBT}$. Then, for each oriented bivariant theory $\bOB \in \Cal {OBT}$ with an orientation $\phi$  there exists a unique oriented Grothendieck transformation
$$\ga_{\bOB} : \bOM^{\Cal C} _{\Cal S} \to \bOB$$
such that for any $f: X \to Y \in \Cal S$ the homomorphism
$\ga_{\bOB} : \bOM^{\Cal C} _{\Cal S}(X  \xrightarrow{f}  Y) \to \bOB(X  \xrightarrow{f}  Y)$
satisfies the normalization condition that $$\ga_{\bOB}([X  \xrightarrow{\op {id}_X}  X; L_1, \cdots, L_r]) = \phi(L_1) \circ \cdots \circ \phi(L_r) (\theta_{\bOB}(f)).$$
\end{enumerate}
\end{thm}

\begin{rem} When we consider algebraic cobordism, the above fibered category $\Cal L$ is the category of line bundles as we deal from now on.
\end{rem}

\section{Bivariant algebraic cobordism with vector bundles}\label{Char0ConstructionSection} 

In this section, motivated by algebraic cobordism $\omega_{*,*}(X)$ of vector bundles studied in \cite{LeeP}, we consider a bivariant analogue $\Omega^{*,*}(X \xrightarrow f Y)$ of $\omega_{*,*}(X)$ in such a way that its covariant part, i.e., 
$\Omega^{-*,*}(X \to pt)$ is the same as Lee--Pandharipande's $\omega_{*,*}(X)$. Derived algebraic geometry is an essential tool for this section, mostly because we need to abuse the functorial properties of homotopy fibre products (see Definition \fref{DerivedFibreProduct}).

\subsection{Construction of $\Omega^{*,*}$}\label{ConstructionOfCobordismOfBundles}

Let us begin with the free $\Laz$-modules $\mathcal{M}^{i,r}_+(X \xrightarrow f Y)$ which are generated by the cobordism cycles of the form
\begin{equation*}
[V \xrightarrow h X, E]
\end{equation*}
where $V$ is connected, $h: V \to X$ is proper and the composite $f \circ h: V \to Y$ is quasi-smooth of virtual relative dimension $-i$, and $E$ is a vector bundle of rank $r$ on $V$. 

\begin{rem}
If we drop the condition \emph{connectedness} on the source variety $V$ in the above definition of $\mathcal{M}^{i,r}_+(X \xrightarrow f Y)$, then we need to modding out the free $\Laz$-modules $\mathcal{M}^{i,r}(X \xrightarrow f Y)$, generated by the cobordism cycles of the form $[V \xrightarrow h X, E]$ where $h: V \to X$ is proper and the composite $f \circ h: V \to Y$ is quasi-smooth, by the additvity relation $\sim_{+}$
$$[V \xrightarrow h X, E] := \sum_j [V_j \xrightarrow {h_j} X, E_j]$$
where $V = \sqcup V_j$ is the disjoint union of irreducible components and $h_j$ is the restriction of $h$ to $V_j$ and $E_j=E|_{V_j}$ is the restriction of $E$ to $V_j$. 
\end{rem}

\begin{rem} If we use the notation simliar to those in Theorem \ref{obt}, $\mathcal{M}^{i,r}_+(X \xrightarrow f Y)$ should be something like
$(\mathbb M^{\mathscr Prop}_{\mathscr {QS}})^{i,r}_+(X \xrightarrow f Y)$ where $\mathscr Prop$ referes to ``proper'' and $\mathscr {QS}$ refers to ``quasi-smooth'', but in order to avoid messy notation we just denote it simply as above.
\end{rem}

These groups form two kinds of bivariant theories (due to defining two kinds of bivariant products) by linearly extending the following operations:
\begin{enumerate}
\item {\bf Pushforward}: 
Let $f: X \to X'$ and $g:X' \to Y$ where $f$ is proper. We define the pushforward map $f_*: \mathcal{M}^{i,r}_+(X \xrightarrow {g \circ f} Y) \to \mathcal{M}^{i,r}_+(X' \to Y)$ by
\begin{equation*}
f_*([V \xrightarrow h X , E]) = [V \xrightarrow {f \circ h} X', E].
\end{equation*}
\item {\bf Pullback}: 
Suppose that we have a map $g:Y' \to Y$ and let $X'$ be the homotopy fibre product $Y' \times_Y^\R X$.
 Then we define the pullback map 
$$g^*: \mathcal{M}^{i,r}_+(X \xrightarrow f Y) \to \mathcal{M}^{i,r}_+(X' \xrightarrow {f'} Y')$$ 
by
\begin{equation*}
g^*([V \xrightarrow h X, E]) = [V' \xrightarrow {h'}, E'],
\end{equation*}
where $V' = Y' \times_X V$ and $E'=(g'')^*E$ is the pullback of $E$ 
by the induced map $g'':V' \to V$:
\begin{equation*}
\begin{tikzcd}
V' \arrow[]{r}{g''} \arrow[]{d}{h'} & V \arrow[]{d}{h} \\
X' \arrow[]{r}{g'} \arrow[]{d}{f'} & X \arrow[]{d}{f} \\
Y' \arrow[]{r}{g} & Y
\end{tikzcd}
\end{equation*}

\item {\bf Bivariant products $\bullet_{\oplus}$ and $\bullet_{\otimes}$}: We define the following two bivariant products
\begin{equation*}
\bullet_{\oplus}: \mathcal{M}^{i,r}_+(X \xrightarrow f Y) \times \mathcal{M}^{j,s}_+(Y \xrightarrow g Z) \to \mathcal{M}^{i+j,r+s}_+(X \xrightarrow {g \circ f} Z)
\end{equation*}
\begin{equation*}
\bullet_{\otimes}: \mathcal{M}^{i,r}_+(X \xrightarrow f Y) \times \mathcal{M}^{j,s}_+(Y \xrightarrow g Z) \to \mathcal{M}^{i+j,rs}_+(X \xrightarrow {g \circ f} Z)
\end{equation*}
as follows. Suppose we have cobordism cycles $[V \xrightarrow h X, E]$ and $[W \xrightarrow k Y, F]$ and form the following homotopy Cartesian diagram
\begin{equation*}
\begin{tikzcd}
V' \arrow[]{r}{h'} \arrow[]{d}{k''} & X' \arrow[]{r}{f'} \arrow[]{d}{k'} & W \arrow[]{d}{k} \\
V \arrow[]{r}{h} & X \arrow[]{r}{f} & Y \arrow[]{r}{g} & Z.
\end{tikzcd}
\end{equation*}


Now we define these two products by
\begin{equation*}
[V \xrightarrow h X, E] \bullet_{\oplus} [W \xrightarrow k Y , F] = [V' \xrightarrow {h \circ k''} X, E' \oplus F'],
\end{equation*}
\begin{equation*}
[V \xrightarrow h X, E] \bullet_{\otimes} [W \xrightarrow k Y , F] = [V' \xrightarrow {h \circ k''} X, E' \otimes F'],
\end{equation*}
where $E'=(k'')^*E$ and $F'=(f' \circ h')^*F$ are pullbacks onto $V'$ of $E$ and $F$ respectively.
\end{enumerate} 

\begin{prop} $\mathcal{M}^{*,*}_+(X \to  Y) $ is a commutative bivariant theory with respect to both products $\bullet_{\oplus}$ and $\bullet_{\otimes}$.
\end{prop}
\begin{proof} Its proof is basically the same as that of Theorem \ref{obt} (\cite[Theorem 4.6]{Yo1}) and the only difference between these two setups is that in Theorem \ref{obt} we consider a finite set $\{L_1, L_2, \cdots, L_r\}$ of line bundles instead of one vector bundle $E$ in the present setup.
\end{proof}

We can also choose natural orientations along quasi-smooth morphisms: indeed, if $f: X \to Y$ is quasi-smooth of relative virtual dimension $-i$, then one may 
define the two orientations
\begin{align*}
\theta_\oplus(f) &:= [X \xrightarrow {\op{id}_X} X,0] \in \mathcal{M}^{i,0}_+ (X \to Y) \\
\theta_\otimes(f) &:= [X \xrightarrow {\op{id}_X} X,\OO_X] \in \mathcal{M}^{i,1}_+ (X \to Y)
\end{align*}
for the bivariant theories $(\mathcal{M}^{*,*}_+, \bullet_\oplus)$ and $(\mathcal{M}^{*,*}_+, \bullet_\otimes)$ respectively. Note that the choice of the vector bundle is --- in both cases --- essentially forced upon us by the requirement that the orientation of the identity morphism $X \to X$ should be the multiplicative identity of the ring $\mathcal{M}^{*,*}_+(X \to X)$.

Like for any other bivariant theory, we have the associated homology and cohomology theories defined by
\begin{equation*}
\mathcal{M}_{*,*}^+(X) := \mathcal{M}^{-*,*}_+(X \to pt)
\end{equation*}
and
\begin{equation*}
\mathcal{M}^{*,*}_+(X) := \mathcal{M}^{*,*}_+(X \stackrel{\mathrm{id}}{\to} X)
\end{equation*}
respectively. 
For a quasi-smooth morphism $f: X \to Y$ of relative virtual dimension $-i$ we 
have the Gysin pullback homomorphism $f^!=\theta_\oplus(f) \bullet_{\oplus}:\mathcal{M}^{+}_{k,r}(Y) \to \mathcal{M}^{+}_{k-i,r}(X)$ and  
the Gysin pushforward homomorphism $f_!=f_{*}(- \bullet_{\oplus} \theta_\oplus(f)):\mathcal{M}_{+}^{k,r}(X) \to \mathcal{M}_{+}^{k+i,r}(Y)$.

\begin{rem}
If we use the bivariant product $\bullet_{\otimes}$, then the corresponding Gysin homomorphisms are respectively
 $f^!=\theta_\otimes(f) \bullet_{\otimes }:\mathcal{M}^{+}_{k,r}(Y) \to \mathcal{M}^{+}_{k-i,r}(X)$ and $f_!=f_{*}(- \bullet_{\otimes} \theta_\otimes(f)):\mathcal{M}_{+}^{k,r}(X) \to \mathcal{M}_{+}^{k+i,r}(Y)$.
\end{rem}

\begin{defn}\label{TopChernClass}
Given a vector bundle $E$ of rank $r$ on $X$, we may define its \emph{top Chern class} as
\begin{equation*}
c_r(E) = s^*s_!(1_X) \in \mathcal{M}_{+}^{r,0}(X) 
\end{equation*}
where $s: X \to E$ is the zero section. One can use this to define the top Chern class operator on $\mathcal{M}^{+}_{*,*}(X)$ as $c_r(E) \bullet_\oplus:\mathcal{M}^{+}_{k,s}(X) \to \mathcal{M}^{+}_{k-r,s}(X)$.
\end{defn}

\begin{rem}
Unwinding the above definition, we see that the top Chern class $c_r(E)$ of a vector bundle $E$ of rank $r$ is the class $[V(0) \to X, 0]$ of the inclusion of the derived vanishing locus of the zero section of $E$. Hence, by Proposition \fref{DerivedVanishingLocus}, it can be identified with the class of
$$\sSpec(\Sym_X^*(E^\vee[-1])) \to X.$$
Unless $E$ is the vector bundle of rank $0$, this does \emph{not} correspond to the identity element in $\mathcal{M}^{0,0}(X)$: for example, the relative virtual dimension (see Section \fref{CotangentComplexSection}) of the above morphism is $-r$ (this follows from Example \fref{CotangentComplexOfVanishingLoci}), so it is not even in the right degree!
\end{rem}

From here up to the section where we deal with the bivariant product $\bullet_{\otimes}$,  we deal with the bivariant product $\bullet_{\oplus}$ and for the sake of simplicity we denote just $\bullet$ without the suffix $\oplus$.

To obtain 
our bivariant algebraic cobordism of vector bundles 
$\Omega^{*,*}(X \xrightarrow f Y)$ from the above $\mathcal{M}^{*,*}_+(X \xrightarrow f Y)$, we use similar relations as in \cite{An}, whose origin lies in the paper \cite{LS} of Lowrey and Sch\"urg. Let us quickly recall the construction of the bivariant derived algebraic cobordism $\Omega^*$ from \cite{An} using the language introduced in this paper. We start with the bivariant theory $\mathcal{M}_+^*$ so that $\mathcal{M}_+^*(X \to Y)$ is the free graded $\Laz$-module on generators $[V \to X]$, with $V$ connected, $V \to X$ proper and the composition $V \to Y$ quasi-smooth. In other words, the bivariant theory $\mathcal{M}_+^*$ can be naturally identified with the bivariant theory $\mathcal{M}_+^{*,0}$ as there is only one vector bundle of rank 0. Moreover, the algebraic bordism groups $d\Omega_*$ of Lowrey and Sch\"urg are, by construction, expressible as quotients of the homology $\Laz$-modules $\mathcal{M}^+_*(X) := \mathcal{M}_+^{-*}(X \to pt)$; we denote by $\mathcal{R}^\mathrm{LS}$ the bivariant subset of 
$\mathcal{M}^*_+$ 
so that $\mathcal{R}^\mathrm{LS} (X \to pt)$ is defined to be the kernel of $\mathcal{M}_+^{-*}(X \to pt)=
\mathcal{M}^+_*(X) \to d\Omega_{*}(X)=\Omega^{-*}(X \to pt)$ and $\mathcal{R}^\mathrm{LS} (X \to Y)$ is defined to be empty whenever $Y$ is not a point. The bivariant algebraic cobordism $\Omega^*$ 
constructed 
in \cite{An} is the quotient theory $\mathcal{M}_+^* / \langle \mathcal{R}^\mathrm{LS} \rangle_{\mathcal{M}_+^*}$, where we recall that $\langle \mathcal{R}^\mathrm{LS} \rangle_{\mathcal{M}_+^*}$ is the bivariant ideal generated by the bivariant set $\mathcal{R}^\mathrm{LS}$ inside $\mathcal{M}_+^*$.

\begin{defn}\label{CobOfBundlesDef}
We define the \emph{bivariant algebraic cobordism with vector bundles} as
$$\Omega^{*,*} := \mathcal{M}_+^{*,*}/\langle \mathcal{R}^\mathrm{LS} \rangle_{\mathcal{M}_+^{*,*}},$$
where we regard $\mathcal{R}^\mathrm{LS}$ as a bivariant subset of $\mathcal{M}_+^{*,*}$ via the identification $\mathcal{M}_+^* = \mathcal{M}_+^{*,0}$.
\end{defn}

Another way of phrasing Definition \fref{CobOfBundlesDef}, after recalling the details of how $d\Omega_*$ is obtained from $\mathcal{M}^+_*$ in \cite{LS}, consists of the following three steps:

\begin{enumerate}
\item \emph{Homotopy fibre relation:} Let $\mathcal{R}^\mathrm{fib}$ be the bivariant ideal of $\mathcal{M}^{*,*}_+$ generated by the elements
\begin{equation*}
[W_0 \to X, 0] - [W_\infty \to X, 0] \in 
\mathcal{M}^{*,0}_+(X \to pt),
\end{equation*}
where $W_0 \to X$ and $W_\infty \to X$ are obtained from a proper map $W \to \Proj^1\times X$ with $W$ quasi-smooth as homotopy fibres over the constant inclusions $0 \times X \hookrightarrow \Proj^1 \times X$ and $\infty \times X \hookrightarrow \Proj^1 \times X$ respectively. We denote the bivariant theory $\mathcal{M}^{*,*}_+ / \mathcal{R}^\mathrm{fib}$ by $\Omega_\mathrm{naive}^{*,*}$.

\item \emph{Formal group law relation:} Let $\mathcal{R}^\mathrm{fgl}$ be the bivariant ideal of $\Omega_\mathrm{naive}^{*,*}$ generated by
\begin{equation*}
c_1(\Li_1 \otimes \Li_2) 1_X - F(c_1(\Li_1),c_1(\Li_2)) 1_X \in \Omega^{*,*}_\mathrm{naive}(X \to pt)
\end{equation*}
where $X$ is a smooth scheme, $\Li_1$ and $\Li_2$ are globally generated line bundles on $X$, and $F$ is the universal formal group law of the Lazard ring. The careful reader should be worried at this point, as a 
formal power series $F$ gives an infinite sum. However, the generating relations above are well defined elements, as the first Chern classes of globally generated line bundles act nilpotently on $\Omega^\mathrm{naive}_{*,*}$ (for a proof, see Proposition 3.15 in \cite{LS}). We denote the bivariant theory $\Omega_\mathrm{naive}^{*,*}/\mathcal{R}^\mathrm{fgl}$ by $\underline{\Omega}^{*,*}$

\item \emph{Strict normal crossing relation:} Suppose $D$ is a strict normal crossing divisor on a smooth scheme $W$ with prime divisors $D_1,...,D_r$ with multiplicities $n_1,...,n_r$. The formal group law relation allows us to express $[D \to W, 0]$ as a $\Laz$-linear combination of inclusions of the prime divisors $D_i$ and their intersections. The latter expression has an obvious lift $\xi_{D,W} \in \underline{\Omega}_{*,*}(D)$. Denote by $\mathcal{R}^\mathrm{snc}$ the bivariant ideal of $\underline{\Omega}^{*,*}$ generated by
\begin{equation*}
1_D - \xi_{D,W} \in \underline{\Omega}^{*,*}(D \to pt)
\end{equation*} 
for all strict normal crossing divisor inclusions $D \subset W$ with $W$ smooth. The kernel $\mathcal R '$ of the natural Grothendieck transformation $\mathcal{M}_+^{*,*} \to \underline{\Omega}^{*,*}/\mathcal{R}^\mathrm{snc}$ is by construction $\langle \mathcal{R}^\mathrm{LS} \rangle$, and therefore we obtain $\Omega^{*,*}$ in three steps. We denote the Grothendieck transformation $\mathcal{M}_+^{*,*} \to \Omega^{*,*}$ (which is the quotient homomorphism $\mathcal{M}_+^{*,*}(X \to Y) \to \Omega^{*,*}(X \to Y)$ for each $X \to Y$) by $\Theta$ for later use.
\end{enumerate}

\begin{rem}\label{RelationsWithBundlesRem}
One might protest that these are obviously wrong relations as none of the generating relations above use the vector bundles in the cycles in any nontrivial way. However, let us consider the relations imposed on $\Omega^{*,*}(X \to Y)$ by the homotopy fibre relation above. Given a proper map $W \to \Proj^1 \times X$ so that the composition $W \to \Proj^1 \times Y$ is quasi-smooth, and a vector bundle $E$ on $W$, we would like to have 
$$[W_0 \to X, E\vert_{W_0}] = [W_\infty \to X, E\vert_{W_\infty}] \in \Omega^{*,*}(X \to Y),$$
where $W_0$ and $W_\infty$ are the fibres of $W \to \Proj^1$ over $0$ and $\infty$ respectively. As in \cite{An} Proposition 3.7, one can prove that the relation
$$[W_0 \to W, 0] = [W_\infty \to W, 0] \in \Omega^{*,0}(W \stackrel {\mathrm{id}_W} \longrightarrow W)$$
holds, and one may multiply them by the class $[W \xrightarrow {\mathrm{id}_W} W, E] \in \Omega^{*,*}(W \xrightarrow {\mathrm{id}_W} W)$, then we have
\begin{align*}
[W_0 \to X, E\vert_{W_0}] &= [W \xrightarrow {\mathrm{id}_W} W, E] \bullet [W_0 \to W, 0], \\
[W_\infty \to X, E\vert_{W_\infty}] &= [W \xrightarrow {\mathrm{id}_W} W, E] \bullet [W_\infty \to W, 0],
\end{align*}
and see that 
$$[W_0 \to W, E\vert_{W_0}] = [W_\infty \to W, E\vert_{W_\infty}] \in \Omega^{*,*}(W \stackrel {\mathrm{id}_W} \longrightarrow W).$$
 As the map $W \to Y$ is quasi-smooth, we can multiply the above relation from the right using $\theta(W \to Y) =[W \xrightarrow {\mathrm{id}_W} W, 0] \in \Omega^{*,0}_+(W \to Y)$ to see that
$$[W_0 \to W, E\vert_{W_0}] = [W_\infty \to W, E\vert_{W_\infty}] \in \Omega^{*,*}(W \to Y).$$
Finally, the map $W \to X$ is proper, and hence we can push forward the above relation to obtain
$$[W_0 \to X, E\vert_{W_0}] = [W_\infty \to X, E\vert_{W_\infty}] \in \Omega^{*,*}(X \to Y)$$
which is exactly what we wanted.

In fact, linear combinations of cycles of the above form can be checked to be stable under all bivariant operations: that is to say, the bivariant ideal $\mathcal{R}^{\mathrm{fib}}$ admits the description that the $\Laz$-modules $\mathcal{R}^{\mathrm{fib}}(X \to Y)$ are generated by cycles of form $[W_0 \to X, E\vert_{W_0}] - [W_\infty \to X, E\vert_{W_\infty}]$.
There is a similar simple description of the bivariant ideal $\mathcal{R}^{\mathrm{fgl}}$ that can be proven using similar tricks. However, the third snc relations do not seem to admit any simpler description than the general one given by 
Proposition \fref{GeneratedBivariantIdeal}.
\end{rem}

\begin{rem}
We should note that the homology groups $\Omega_\mathrm{naive}^{*,0}(X \to pt)$ do \emph{not} agree with the naive bordism groups $d \Omega^\mathrm{naive}_{*}(X)$ of Lowrey and Sch\"urg. The problem is that the naive bordism groups $d \Omega^\mathrm{naive}_{*}(X)$ are merely Abelian groups, not $\Laz$-modules. In fact, it can be shown (but we will not do it here) that $\Omega_\mathrm{naive}^{-*,0}(X\to pt) \cong \Laz \otimes_\Z \Omega^\mathrm{naive}_{*}(X)$. On the other hand, as $d \Omega^\mathrm{pre}_*$ is defined as the quotient of $\Laz \otimes_\Z \Omega^\mathrm{naive}_{*}$ by the formal group law relations, it is not hard to show that there is 
a natural isomorphism $d \Omega^\mathrm{pre}_*(X) \cong \underline{\Omega}^{-*,0}(X \to pt)$. Finally, we note that the groups $d\Omega_*(X)$ and $\Omega^{-*,0}(X \to pt)$ are isomorphic, which will follow from the equivalence $\Omega^{*,0} \cong \Omega^*$ of bivariant theories.
\end{rem}

Before continuing, we will study in greater detail the bivariant ideal $\langle \mathcal{R}^\mathrm{LS} \rangle$ of $\mathcal{M}^{*,*}_+$ defining $\Omega^{*,*}$. Recall that $\mathcal{R}^\mathrm{LS}$ is the subset consisting of the kernel of the induced map $\mathcal{M}^+_{*,0} \to d\Omega_*$ of homology theories, where the right hand side is the derived algebraic bordism of Lowrey-Sch\"urg. As all the homotopy Cartesian squares are independent, we obtain from Proposition \fref{GeneratedBivariantIdeal} an easy characterization of the elements of the generated ideal: they are linear combinations of elements of form $f_*(\alpha \bullet g^*(s) \bullet \beta)$, where $s \in \mathcal{R}^\mathrm{LS}$, $\alpha$ and $\beta$ are elements of the bivariant theory $\mathcal{M}^{*,*}_+$, $g$ is an arbitrary morphism and $f$ is a proper morphism of derived schemes. In fact, we can go even further. Let us denote by $[E]$ the class $[X \xrightarrow {\op{id}_X} X, E]$ in $\mathcal{M}_+^{*,*}(X \xrightarrow {\op{id}_X} X)$ (or the image of this class in any quotient theory of $\mathcal{M}_+^{*,*}$).

\begin{lem}\label{VectorBundleLemma}
Let $E$ be a vector bundle on $Y$, and let $\alpha \in \mathcal{M}_+^{*,*}(X \stackrel f \to Y)$. Then 
$$\alpha \bullet [E] = [f^*E] \bullet \alpha.$$
\end{lem}
\begin{proof}
This is a trivial consequence of the definition of the bivariant product.
\end{proof}

\begin{lem}\label{VectorBundleLemma 2}
The theory $\mathcal{M}_+^{*,*}$ consists of linear combinations of elements of the form 
$$f_*([E] \bullet \alpha),$$
where $\alpha \in \mathcal{M}_+^{*,0}$.
\end{lem}
\begin{proof}
This is clear, because 
$[X \stackrel f \to Y, E] \in \mathcal{M}_+^{*,*}(Y \xrightarrow g Z)$ can be written as 
$$f_*([E] \bullet [X \xrightarrow {\op{id}_X} X, 0]),$$
where $[E] \in \mathcal{M}_+^{*,*}(X \xrightarrow {\op{id}_X} X)$ and $[X \xrightarrow {\op{id}_X} X, 0] \in \mathcal{M}_+^{*,0}(X \xrightarrow {g \circ f} Z)$.
\end{proof}

\begin{prop}(cf. Proposition \fref{GeneratedBivariantIdeal})\label{RelationsInCobordismOfBundles}
The kernel of the Grothendieck transformation $\mathcal{M}^{*,*}_+ \to \Omega^{*,*}$ consists of linear combinations of elements of the form
$$f_*([E] \bullet \alpha_0 \bullet g^* (s) \bullet \beta_0)$$
where $f, g, \alpha_0 \in \mathcal{M}^{*,0}_+, \beta_0 \in \mathcal{M}^{*,0}_+$ and $s \in \mathcal{R}^\mathrm{LS}$ are as in the following diagram
\begin{center}
\begin{tikzpicture}[scale=2]
\node (A1) at (1,0) {$A$};
\node (A2) at (2,0) {$B=pt$};
\node (B0) at (0,1) {$A''$};
\node (B1) at (1,1) {$A'$};
\node (B2) at (2,1) {$B'$};
\node (B3) at (3,1) {$Y$};
\node (X) at (1,2) {$X$};
\path[every node/.style={font=\sffamily\small}]
(A1) edge[->] node[yshift=0.4cm,circle,draw,inner sep=0pt, minimum size=0.5cm]{$s$} (A2)
(B0) edge[->] node[yshift=0.4cm,circle,draw,inner sep=0pt, minimum size=0.5cm]{$\alpha_0$} (B1)
(B1) edge[->] (B2)
(B2) edge[->] node[yshift=0.4cm,circle,draw,inner sep=0pt, minimum size=0.5cm]{$\beta_0$} (B3)
(B1) edge[->] (A1)
(B2) edge[->] node[right]{$g$} (A2)
(B0) edge[bend left,->] node[above]{$f$} (X)
(X) edge[bend left,->] node[above]{$h$} (B3)
;
\end{tikzpicture}
\end{center}
and $E$ a vector bundle on $A''$. In the above diagram, the bottom square is assumed to be homotopy Cartesian and $f$ to be confined.
\end{prop}
\begin{proof}
It is clear that elements of the above form lie in $\langle \mathcal{R}^\mathrm{LS} \rangle_{\mathcal{M}_+^{*,*}}$, therefore it is enough to show that these elements form a bivariant ideal, i.e., they satisfy the three conditions (1), (2) and (3) of Definition \ref{BI}. The proof will follow closely that of Proposition \fref{GeneratedBivariantIdeal}

\begin{enumerate}
\item This is proven exactly as (1) in the proof of Proposition \fref{GeneratedBivariantIdeal}.

\item This is proven exactly as (2) in the proof of Proposition \fref{GeneratedBivariantIdeal}.

\item We are left to check that elements of the above form are closed under left and right multiplication. 
\begin{enumerate}
\item \emph{Right multiplication:} Consider $\gamma \in \mathcal{M}_+^{*,*}(Y \xrightarrow{h} Y')$. Using the bivariant axiom $(A_{12})$, one observes that
$$f_*([E] \bullet \alpha_0 \bullet g^* (s) \bullet \beta_0) \bullet \gamma = f_* \bigl( ([E] \bullet \alpha_0 \bullet g^* (s) \bullet \beta_0) \bullet \gamma \bigr).$$
Moreover, by Lemma \fref{VectorBundleLemma 2}, $\gamma$ can be assumed to be of form $i_*([F] \bullet \gamma_0)$, where $\gamma_0 \in \mathcal{M}_+^{*,0}(W \xrightarrow{h \circ i} Y')$, $i: W \to Y$ is proper and $F$ is a vector bundle on $W$. Consider now the diagram
\begin{center}
\begin{tikzpicture}[scale=2]
\node (T1) at (0,1) {$C''$};
\node (T2) at (1,1) {$C''$};
\node (T3) at (2,1) {$C'$};
\node (T4) at (3,1) {$D'$};
\node (T5) at (4,1) {$W$};
\node (B1) at (0,0) {$A''$};
\node (B2) at (1,0) {$A''$};
\node (B3) at (2,0) {$A'$};
\node (B4) at (3,0) {$B'$};
\node (B5) at (4,0) {$Y$};
\node (B6) at (5,0) {$Y'$};
\path[every node/.style={font=\sffamily\small}]
(T1) edge[->] node[yshift=0.4cm,circle,draw,inner sep=0.2pt, minimum size=0.5cm]{$[E']$} node[below]{$\mathrm{Id}_{C''}$} (T2)
(T2) edge[->] node[yshift=0.4cm,circle,draw,inner sep=0pt, minimum size=0.5cm]{$\beta'_0$} (T3)
(T3) edge[->] node[yshift=0.4cm,circle,draw,inner sep=0pt, minimum size=0.5cm]{$g'^*s$} (T4)
(T4) edge[->] node[yshift=0.4cm,circle,draw,inner sep=0pt, minimum size=0.5cm]{$\alpha'_0$} (T5)
(T1) edge[->] node[left]{$i'''$} (B1)
(T2) edge[->] node[left]{$i'''$} (B2)
(T3) edge[->] node[left]{$i''$} (B3)
(T4) edge[->] node[left]{$i'$} (B4)
(T5) edge[->] node[left]{$i$} (B5)
(B1) edge[->] node[yshift=-0.4cm,circle,draw,inner sep=0.2pt, minimum size=0.5cm]{$[E]$} node[above]{$\mathrm{Id}_{A''}$} (B2)
(B2) edge[->] node[yshift=-0.4cm,circle,draw,inner sep=0pt, minimum size=0.5cm]{$\beta_0$} (B3)
(B3) edge[->] node[yshift=-0.4cm,circle,draw,inner sep=0pt, minimum size=0.5cm]{$g^*s$} (B4)
(B4) edge[->] node[yshift=-0.4cm,circle,draw,inner sep=0pt, minimum size=0.5cm]{$\alpha_0$} (B5)
(B5) edge[->] node[above]{$h$} (B6)
;
\end{tikzpicture}
\end{center}
where every square is homotopy Cartesian, where the bivariant elements of the top row are bivariant pullbacks of the bivariant elements of the bottom row, and where $g' = i' \circ g$. We can now use the bivariant projection formula $(A_{123})$ (together with $(A_{13})$) to conclude that
$$([E] \bullet \alpha_0 \bullet g^* (s) \bullet \beta_0) \bullet i_*([F] \bullet \gamma_0) = i'''_* \bigl( [E'] \bullet \alpha_0' \bullet g'^*(s) \bullet \beta_0' \bullet ([F] \bullet \gamma_0) \bigr)$$
and using Lemma \fref{VectorBundleLemma} the right hand side of the previous equation is equal to 
$$i'''_* \bigl( [F'] \bullet [E'] \bullet \alpha_0' \bullet g'^*(s) \bullet \beta_0' \bullet \gamma_0 \bigr)$$
where we denote by $F'$ the vector bundle obtained as the pullback of $F$ along the composition 
$$C'' \to C' \to D' \to W.$$
To conclude, we have shown that
\begin{align*}
f_*([E] \bullet \alpha_0 \bullet g^* (s) \bullet \beta_0) \bullet \gamma &= f_* \Bigl( i'''_* \bigl( [F'] \bullet [E'] \bullet \alpha_0' \bullet g'^*(s) \bullet \beta_0' \bullet \gamma_0 \bigr) \Bigr) \\
&= (f \circ i''')_* \Bigl( [F' \oplus E'] \bullet \alpha_0' \bullet g'^*(s) \bullet (\beta_0' \bullet \gamma_0) \Bigr)
\end{align*}
which is of the desired form.
\item \emph{Left multiplication:} Consider now $\gamma \in \mathcal{M}^{*,*}_+(X' \xrightarrow{h} X)$. Recall that we may apply the bivariant projection formula $(A_{123})$ to the derived fibre diagram
\begin{equation*}
\begin{tikzcd}
C'' \arrow[]{r}{f'} \arrow[]{d}{h'} & X' \arrow[]{d}{h} \\
A'' \arrow[]{r}{f} & X \arrow[]{r} & Y
\end{tikzcd}
\end{equation*}
to conclude that
$$\gamma \bullet f_*([E] \bullet \alpha_0 \bullet g^*(s) \bullet \beta_0) = f'_*(f^*(\gamma) \bullet [E] \bullet \alpha_0 \bullet g^*(s) \bullet \beta_0).$$
By Lemma \fref{VectorBundleLemma 2}, we may assume $f^*(\gamma)$ to be of form $i_*([F] \bullet \gamma_0)$, where $i: D \to C''$ is proper, $F$ is a vector bundle on $D$ and $\gamma_0 \in \mathcal{M}_+^{*,0}(D \xrightarrow{h' \circ i} A'')$. But
\begin{align*}
f'_*\Bigl( i_*([F] \bullet \gamma_0) \bullet [E] \bullet \alpha_0 \bullet g^*(s) \bullet \beta_0 \Bigr) &= f'_*\Bigl( i_*\bigl([F] \bullet \gamma_0 \bullet [E] \bullet \alpha_0 \bullet g^*(s) \bullet \beta_0 \bigr) \Bigr) \quad \text{(by $(A_{12})$)}\\
&= (f' \circ i)_*\Bigl( [F] \bullet \gamma_0 \bullet [E] \bullet \alpha_0 \bullet g^*(s) \bullet \beta_0 \Bigr) \\
&= (f' \circ i)_*\Bigl( [F \oplus E'] \bullet (\gamma_0 \bullet \alpha_0) \bullet g^*(s) \bullet \beta_0 \Bigr)
\end{align*}
where the vector bundle $E'$ is the pullback of $E$ along $h' \circ i: D \to A''$. The element is of the desired form, and we are therefore done.
\end{enumerate}
\end{enumerate}
\end{proof}

\begin{rem}\label{CharOfRelations}
One may give a simpler (and a less precise) version of Proposition \ref{RelationsInCobordismOfBundles}. Indeed, it is now trivial that the kernel of the Grothendieck transformation $\mathcal{M}^{*,*}_+ \to \Omega^{*,*}$ consists of linear combinations of elements of form $f_*([E] \bullet r)$, where $r$ is in the kernel of $\mathcal{M}^{*,0}_+ \to \Omega^{*,0}$. This kernel can be naturally identified with the kernel of the quotient morphism $\mathcal{M}^{*}_+ \to \Omega^{*}$ in \cite{An}. This formulation is especially useful for comparing the bivariant theories $\Omega^{*,*}$ and $\Omega^*$.
\end{rem}

\subsection{Basic properties}

The following proposition is a direct consequence of the construction of $\Omega^{*,*}$.
\begin{prop}\label{prop-strong-orient}
The bivariant theory $\Omega^{*,*}$ is commutative and has strong orientations along smooth morphisms.
\end{prop}
\begin{proof}
The first claim is trivial. The latter claim follows from the fact that \emph{any commutative bivariant theory} on the category of (quasi-projective) derived $k$-schemes having proper maps as confined morphisms, all derived Cartesian squares as independent squares and orientations along quasi-smooth morphisms has strong orientations along smooth morphisms. This is a special case of Proposition 2.9 of \cite{An}, where the smooth maps are playing the role of \emph{specialized projections} (the essential fact needed here is that a section of a smooth morphism is quasi-smooth, which follows from Theorem \fref{QSProperties} (3)).

\end{proof}

\begin{defn}
We define the following two Grothendieck transformations: 
\begin{itemize}
\item attaching the zero bundle:
$$\mathscr Z:\Omega^*(X \xrightarrow f Y) \to \Omega^{*,0}(X \xrightarrow f Y)$$
 defined by $\mathscr Z([V \to X]):=[V \to X, 0]$,
\item the \emph{forgetful Grothendieck transformation} forgetting the vector bundle:
$$\mathscr F: \Omega^{*,*}(X \xrightarrow f Y) \to \Omega^*(X \xrightarrow f Y)$$ defined by
$\mathscr F([V \to X, E])= [V \to X].$
\end{itemize}
\end{defn}

Note that the Grothendieck transformations above restrict to give an isomorphism
\begin{equation}\label{CobOfBundlesExtendsCob}
\Omega^{*} \cong \Omega^{*,0}.
\end{equation} 
From this, it follows that the theory of Chern classes carries over from $\Omega^{*,0}$ to the whole theory (cf. \cite{An} Theorem 3.17):

\begin{prop}
Let $E$ be a vector bundle of rank $r$ over $X$. We now have the \emph{total Chern class} 
$$c(E) = 1_X + c_1(E) + \cdots + c_r(E) \in \Omega^{*,0}(X) = \Omega^{*,0}(X \xrightarrow {\op{id}_X} X),$$
where $c_i(E) \in \Omega^{i,0}(X)$. These classes satisfy the following basic properties:
\begin{enumerate}
\item \emph{Naturality in pullbacks:} given $f: X' \to X$  we have $f^* c(E) = c(f^*E)$. 
\item \emph{Whitney sum formula:} given a short exact sequence of vector bundles
$$0 \to E' \to E \to E'' \to 0$$
we have $c(E) = c(E') \bullet c(E'')$.
\item \emph{Normalization:} the top Chern class of Definition \fref{TopChernClass} agrees with $c_r(E)$.
\item \emph{Formal group law:} given \emph{any} derived scheme $X$, and \emph{any} two line bundles $\Li$ and $\Li'$ on $X$, we have
$$c_1(\Li \otimes \Li') = F(c_1(\Li), c_1(\Li'))$$
where $F$ is the universal formal group law of $\Laz$.
\item \emph{Nilpotency:} the Chern classes $c_i(E)$ are nilpotent with respect to the bivariant product.
\end{enumerate}
\end{prop}
Lastly, we recall the fact that top Chern classes admit a description as a derived zero locus of any global section of the vector bundle.
\begin{prop}[cf. Lemma 3.9 \cite{LS}]\label{SectionAxiom}
Let $E$ be a vector bundle of rank $r$ over $X$, and let $s$ be a global section of $E$. Then the top Chern class $c_r(E)$ of $E$ coincides with the class $[V(s) \to X, 0]$, where $V(s) \hookrightarrow X$ is the inclusion of the derived vanishing locus of $s$ to $X$.
\end{prop}
\begin{proof}
It is enough to show that the class does not depend on the section $s$, as Definition \fref{TopChernClass} of $c_r$ can be rephrased as taking the derived vanishing locus of the zero section. Suppose now we have two global sections $s_1$ and $s_2$ on $E$, and consider the vanishing locus of $x_0 s_1 + x_1 s_2$ of the vector bundle $E(1)$ on $\Proj^1 \times X$. The bivariant homotopy fibre relation identifies the pullbacks of this class to $0 \times X$ and to $\infty \times X$ with each other, which is exactly what we want as the first one coincides with $[V(s_2) \to X, 0]$ and the second one with $[V(s_1) \to X, 0]$.
\end{proof}
\begin{rem}\label{SectionAxiomForDivisors}
By Proposition 3.2.6. of \cite{Khan} any quasi-smooth closed embedding $D \hookrightarrow X$ appears (essentially uniquely) as the derived vanishing locus of a global section $s$ of a line bundle $\OO(D)$. Hence, as a special case of the above Proposition \ref{SectionAxiom}, we obtain the following equality 
$$[D \to X, 0] = c_1(\OO(D)) \in \Omega^{*,0}(X).$$ 
This is analogous to (although more general than) the \emph{section axiom} used in the original construction of algebraic cobordism by Levine and Morel \cite{LM}.
\end{rem}

\subsection{Relation to Lee-Pandharipande's algebraic cobrodism with vector bundles $\omega_{*,*}(X)$}

The purpose of this section is to show that for all quasi-projective derived schemes $X$, we have
$$\Omega_{*,*}(X) \cong \omega_{*,*}(tX),$$
showing that the bivariant theory $\Omega^{*,*}$ we have constructed in the beginning of this section is an extension of the homological theory $\omega_{*,*}$ of Lee-Pandharipande. We first show that there is a well defined morphism $\iota: \omega_{*,*}(tX) \to \Omega_{*,*}(X)$, which can be described on the level of cycles as $[V \to X, E] \mapsto [V \to X, E]$. That this map is well defined follows from the next lemma (after applying Poincar\'e duality to identify homology and cohomology of smooth schemes), which will also have some use later.

\begin{lem}[Derived double point cobordism, cf. \cite{LP} Section 3.3]\label{DerivedDoublePointRelations}
Let $X \to Y$ be a morphism of derived schemes, $W \to \Proj^1 \times X$ a projective morphism so that the composition $W \to \Proj^1 \times Y$ is quasi-smooth and let $E$ be a vector bundle on $W$. Denote by $W_0$ the fibre over $0$, and suppose that the fibre over $\infty$ is the sum of two virtual Cartier divisors $A$ and $B$. Then
\begin{equation}\label{DerivedDoublePointCobEq}
[W_0 \to X, E \vert_{W_0}] = [A \to X, E \vert_A] + [B \to X, E \vert_B] - [\Proj_{A \cap B} (\OO_W(A) \oplus \OO) \to X, E \vert_{A \cap B}]
\end{equation} 
in $\Omega^{*,*}(X \to Y)$.
\end{lem}

\begin{rem}\label{DerivedDoublePointRemark}
Note that the restriction of $\OO(A + B)$ to $W_\infty$ is trivial, since it $\OO(A + B)$ is the pullback of $\OO(1)$ on $\Proj^1$. Hence the restrictions of $\OO_W(A)$ and $\OO_W(B)$ to $A \cap B \subset A + B$ are inverses of each other. Recalling that the equivalence class of a projective bundle $\Proj(E)$ does not change if $E$ is tensored with a line bundle (see Proposition \fref{ProjectiveBundlesTensoredWithLineBundles}), the above implies that 
\begin{align*}
\Proj_{A \cap B} \bigl(\OO_W(A) \oplus \OO \bigr) &\simeq \Proj_{A \cap B} \bigl(\OO  \oplus \OO_W(B)\bigr) \\
&\simeq \Proj_{A \cap B} \bigl(\OO_W(B) \oplus \OO \bigr) 
\end{align*}
which shows that the equivalence class of the projective bundle appearing in (\fref{DerivedDoublePointCobEq}) does not depend on the choice of labels $A$ and $B$.
\end{rem}

\begin{proof}
Let us start with immediate reductions. It is enough to show that
\begin{align}\label{DerivedDoublePointCobRed1}
[W_0 \to W, E \vert_{W_0}] &= [A \to W, E \vert_A] + [B \to W, E \vert_B] \\ 
&- [\Proj_{A \cap B} (\OO_W(A) \oplus \OO) \to W, E \vert_{A \cap B}] \in \Omega^{*,*}(W): \notag
\end{align}
the formula (\fref{DerivedDoublePointCobEq}) then follows in a standard way using bivariant operations. But of course, (\fref{DerivedDoublePointCobRed1}) can be obtained from the formula
\begin{align}\label{DerivedDoublePointCobRed2}
[W_0 \to W] &= [A \to W] + [B \to W] - [\Proj_{A \cap B} (\OO_W(A) \oplus \OO) \to W] \in \Omega^*(W)
\end{align}
by multiplying with $[W \to W, E]$, so it is enough to look at the theory $\Omega^*$. Finally, we can use the formal group law relation to $W_\infty = A + B$ in order to deduce that
\begin{align*}
[W_0 \to W] &= [W_\infty \to W] \\
&= \sum_{i,j} a_{ij} c_1(\OO_W(A))^i c_1(\OO_W(B))^j \\
&= [A \to W] + [B \to W] + [A \cap B \to W] \bullet \sum_{i,j \geq 1} a_{ij} c_1(\OO_W(A))^{i-1} \bullet c_1(\OO_W(B))^{j-1} \\
&= [A \to W] + [B \to W] + i_!\Bigl( \sum_{i,j \geq 1} a_{ij} c_1(\OO_W(A)\vert_{A \cap B})^{i-1} \bullet c_1(\OO_W(B)\vert_{A \cap B})^{j-1} \Bigr) 
\end{align*}
where $i$ is the inclusion $A \cap B \hookrightarrow W$ and the last equality follows from the projection formula. Since the line bundles $\OO_W(A)\vert_{A \cap B}$ and $\OO_W(B)\vert_{A \cap B}$ are duals of each other (Remark \fref{DerivedDoublePointRemark}) we have reduced the proof of the theorem to proving the following claim: if $Z$ is a derived scheme and $\Li$ is a line bundle on $Z$, then
\begin{equation}\label{DerivedDoublePointCobRed3}
[\Proj_Z(\Li \oplus \OO) \to Z] = - \sum_{i,j \geq 1} a_{ij} c_1(\Li)^{i-1} \bullet c_1(\Li^\vee)^{j-1} \in \Omega^*(Z).
\end{equation}

In order to prove (\fref{DerivedDoublePointCobRed3}), choose a smooth scheme $T$ together with a line bundle $\Li'$, and a morphism $f: Z \to T$ so that $f^* \Li' \simeq \Li$. Since $\Omega^*(T)$ is known to be isomorphic to the algebraic bordism group $\Omega_{\dim(T) - *} (T)$ of Levine-Morel, we can compute that
\begin{align*}
[\Proj_Z(\Li \oplus \OO) \to Z] &= f^* ([\Proj_T(\Li' \oplus \OO) \to T]) \\
&= - f^* \Bigl( \sum_{i,j \geq 1} a_{ij} c_1(\Li')^{i-1} \bullet c_1(\Li'^\vee)^{j-1} \Bigr) \quad (\text{Lemma 3.3 of \cite{LP}})\\
&= \sum_{i,j \geq 1} a_{ij} c_1(\Li)^{i-1} \bullet c_1(\Li^\vee)^{j-1}
\end{align*}
proving the claim.
\end{proof}
The easier part of the proof is the surjectivity:

\begin{lem}\label{LeePSurj}
The map $\iota:  \omega_{*,*} (tX) \to \Omega_{*,*}(X)$ is surjective.
\end{lem}
\begin{proof}
Consider the cobordism cycle $[V \to V, E]$ in $\Omega_{*,*}(V)$. We have already remarked that $\Omega_{*,0}(V)$ agrees with the Lowrey-Sch\"urg algebraic bordism group of $V$, and as this is equivalent to the Levine-Morel (Levine-Pandharipande) algebraic bordism of the truncation $tV$, we know that the cycle $[V \to V, 0] \in \Omega_{*,0}(V)$ is equivalent to a cycle $\alpha$ which is a linear combination of smooth schemes mapping to $V$. Therefore $[V \to V, E] = [V \to V, E] \bullet \alpha$, and the right hand side is clearly contained in the image of $\omega_{*,*} (tX) \to \Omega_{*,*}(X)$. 
\end{proof}

We will show the injectivity of $\iota$ essentially in the same way as the injectivity of $\omega_{*,*}(pt) \otimes_\Laz \omega_*(X) \to \omega_{*,*}(X)$ is shown in \cite{LeeP}. Let us begin with a definition.

\begin{defn}[cf. Section 4.2 of \cite{LeeP}]\label{PairingRhoDef1}
Let $\Psi_{r,d}$ be the space of homogeneous degree $d$ integral coefficient polynomials in formal variables $c_1,...,c_r$, where $c_i$ is of degree $i$. Define a bilinear pairing
\begin{equation*}
\rho: \Psi_{r,d} \times \mathcal{M}^+_{*,r}(X) \to \Omega_{*-d}(X) 
\end{equation*}
by the formula  
$$(\Phi,[V \stackrel h \to X, E]) \mapsto h_* \Bigl (\Phi(c_1(E),...,c_r(E)) \bullet \Theta ( [V \xrightarrow {\op{id}_V} V] ) \Bigr ),$$
where $\Phi(c_1(E),...,c_r(E)) \in \Omega^{*,0}(V)$ is the value of the polynomial $\Phi$ evaluated at the Chern classes of $E$ and where $[V \xrightarrow {\op{id}_V} V] \in \Omega_{*}(V)$. (Recall that $\Theta$ is the surjection $\mathcal{M}_+^{*,*} \to \Omega^{*,*}$).
\end{defn}

Our strategy is to show that the pairing $\rho$ descends to a pairing
\begin{equation}\label{PairingRhoDef2}
\rho: \Psi_{r,d} \times \Omega_{*,r}(X) \to \Omega_{*-d}(X),
\end{equation}
which will then allow us to copy the rest of the proof from \cite{LP}. The desired descent result follows from a more general fact. To state it, we need another definition.

\begin{defn}
Let $\Phi \in \Psi_{r,d}$ be a polynomial. Consider the \emph{differentiation morphism}
$$\partial_\Phi: \mathcal{M}_+^{*,r}(X \xrightarrow f Y) \to \Omega^{*+d,r}(X \xrightarrow f  Y)$$
defined by
$$\partial_\Phi([V \xrightarrow h X, E]):= h_* \Bigl (\Phi(c_1(E),...,c_r(E)) \bullet \Theta ([V \xrightarrow {\op{id}_V} V, E] )\Bigr ),$$
where $\Phi(c_1(E),...,c_r(E)) \in \Omega^{*,0}(V)$ is as in Definition \fref{PairingRhoDef1}, and where $[V \xrightarrow {\op{id}_V} V, E] \in \mathcal{M}_+^{*,r}(V \xrightarrow {f \circ h} Y)$ and  $\Theta ([V \xrightarrow {\op{id}_V} V, E] )  \in \Omega^{*,r}(V \xrightarrow {f \circ h}  Y)$.
\end{defn}

If the maps $\partial_\Phi$ defined above descend to give operators $\Omega^{*,r}(X \xrightarrow f Y) \to \Omega^{*+d,r}(X \xrightarrow f  Y)$, the pairing $\rho$ defined earlier in (\fref{PairingRhoDef2}) is the homological special case of the pairing
$$\rho: \Psi_{r,d} \times \Omega^{*,r}(X \xrightarrow f  Y) \to \Omega^{*+d}(X \xrightarrow f  Y)$$
on bivariant groups defined by the formula
\begin{equation}\label{DefOfPairing}
(\Phi, \alpha) \mapsto \Cal F(\partial_\Phi(\alpha)),
\end{equation}
where we recall that $\Cal F: \Omega^{*+d,r}(X \xrightarrow f  Y) \to \Omega^{*+d}(X \xrightarrow f  Y)$ is the forgetful Grothendieck transformation. This, and a lot more, follows from the next result.

\begin{prop}\label{diff}
The differential operations 
\begin{equation*}
\partial_\Phi: \Omega^{*,*} \to \Omega^{*,*} 
\end{equation*}
are well defined, commute with pushforwards and pullbacks, and are linear over the subtheory $\Omega^{*,0} = \Omega^{*}$. Moreover, $\partial_{\Phi_1} \partial_{\Phi_2} = \partial_{\Phi_1 \Phi_2}$, the generating differentials $\partial_{c_n}$ satisfy a generalized Leibniz rule
\begin{equation*}
\partial_{c_n} (\alpha \bullet \beta) = \sum_{i+j=n} \partial_{c_i}(\alpha) \bullet \partial_{c_j}(\beta).
\end{equation*}
\end{prop} 
\begin{proof}
\begin{enumerate}
\item The well-definedness of $\partial_\Phi: \Omega^{*,*} \to \Omega^{*,*}$:
Using the characterization of the relations imposed on $\mathcal{M}_+^{*,*}$, we need only to show that 
the kernel of the Grothendieck transformation $\Theta: \mathcal{M}_+^{*,*} \to \Omega^{*,*}$ is sent to zero  by $\partial_{\Phi}$ ($\partial_\Phi$ is clearly $\Laz$-linear). Indeed, as is noted in Remark \ref{CharOfRelations}, the kernel of $\Theta$
is a linear combination of elements of the form $g_*([E] \bullet r)$ where $r$ is in the kernel of $\Theta: \mathcal{M}_+^{*,0} \to \Omega^{*,0}$. Hence  we have
\begin{align*}
\partial_\Phi(g_*([E] \bullet r)) & = g_*\Bigl (\Phi (E) \bullet \Theta ([E] \bullet r) \Bigr).\\
& = g_*\Bigl (\Phi (E) \bullet \Theta ([E]) \bullet \Theta (r) \Bigr)\\
& = 0 \qquad \text {(since $\Theta (r) =0$)}
\end{align*}
Thus $\partial_{\Phi}(\op{ker}\Theta)=0$. Therefore $\partial_\Phi: \Omega^{*,*} \to \Omega^{*,*}$ defined by 
$$\partial_\Phi \Bigl ( \Theta \bigl (g_*([E] \bullet \alp) \bigr ) \Bigr ):= g_*\Bigl (\Phi (E) \bullet \Theta ([E] \bullet \alp) \Bigr)$$
is well-defined. Here $\alp \in \mathcal{M}_+^{*,0}$ (see Lemma \ref {VectorBundleLemma 2}). 

\item $\partial_\Phi$ commutes with $f_*$: it suffices to show that for a generating element.
{\begin{align*}
(\partial_\Phi \circ f_*) \Bigl (\Theta \bigl (g_*([E] \bullet \alpha) \bigr )\Bigr ) & =  \partial_\Phi \Bigl (f_* \Bigl ( \Theta \bigl (g_*([E] \bullet \alpha) \bigr ) \Bigr ) \Bigr ) \\
& = \partial_\Phi \Bigl (\Theta  \Bigl ( f_*\bigl (g_*([E] \bullet \alpha) \bigr ) \Bigr ) \Bigr )  \quad  \text{(since $f_*\circ \Theta = \Theta \circ f_*$)} \\
& = \partial_\Phi \Bigl (\Theta  \Bigl ( (f \circ g)_*([E] \bullet \alpha) \bigr ) \Bigr ) \Bigr ) \\
& = (f \circ g)_* \Bigl (\Phi (E) \bullet \Theta ([E] \bullet \alpha) \Bigr ) \\
& = f_* \Bigl ( g_* \Bigl (\Phi (E) \bullet \Theta ([E] \bullet \alpha) \Bigr )  \Bigr) \\
&= f_* \Bigl (\partial_\Phi \Bigl ( \Theta \bigl (g_*([E] \bullet \alp) \bigr ) \Bigr ) \Bigr )\\
& = (f_* \circ \partial_\Phi) \Bigl (\Theta \bigl (g_*([E] \bullet \alp) \bigr ) \Bigr )
\end{align*}
For the sake of simplicity, from now on we drop off the symbol $\Theta$.
}

\item Commutativity with pullbacks is proven in a similar fashion. 

\item It is also trivial that $\partial_{\Phi_1} \partial_{\Phi_2} = \partial_{\Phi_1 \Phi_2}$.

\item $\partial_{c_n} (\alpha \bullet \beta) = \sum_{i+j=n} \partial_{c_i}(\alpha) \bullet \partial_{c_j}(\beta)$: This also gives an easy proof for the generalized Leibniz formula: recalling that
\begin{align*}
& f_*([E] \bullet \alpha) \bullet g_* ([F] \bullet \beta) \\
&= f_* \Bigl ([E] \bullet \alpha \bullet g_* ([F] \bullet \beta) \Bigr ) \\
&= f_* \Bigl (g'_* \Bigl ([E' \oplus F'] \bullet g^*(\alpha) \bullet \beta \Bigr ) \Bigr ) \quad \text {(by $A_{123}$ for $[E] \bullet \alpha \bullet g_* ([F] \bullet \beta) $ )}
\end{align*}
we notice that
\begin{align*}
& \partial_{c_n} \Bigl (f_*([E] \bullet \alpha) \bullet g_* ([F] \bullet \beta) \Bigr ) \\
&= \partial_{c_n} \Bigl (f_* \Bigl (g'_* \Bigl ([E' \oplus F'] \bullet g^*(\alpha) \bullet \beta \Bigr ) \Bigr ) \Bigr ) \\
&= f_* \Bigl ( \partial_{c_n} \Bigl (g'_* \Bigl ([E' \oplus F'] \bullet g^*(\alpha) \bullet \beta \Bigr ) \Bigr ) \Bigr ) \quad \text{ (since $\partial_{c_n} \circ f_* = f_* \circ \partial_{c_n}$)} \\
&= f_* \Bigl (g'_* \Bigl (c_n(E' \oplus F') \bullet [E' \oplus F'] \bullet g^*(\alpha) \bullet \beta \Bigr ) \Bigr ) \\
&= \sum_{i,j} f_* \Bigl (g'_* \Bigl (c_i(E') \bullet c_j(F') \bullet [E' \oplus F'] \bullet g^*(\alpha) \bullet \beta \Bigr ) \Bigr ) \\
&= \sum_{i,j} f_* \Bigl ( \bigl (c_i(E) \bullet [E] \bullet \alpha  \bigr ) \bullet g_* \bigl (c_i(F) \bullet [F] \bullet \beta  \bigr ) \Bigr ) \\
&= \sum_{i,j} f_* \bigl (c_i(E) \bullet [E] \bullet \alpha \bigr ) \bullet g_* \bigl (c_i(F) \bullet [F] \bullet \beta \bigr ) \\
&= \sum_{i,j} \partial_{c_i} \bigl (f_*([E] \bullet \alpha) \bigr ) \bullet \partial_{c_j} \bigl (g_*([F] \bullet \beta) \bigr )
\end{align*}
which is exactly what we wanted to prove. This also shows that if $\partial_{c_i}(\alpha) = 0$ for all $i>0$ (e.g. $\alpha \in \Omega^*$) then $\partial_\Phi(\alpha \bullet \beta) = \alpha \bullet \partial_\Phi(\beta)$ and $\partial_\Phi(\beta \bullet \alpha) = \partial_\Phi(\beta) \bullet \alpha$ for all $\Phi$ (they are polynomials of $c_i$, and $\partial_{c_0} = id$). This proves the claim about $\partial_\Phi$ being linear over the subtheory $\Omega^*$.
\end{enumerate}
\end{proof}

Now that we have shown that the action
\begin{equation*}
\Psi_{r,*} \times \Omega^{*,r} \to \Omega^{*}
\end{equation*}
is well defined, we would like to make sure that it coincides with the action of Lee and Pandharipande whenever both make sense.
\begin{lem}\label{CompatibilityLemma}
Let $X$ be a quasi-projective derived scheme. Now the natural map $\iota: \omega_{*,r}(tX) \to \Omega_{*,r}(X)$ induces a commutative square
\begin{equation*}
\begin{tikzcd}
\Psi_{r,d} \times \omega_{*,r}(tX) \arrow[]{r} \arrow[]{d} & \omega_{*-d}(tX)  \arrow[]{d}{\cong} \\
\Psi_{r,d} \times \Omega_{*,r}(X) \arrow[]{r} & \Omega_{*-d}(X)
\end{tikzcd}
\end{equation*}
where the upper horizontal map is the bilinear pairing $\rho^X$ defined in the Section 4.2 of \cite{LeeP}, and where the lower horizontal map is the pairing defined in equation (\fref{DefOfPairing}). In other words,
$$\rho^X(\Phi, \alpha) = \Cal F \bigl(\partial_\Phi(\iota(\alpha))\bigr)$$
for all $\alpha \in \omega_{*,r}(tX)$ and all $\Phi \in \Psi_{r,d}$.
\end{lem}
\begin{proof}
It is enough to show this for  $\alpha = [V \stackrel h \to tX, E] \in \omega_*(tX)$ as these generate the theory. Both sides of the square send the pair $(\Phi, \alpha)$ to
$$h_*\bigl(\Phi\bigl(c_1(E),...,c_r(E)\bigr) \cap 1_V\bigr).$$
The claim now follows from the fact that the action of the bivariant Chern classes on $\Omega^{*}(V \to pt)$ agrees with the action of the more classical Chern class operators of Levine and Morel that are used in \cite{LeeP} (this follows from \cite{An} Proposition 3.19).  
\end{proof}

We are finally ready to prove the the following injectivity result.

\begin{lem}\label{LeePInj}
The natural (cross product) map $\mathscr P: \omega_{*,*}(pt) \otimes_\Laz \Omega^*(X \to Y) \to \Omega^{*,*}(X \to Y)$ defined by 
$$[V \to pt, E] \otimes [W \to X] \mapsto [V \times W \to X, \op{pr}_1^* E]$$
is injective. Here $\op{pr}_1: V \times W \to V$ is the projection to the first factor.
\end{lem}

\begin{rem}\label{CrossProduct}
The above morphism admits the following alternative description which will be useful in the proof of the lemma. Given $\beta \in \omega_{*,*}(pt)$ and $\alpha \in \Omega^*(X \to Y) = \Omega^{*,0}(X \to Y)$, 
$\mathscr P(\beta \otimes \alpha)= \pi_X^*(\iota(\beta)) \bullet \alpha$, where $\pi_X$ is the structure morphism $X \to pt$.

Here, for the sake of later use, we also remark that $\pi_X^*(\iota(\beta)) \bullet \alpha$ is nothing but the \emph{cross product} $\iota(\beta) \times \alpha$ of $\iota(\beta)$ and $\alpha$. The cross product (see \cite[2.4 External products]{FM}) is defined as follows:
\begin{equation*}
\times : \mathbb B^i(X_1 \xrightarrow f Y_1) \otimes \mathbb B^j(X_2 \xrightarrow g Y_2) \to \mathbb B^{i+j}(X_1 \times X_2 \xrightarrow {f \times g} Y_1 \times Y_2)
\end{equation*}
is determined by the formula $\beta \times \alpha :=\op{pr}_1^*\beta \bullet \op{pr}_2^*\alpha$ where $\op{pr}_1:Y_1 \times X_2 \to Y_1$ is the projection to the first factor and $\op{pr}_2:Y_1 \times X_2 \to X_2$ is the projection to the second factor. Refer to the following commutative diagram:
$$
\xymatrix{
&&  Y_1 \times Y_2 \ar[rr]  && Y_2 \\
X_1 \times X_2 \ar[rr]_{f \times \op{id}_{X_2}} \ar[d] \ar[urr]^{f\times g} && Y_1 \times X_2 \ar[rr]_{\op{pr}_2} \ar[u]_{\op{id}_{Y_1} \times g} \ar[d]^{\op{pr}_1} && X_2 \ar[u]_g\\
X_1  \ar[rr]_{f} && Y_1   \\
}
$$
Note that in the above Lemma \ref{LeePInj} the structure morphism $\pi_X: X \to pt$ is nothing but the first factor projection $\op{pr}_1:X =pt \times X \to pt$ and the second factor projection $\op{pr}_2:X =pt \times X \to X$ is the same as the identity map $\op{id}_X$, and we consider the following cross product:
\begin{equation*}
\times : \Omega^{i,r}(pt \xrightarrow {} pt) \otimes \Omega^{j,0}(X \xrightarrow g Y) \to \Omega^{i+j,r}(X = pt \times X \xrightarrow g pt \times Y=Y)
\end{equation*}
defined by $\iota(\beta) \times \alpha$.
\end{rem}

\begin{proof}
This is basically the same as the proof of the Proposition 15 in \cite{LeeP}. Recall that the proof is based on the fact that $\omega_{*,r}(pt)$ is a free $\Laz$-module with a certain basis $e_1,e_2,...$, (for more details see \S 0.7 Basis and Theorem 2 of \cite{LeeP}) and that the $\Laz$-module $\Psi^\Laz_{r,*}$ of polynomials in variables $c_1,...,c_r$ has a dual basis $e_1^\vee,e_2^\vee,...$ for the pairing 
$$\rho = \rho^{pt}: \Psi^\Laz_{r,*} \times \omega_{*,r}(pt) \to \omega_{*}(pt) \cong \Laz$$ 
of Lee-Pandharipande defined by the formula
$$(\Phi, [V \stackrel f \to pt,E]) \mapsto f_*(\Phi(c_1(E),...,c_r(E)) \cap 1_V).$$ 
In other words,
$$\rho(e^\vee_j, e_i) = \delta_{ji}.$$

Hence, if $X$ is an arbitrary quasi-projective derived scheme, and if the image $\mathscr P(\alpha)$ of
$$\alpha := \sum_i e_i \otimes \alpha_i \in \omega_{*,*}(pt) \otimes_\Laz \Omega^*(X \to Y)$$
in $\Omega^{*,*}(X \to Y)$ vanishes, then as
\begin{align*}
\Cal F \bigl( \partial_{e_j^\vee}(\mathscr P(\alpha))\bigr) &= \Cal F \bigl( \partial_{e_j^\vee}(\mathscr P\Bigl(\sum_i e_i \otimes \alpha_i \Bigr) \bigr) \\
&= \Cal F \Bigl( \sum_i \partial_{e_j^\vee} \bigl( \pi_X^*(\iota(e_i)) \bullet \alpha_i \bigr)\Bigr) \\
&= \Cal F \Bigl( \sum_i \partial_{e_j^\vee} \bigl( \pi_X^*(\iota(e_i)) \bigr) \bullet \alpha_i \Bigr) \quad (\alpha_i \in \Omega^{*,0}) \\
&= \sum_i  \Cal F  \Bigl( \partial_{e_j^\vee} \bigl( \pi_X^*(\iota(e_i)) \bigr)\Bigr)  \bullet \alpha_i \quad (\Omega^{*,0} = \Omega^*) \\
&= \sum_i \pi_X^* \Bigl( \Cal F \bigl( \partial_{e_j^\vee}(\iota(e_i)) \bigr)\Bigr)  \bullet \alpha_i \quad (\text{$\pi_X^*$ commutes with $\partial_{e_j^\vee}$ and $\Cal F$}) \\
&= \sum_i \pi_X^* \bigl( \rho(e_j^\vee, e_i)\bigr)  \bullet \alpha_i \quad (\text{Lemma \fref{CompatibilityLemma}}) \\
&= \sum_i \pi_X^* \bigl( \delta_{ji} \bigr)  \bullet \alpha_i \\
&= \alpha_j
\end{align*}
$\alpha_j$ must be zero for all $j$. Hence $\alpha = 0$, and we are done.
\end{proof}

Combining Lemma \fref{LeePSurj} and Lemma \fref{LeePInj}, we obtain

\begin{thm}\label{LeePIso}
The map
$$\omega_{*,*}(tX) \to \Omega_{*,*}(X)$$
is an isomorphism for all quasi-projective derived schemes $X$.
\end{thm}
\begin{proof}
As a special case of the preceding Lemma \ref{LeePInj}, we know that the natural map
$$\omega_{*,*}(pt) \otimes_\Laz \Omega_*(X) \to \Omega_{*,*}(X)$$
is injective. On the other hand, as $\Omega_*(X)$ is naturally equivalent to $\Omega_*(tX)$, which in turn is naturally equivalent to $\omega_*(tX)$, we conclude that the map
$$\omega_{*,*}(pt) \otimes_\Laz \omega_*(tX) \to \Omega_{*,*}(X)$$
is injective. On the other hand, by the results of Lee-Pandharipande, this map is naturally equivalent to the map
$$\omega_{*,*}(tX) \to \Omega_{*,*}(X),$$
which we already know to be a surjection by Lemma \fref{LeePSurj}, and the claim follows.
\end{proof}

\section{Bivariant precobordism with line bundles and the weak projective bundle formula}\label{CobordismOfLines}

$\Omega^{*,*}(X \to Y)$ is a bivariant theory with the product $\bullet_{\otimes}$. We can see that by the definition of $\bullet_{\otimes}$ the subtheory $\Omega^{*,1}(X \to Y)$ of cobordism cycles of line bundles becomes a bivariant theory, thus we call it \emph{the bivariant cobordism theory with line bundles}. This theory comes with a natural inclusion $(\Omega^*, \bullet) \to (\Omega^{*,1}, \bullet_\otimes)$ which equips a cycle $[V \to X]$ with the trivial line bundle $\OO_V$. This also induces the orientation $\theta_\otimes$ on $\Omega^{*,1}$.

The purpose of this section is to study the structure of bivariant theories of line bundles, and use the gained knowledge in order to compute the cobordism group of $\Proj^n \times X \to Y$ in terms of that of $X \to Y$ (weak projective bundle formula, see Theorem \fref{WPBF}). In Section \fref{PrecobordismSect} we introduce the notion of a \emph{bivariant precobordism theory} $\B^{*}$ and the \emph{associated theory with line bundles} $\B^{*,1}$, which form a natural class of bivariant theories for which the results of this section hold. We note that they are quite general: our results hold over an arbitrary Noetherian base ring $A$ of finite Krull dimension. In Section \fref{StructOfPrecobOfLinesSect} we show that, additively (that is, disregarding the bivariant product), $\B^{*,1}$ is just a direct sum of copies of $\B^*$. In Section \fref{WPBFSect} we connect the structure of $\B^{*,1}(X \to Y)$ to the structure of $\B^{*}(\Proj^n \times X \to Y)$ using the auxiliary theories $\B^{*,1}_\mathrm{gl}$ and $\B^{*}_{\Proj^\infty}$, and use the structural results obtained in Section \fref{StructOfPrecobOfLinesSect} in order to prove the weak projective bundle formula.

Throughout the section, unless otherwise specified, $A$ will be a Noetherian ring of finite Krull dimension.

\subsection{Bivariant precobordism theories}\label{PrecobordismSect}

In this subsection, we are going to introduce the theories for which the results of this section apply. Let us denote by $\mathcal{M}^*_+$ the universal bivariant theory of Yokura (as recalled in Section \fref{FM-BT}) applied to the homotopy category of the $\infty$-category of quasi-projective derived $A$-schemes with proper morphisms as confined morphisms, quasi-smooth morphisms as specialized morphisms, and all homotopy Cartesian squares as independent squares. We recall that $\mathcal{M}^d_+(X {\xrightarrow f} Y)$ is the free Abelian group on homotopy classes of proper maps $h:V \to X$ so that the composition $f\circ h:V \to Y$ is a quasi-smooth morphism of relative virtual dimension $-d$, modulo the relation identifying disjoint union with summation. Recall also that a quasi-smooth morphism $f: X \to Y$ of relative virtual dimension $-d$ has a canonical orientation 
$$\theta(f) := [X \xrightarrow{\mathrm{Id}_X} X] \in \mathcal{M}_+^d(X \to Y),$$
and these are stable under pullback. 
\begin{defn}\label{PrecobordismDef}
Let $\B^*$ be a quotient theory of $\mathcal{M}^*_+$. Then we say that
\begin{enumerate}
\item $\B^*$ is a \emph{naive cobordism theory} if, given $W \to \Proj^1 \times X$ is a projective morphism so that the composition $W \to {\Proj^1 \times X \xrightarrow {\op{id}_{\Proj^1} \times f}} \Proj^1 \times Y$ is quasi-smooth of relative dimension $d$, then 
$$[W_0 \to X] = [W_\infty \to X] \in \B^{-d}(X \xrightarrow f Y),$$
where $W_0$ and $W_\infty$ the homotopy fibres of $W \to \Proj^1 \times X$ lying over $0 \times X$ and $\infty \times X$ respectively;

\item $\B^*$ is a \emph{precobordism theory} if it is a naive cobordism theory, and if given line bundles $\Li_1$ and $\Li_2$ on $X$, we have
\begin{align}\label{GeometricFGL}
c_1(\Li_1 \otimes \Li_2) &= c_1(\Li_1) + c_1(\Li_2) - c_1(\Li_1) \bullet c_1(\Li_2) \bullet [\Proj_1 \to X] \\
&- c_1(\Li_1) \bullet c_1(\Li_2) \bullet c_1(\Li_1 \otimes \Li_2) \bullet ([\Proj_2 \to X] - [\Proj_3 \to X]) \notag
\end{align}
in $\B^1(X) :=\B^1(X \xrightarrow {\op{id}_X} X)$, where 
\begin{align*}
\Proj_1 &:= \Proj_X(\Li_1 \oplus \OO); \\
\Proj_2 &:= \Proj_X(\Li_1 \oplus (\Li_1 \otimes \Li_2) \oplus \OO); \\
\Proj_3 &:= \Proj_{\Proj_X(\Li_1 \oplus (\Li_1 \otimes \Li_2))}(\OO(-1) \oplus \OO).
\end{align*}
\end{enumerate}
\end{defn}

The following result follows easily from the definition.
\begin{lem}\label{NilpPrecobordismChernClass}
Suppose $\B^*$ is a bivariant precobordism theory and $\Li$ is a line bundle on a quasi-projective derived $A$-scheme $X$. Then the first Chern class $c_1(\Li) \in \B^1(X)$ is nilpotent (with respect to $\bullet$).
\end{lem}
\begin{proof}
The proof splits to three parts.
\begin{enumerate}
\item Suppose $\Li$ is globally generated. By the Noetherianity hypothesis on $A$, we can find finitely many global sections $s_1,...,s_n$ of $\Li$ that generate. These sections correspond to virtual Cartier divisors $D_1,...,D_n$ whose total derived intersection is empty. As $c_1(\Li) = [D_i \to X] \in \B^1(X)$ for any $i$ by Proposition \fref{SectionAxiom} (whose proof only uses homotopy fibre relation), we see that $c_1(\Li)^n = 0$. 

\item Suppose that the dual bundle $\Li^\vee$ is globally generated. We can now use (\fref{GeometricFGL}) with $\Li_1 = \Li$, $\Li_2 = \Li^\vee$ to conclude that
$$c_1(\Li) + c_1(\Li^\vee) - c_1(\Li) \bullet c_1(\Li^\vee) \bullet [\Proj_X(\Li \oplus \OO) \to X] = 0$$
(recall that $c_1(\OO_X) = 0$). It follows that
$$c_1(\Li^\vee) = - {c_1(\Li) \over 1 - c_1(\Li) \bullet [\Proj_X(\Li \oplus \OO) \to X]}$$
which is clearly nilpotent (and well defined) since $c_1(\Li)$ is nilpotent by the first part.

\item In general, as $X$ is quasi-projective, any line bundle $\Li$ is equivalent to $\Li_1 \otimes \Li_2^\vee$ with $\Li_1, \Li_2$ globally generated. As both $c_1(\Li_1)$ and $c_1(\Li_2^\vee)$ are nilpotent, the nilpotency of $c_1(\Li)$ follows from the formula (\fref{GeometricFGL}). \qedhere 
\end{enumerate}
\end{proof}

Any theory satisfying the homotopy fibre relation and double point cobordism of \cite{LP} is a bivariant precobordism theory.
\begin{prop}[cf. \cite{LP} Section 0.3]\label{DPCProp}
Suppose $\B^*$ is a naive cobordism theory in the sense of Definition \fref{PrecobordismDef}. Then $\B^*$ is a precobordism theory if and only if, for any quasi-smooth $W \to \Proj^1 \times X$ with fibres $W_0$ over $\{ 0 \} \times X$ and $W_\infty$ over $\{ \infty \} \times X$, so that $W_\infty$ equivalent to the sum of divisors $A \to W$ and $B \to W$, the double point cobordism formula
\begin{equation}\label{DPCFormula}
[W_0 \to X] = [A \to X] + [B \to X] - [\Proj \to X]
\end{equation}
holds in $\B^1(X)$, where
$$\Proj := \Proj_Z(\OO(A) \oplus \OO)$$
with $Z$ the derived intersection of $A$ and $B$ in $W$.
\end{prop}
\begin{proof}
The proof is an easy imitation of the proof of \cite{LP} Lemma 5.2 in the derived setting, which we are going to give for the sake of completeness.   Suppose $\Li_1$ and $\Li_2$ are line bundles on quasi-projective derived $A$-scheme $X$, and let $A$, $B$ and $C$ be divisors in the linear systems of $\Li_1$, $\Li_2$ and $\Li_1 \otimes \Li_2$ respectively. Note that if $A,B,C=W_0$ are as in the statement, then (\fref{DPCFormula}) is a special case of (\fref{GeometricFGL}) as by assumption $c_1(\Li_2) \bullet c_1(\Li_1 \otimes \Li_2)$ vanishes ($B$ and $W_0$ do not meet). Therefore the only nontrivial part is to show that the double point cobordism formula implies (\fref{GeometricFGL}).

Let $X_1$ be the derived blow up of $X$ at $A \cap C$, and let $A_1, C_1 \to X_1$ be the strict transforms of $A$ and $C$ respectively. Note that by Proposition \fref{BlowUpOfIntersection} $A_1$ and $C_1$ do not meet inside $X_1$. Let $B_1 \to X_1$ be the natural map
$$\bl_{A \cap B \cap C}(B) \to \bl_{A \cap C}(X),$$
which fits by the contravariant naturality of blow-ups (see \cite{Khan} Theorem 4.1.5 (ii)) in a derived Cartesian square

\begin{equation}\label{B1Diagram}
\begin{tikzcd}
B_1 \arrow[]{r} \arrow[]{d} & X_1 \arrow[]{d} \\
B \arrow[]{r} & X
\end{tikzcd}
\end{equation}
and is therefore a virtual Cartier divisor. Note that the divisors $A_1 + B_1$ and $C_1$ are rationally equivalent inside $X_1$.

Next, blow up $X_1$ at the intersection $B_1 \cap C_1$ to obtain $X_2$ and the strict transformations $A_2, B_2, C_2 \to X_2$, which are naturally equivalent to $A_1, B_1$ and $C_1$ respectively. As $C_2$ does not meet either $A_2$ or $B_2$ inside $X$, and as $A_2 + B_2$ is rationally equivalent to $C_2$, we get two sections $s_1$ and $s_2$ of $\Li := \OO_{X_2}(A_2 + B_2) \simeq \OO_{X_2}(C_2)$ generating $\Li$. Together with the natural map $X_2 \to X$, this induces a morphism 
$$\pi: X_2 \to \Proj^1 \times X,$$
whose fibre over $\{0\} \times X$ is equivalent to $C$, and whose fibre over $\{\infty \} \times X$ is sum of the divisors $A_2 \simeq A$ and $B_2$ in $X_2$. As the intersection of $A_2$ and $B_2$ inside $X_2$ is equivalent to $A \cap B =: Z$, the double point cobordism formula (\fref{DPCFormula}) implies that 
\begin{align}\label{AlmostGeometricFGL}
c_1(\Li_1 \otimes \Li_2) &= [C \to X] \\
&= [A \to X] + [B_2 \to X] - [\Proj_Z(\Li \oplus \OO) \to X] \notag \\
&= c_1(\Li_1) + c_1(\Li_2) \bullet [X_1 \to X] \notag \\
&-c_1(\Li_1) \bullet c_1(\Li_2) \bullet [\Proj_X(\Li \oplus \OO) \to X]. \quad (\text{$B_1$ and $B_2$ are equivalent, (\fref{B1Diagram})}) \notag
\end{align}
By the blow up relation, we have
\begin{align}\label{BlowUpFormula}
[X_1 \to X] &= [X \to X] - c_1(\Li_1) \bullet c_1(\Li_1 \otimes \Li_2) \bullet [\Proj(\Li_1 \oplus (\Li_1 \otimes \Li_2) \oplus \OO)] \\
&+ c_1(\Li_1) \bullet c_1(\Li_1 \otimes \Li_2) \bullet [\Proj_{\Proj(\Li_1 \oplus (\Li_1 \otimes \Li_2))}(\OO(-1) \otimes \OO) \to X] \notag
\end{align}
(see \cite{LP} Lemma 5.1, the proof in the derived case goes through word by word). Combining (\fref{AlmostGeometricFGL}) with (\fref{BlowUpFormula}), we get (\fref{GeometricFGL}), finishing the proof.
\end{proof}

\begin{ex}
The bivariant algebraic cobordism $\Omega^*$ constructed in \cite{An} over a field $k$ of characteristic zero is a precobordism theory in the above sense. Indeed, $\Omega^*$ is a naive cobordism theory by construction, and the double point cobordism relation holds for $\Omega^* = \Omega^{*,0}$ by Lemma \fref{DerivedDoublePointRelations}.
\end{ex}

\begin{cons}[Universal precobordism theory $\underline \Omega^*$]\label{UniversalPrecobordism}
It is clear that there is a universal precobordism theory $\underline \Omega^*$ over $A$ constructed from $\mathcal{M}^*_+$ by enforcing the homotopy fibre relation and either the formulas (\fref{GeometricFGL}) or (\fref{DPCFormula}) for all line bundles satisfying restrictions of each formula (compare to the construction of Section \fref{ConstructionOfCobordismOfBundles}). It is clear that any other precobordism theory is a quotient of $\underline \Omega^*$.
\end{cons}

Next we are going to construct the associated bivariant theory with line bundles. Let $\mathcal{M}^{*,1}_+$ be as in Section \fref{ConstructionOfCobordismOfBundles}: $\mathcal{M}^{d,1}_+(X {\xrightarrow f} Y)$ is the free Abelian group on equivalence classes
$$[V {\xrightarrow h} X, \Li]$$
where $h:V \to X$ is proper, the composition $f\circ h:V \to Y$ is quasi-smooth of virtual relative dimension $-d$, and $\Li$ is a line bundle on $X$. The bivariant product $\bullet = \bullet_\otimes$ makes $\mathcal{M}^{*,1}_+$ into a bivariant theory. Note that we can identify $\mathcal{M}^{*}_+$ with the subtheory of $\mathcal{M}^{*,1}_+$ consisting of cycles where the line bundle $\Li$ is trivial. 

\begin{defn}[Bivariant precobordism with line bundles $\B^{*,1}$]\label{PrecobordismOfLines}
Let $\B^* = \mathcal{M}^{*}_+/\bI$ be a precobordism theory. We define the associated \emph{precobordism with line bundles} $\B^{*,1}$ as
$$\B^{*,1} := \mathcal{M}^{*,1}_+ / \langle \bI \rangle_{\mathcal{M}^{*,1}_+}.$$
\end{defn}

\begin{rem}[Double point cobordism with a line bundle]\label{DoublePointCobordismWithBundle}
Let us record the following trivial observation here, which will be the basis for most of our arguments. Suppose we have a morphism $X \to Y$ of derived schemes. Given $W \to \Proj^1 \times X$ projective with the composition $W \to \Proj^1 \times Y$ is quasi-smooth, and a line bundle $\Li$ on $W$. Let $W_0$ and $W_\infty$ be the fibres over $0$ and $\infty$ respectively, and suppose $W_\infty \hookrightarrow W$ is the sum of two divisors $A$ and $B$ in $W$. Then
\begin{align*}
[W_0 \to X, \Li \vert_{W_0}] &= [A \to X, \Li \vert_A] + [B \to X, \Li \vert_B] - [\Proj \to X, \Li \vert_{A \cap B}] \in \B^{*,1}(X \to Y)
\end{align*}
where $\Proj$ is defined as in Proposition \fref{DPCProp}.
\end{rem}

The following results follow easily from earlier considerations.

\begin{prop}\label{RelationsInPrecobordismOfLines}
 Given a line bundle $\Li$ on $X$, denote by $[\Li]$ the element
$$[X \xrightarrow {\op{id}_X} X, \Li] \in \mathcal{M}^{0,1}_+(X).$$
The bivariant ideal $\langle \bI \rangle_{\mathcal{M}^{*,1}_+}$ as in Definition \fref{PrecobordismOfLines} consists of linear combinations of elements of form 
$$f_*([\Li] \bullet \alpha)$$
where $\alpha \in \bI(X \to Y)$, and the map $X \to Y$ factors through a proper morphism $f: X \to X'$.
\end{prop}
\begin{proof}
The proof is the same as that of Proposition \fref{RelationsInCobordismOfBundles}, which does not use any assumptions on the base ring.
\end{proof}

The above Proposition has two immediate consequences.

\begin{defn}\label{ForgetfulPrecobordismTransformation}
By Proposition \fref{RelationsInPrecobordismOfLines}, the forgetful Grothendieck transformation $\Cal F: \mathcal{M}^{*,1}_+ \to \mathcal{M}^*_+$ descends to give the \emph{forgetful transformation} $\Cal F: \B^{*,1} \to \B^*$ for any precobordism theory $\B^*$.
\end{defn}

\begin{prop}[cf. Proposition \fref{diff}]
The differential operator $\partial_{c_1}: \mathcal{M}^{*,1}_+ \to \B^{*+1,1}$ defined by the formula
$$[V \xrightarrow{f} X, \Li] \mapsto f_*(c_1(\Li) \bullet [V \to V])$$
descents to an operator
$$\partial_{c_1}: \B^{*, 1} \to \B^{*+1,1}.$$
\end{prop}
\begin{proof}
The well definefness follows as in the proof of Proposition \fref{diff}.
\end{proof}

\begin{ex}
The associated precobordism with line bundles for the bivariant algebraic cobordism $\Omega^*$ of \cite{An} is $(\Omega^{*,1}, \bullet_\otimes)$ as constructed in Section \fref{ConstructionOfCobordismOfBundles}. This follows from 
Proposition \fref{RelationsInPrecobordismOfLines} above and from Proposition \fref{RelationsInCobordismOfBundles}, after noting that an element of form $[\Li] \bullet_\oplus \alpha$, with $\alpha \in \mathcal{M}_+^{*}$, where $\mathcal{M}_+^*$ is identified with $\mathcal{M}_+^{*,0}$, is the same as the element of form $[\Li] \bullet_\otimes \alpha$, where $\mathcal{M}_+^{*}$ is now identified with the subtheory of $\mathcal{M}_+^{*,1}$ consisting of cycles with trivial line bundles.
\end{ex}

\subsection{Structure of $\B^{*,1}$}\label{StructOfPrecobOfLinesSect}

The purpose of this section is to prove the following two theorems.

\begin{thm}\label{AlgebraicCobordismOfBGm}
Let $\B^{*}$ be a precobordism theory over a Noetherian base ring $A$ of finite Krull dimension. Then the cobordism group $\B^{*,1}(pt)$ admits an $\B^*(pt)$-linear basis $([\Proj^i \to pt, \OO(1)])_{i=0}^\infty$.
\end{thm}

The following natural map is similar to that in Lemma \ref{LeePInj}.

\begin{thm}\label{KunnethFormulaForBGm}
Let $\B^{*}$ be a precobordism theory over a Noetherian base ring $A$ of finite Krull dimension. The natural {(cross product)} map
\begin{equation*}
\B^{*,1} (pt \to pt) \otimes_{\B^*(pt)} \B^{*}(X \to Y) \to \B^{*,1}(X \to Y)
\end{equation*}
defined by
$$([V \to pt, \Li], [W \to X]) \mapsto [V \to pt, \Li] \times 
{[W \to X]}= [V \times W \to X, \mathrm{pr}_1^*\Li]$$
is an isomorphism. Above, the map $\mathrm{pr}_1$ is the natural projection $V \times W \to V$.
\end{thm}

For the rest of the subsection $\B^*$ will be a fixed precobordism theory (e.g. bivariant algebraic cobordism $\Omega^*$), and $\B^{*,1}$ will denote the associated precobordism of line bundles defined in the previous subsection.

Our strategy is to first prove the surjectivity part of Theorem \fref{KunnethFormulaForBGm} (Proposition \fref{KunnethSurjectivity}), which is done by explicitly constructing algebraic cobordisms realizing desirable relations (see Lemma \fref{MainDeformationLemma}). The rest of Theorem \fref{KunnethFormulaForBGm} and Theorem \fref{AlgebraicCobordismOfBGm} then follow with a relatively little amount of effort. Until the end of this subsection, $X$ will be a quasi-projective derived scheme, $\Li$ a line bundle on $X$, and $D \hookrightarrow X$ is a virtual Cartier divisor in the linear system of $\Li$. Note that at least one such $D$ always exists, as the derived vanishing locus of the zero section is a virtual Cartier divisor. 

We begin with the following construction.

\begin{cons}\label{FirstTowerConstruction}
Let $X$, $\Li$ and $D$ be as above. Let $W = W(X, \Li, D)$ be the blow up $\bl_{\infty \times D}(\Proj^1 \times X)$, and let $\wtil\Li$ be the line bundle $\Li(-E)$ on $W$, where $E$ denotes the exceptional divisor of the blow up. The pair $W \to \Proj^1, \wtil \Li$ satisfies the following properties.
\begin{enumerate}
\item The fibre of $W$ over $0$ is equivalent to $X$, and the restriction of $\wtil\Li$ to the fibre is equivalent to $\Li$.

\item The exceptional divisor $E$ is equivalent to $\Proj_D(\Li \oplus \OO_D)$, and $\wtil\Li \vert_E$ is equivalent to $\Li(1)$. This is true because the conormal bundle of $E$ inside the blow up, which can be identified with $\OO(-E) \vert_E$, is equivalent to $\OO(1)$ by Theorem \fref{BasicPropertiesOfBlowUp} (3).

\item The strict transform $\wtil{\infty \times X}$ of $\infty \times X$ is equivalent to $X$, and the restriction of $\wtil \Li$ to it is the trivial line bundle $\OO_X$. Indeed, the restriction of the divisor $E$ to the strict transform is $D$ by Proposition \fref{BlowUpOfIntersection} (2), and therefore the line bundle $\OO_W(E)$ restricts to $\Li$. It follows that $\wtil\Li = \Li(-E)$ restricts to $\OO_X$.
\end{enumerate}
\end{cons}

We note that $W$ from the above construction can be understood as an algebraic cobordism over $X$, witnessing the equivalence between the homotopy fibres over $0 \times X$ and $\infty \times X$. Next we are going to build a tower of projective bundles on $W$.

\begin{cons}\label{SecondTowerConstruction}
Let $X$, $\Li$ and $D$ be as in Construction \fref{FirstTowerConstruction}. Let us moreover denote by $(W_0, \wtil\Li_0)$ the pair $(W, \wtil\Li)$ constructed in Construction \fref{FirstTowerConstruction}.

We define $(W_{i+1}, \wtil \Li_{i+1})$ recursively as
$$W_{i+1} := \Proj_{W_i}(\wtil\Li_i \oplus \OO_{W_i})$$
and 
$$\wtil\Li_{i+1} := \wtil\Li_{i} (1).$$
Moreover, we shall denote by $\pi_i$ be the composition of the natural maps 
$$W_{i} \to W_{i-1} \to \cdots \to W \to \Proj^1 \times X,$$
and by $\pi'_i$ the composition 
$$W_{i} \to W_{i-1} \to \cdots \to W.$$
\end{cons}

Lemma \fref{MainDeformationLemma} and  Proposition \fref{KunnethSurjectivity} follow in a straightforward fashion from the above construction. For the convenience of the reader we will record the following lemmas, concerning the fibres of $\pi_i$ and $\pi'_i$.

\begin{lem}\label{FirstFibreLemma}
Let $X$, $\Li$ and $D$ be as above, and define pairs $(T_i = T_i(X, \Li), \Li_i = \Li_i(X, \Li))$ recursively by
$$(T_0, \Li_0) := (X, \Li)$$
and 
$$(T_{i+1}, \Li_{i+1}) := (\Proj_{T_i}(\Li_{i} \oplus \OO_{T_i}), \Li_i(1)).$$
Let $\pi_i$, $\pi'_i$ be as in Construction \fref{SecondTowerConstruction}. Then:
\begin{enumerate}
\item the fibre of $\pi_i$ over $0 \times X$ is equivalent to $T_i$, and the restriction of $\wtil\Li_i$ is equivalent to $\Li_i$;

\item the fibre of $\pi'_i$ over the exceptional divisor $E$ is equivalent to the derived fibre product $D \times^\R_X T_{i+1}$, and the restriction of the line bundle $\wtil\Li_i$ to it is equivalent to $\Li_{i+1}$.
\end{enumerate}
\end{lem}
\begin{proof}
Both claims are true by inspection by focusing on the fibres of interest in  Construction \fref{SecondTowerConstruction}.
\end{proof}

\begin{lem}\label{SecondFibreLemma}
Let us define recursively
$$(P_0, \Cal M_0) := (pt, \OO_{pt})$$
and 
$$(P_{i+1}, \Cal M_{i+1}) := (\Proj_{P_i}(\Cal M_i \oplus \OO_{P_i}), \Cal M_{i}(1)),$$
and let $X, \Li, W_i, \pi'_i$ be as in Construction \fref{SecondTowerConstruction}. Then
\begin{enumerate}
\item the fibre of $\pi'_i: W_i \to W_0 = W$ over $\wtil{\infty \times X}$ is equivalent to $P_i \times X$, and the restriction of the line bundle $\wtil\Li_i$ to it is equivalent to $\Cal M_i = \Cal M_i \boxtimes \OO_X$; 

\item the fibre of $\pi'_i: W_i \to W_0$ over the intersection of the exceptional divisor $E$ and the strict transform $\wtil{\infty \times X}$ is equivalent to $P_i \times D$, and the restriction of the line bundle $\wtil\Li_i$ to it is equivalent to $\Cal M_i = \Cal M_i \boxtimes \OO_D$.
\end{enumerate}
\end{lem}

\begin{proof}
Again, this is clear by construction by focusing on the fibre of interest. Moreover, the second claim follows trivially from the first.
\end{proof}

We can now use the two previous lemmas to prove the following.

\begin{lem}\label{MainDeformationLemma}
Let $X$, $\Li$, $W_i$ and $\wtil \Li_i$ be as in Construction \fref{SecondTowerConstruction}, $T_i$ and $\Li_i$ as in Lemma \fref{FirstFibreLemma}, and $P_i, \Cal M_i$ as in Lemma \fref{SecondFibreLemma}. Then $(\pi_i: W_i \to \Proj^1 \times X, \wtil \Li_i)$ realizes the relation
\begin{align*}
[T_i \to X, \Li_i] &= [P_i \to pt, \Cal M_i] \times 1_X \\
&+ c_1(\Li) \bullet [T_{i+1} \to X, \Li_{i+1}] \\
&- c_1(\Li) \bullet [P_i \to pt, \Cal M_i] \times [\Proj_X(\Li \oplus \OO_X) \to X, \OO_X]
\end{align*}
in $\B^{*,1}(X)$. Above $[P_i \to pt, \Cal M_i] \in \B^{*,1}(pt)$.
\end{lem}

\begin{proof}
The proof is just an easy application of the double point cobordism formula (Proposition \fref{DPCProp}) to the map $\pi_i: W_i \to \Proj^1 \times X$ together with the help of the two lemmas preceding the statement. For the convenience of the reader, we are going to give detailed explanations for where the various terms come from. We first note that the left hand side of the equation comes from the fibre of $\pi_i$ over $0 \times X$ by Lemma \fref{FirstFibreLemma} (1).

The fibre of $\pi_i$ over $\infty \times X$ is a sum of two divisors: the pre-images of the exceptional divisor $E$ and the strict transform of $\infty \times X$ under the map $\pi'_i: W_i \to W_0 = W$. It is to these that we are going to apply Lemma \fref{DerivedDoublePointRelations}. The first term on the right hand side comes from the divisor over $\wtil{\infty \times X}$ by Lemma \fref{SecondFibreLemma} (1) as
$$[P_i \to pt, \Cal M_i] \times 1_X = [P_i \times X \to X, \Cal M_i \boxtimes \OO_X] \in \B^{*,1}(X).$$
The second term on the right hand side corresponds by Lemma \fref{FirstFibreLemma} (2) to the divisor lying over $E$ as 
$$c_1(\Li) = [D \to X, \OO_D]$$
and as the bivariant product is given by the derived fibre product over $X$.

Finally, we need to identify the third term on the right hand side as the third term in the double point cobordism. Note first that the intersection of the two divisors of last paragraph in $W_i$ is equivalent to $P_i \times D$ by Lemma \fref{SecondFibreLemma} (2). Moreover, as the inclusion $P_i \times D \hookrightarrow P_i \times X$ is merely the homotopy pullback of $D \to X$ along the projection, and the corresponding normal bundle is just $\OO_{P_i} \boxtimes \Li \vert_D$. Therefore the projective bundle over the intersection (as in Lemma \fref{DerivedDoublePointRelations}) is equivalent to $P_i \times \Proj_D(\Li \oplus \OO_D)$, and the corresponding line bundle over it is just $\Cal M_i \boxtimes \OO_{\Proj_D(\Li \oplus \OO_D)}$. That this corresponds to the final term of the equation follows now by recalling the definition of the bivariant exterior product $\times$ and the first Chern class $c_1(\Li)$ as in the previous paragraph. 
\end{proof}

Using the previous lemma, it is now easy to prove the desired surjectivity.

\begin{prop}\label{KunnethSurjectivity}
The natural map 
$$\B^{*,1}(pt \to pt) \otimes_{\B^*(pt)} \B^*(X \to Y) \to \B^{*,1}(X \to Y)$$
is a surjection.
\end{prop}
\begin{proof}
As the theory is generated under bivariant operations by elements of form
$$[V \to V, \Li] \in \B^*(V),$$
it is enough to show that such an element lies in the image of the map. This follows from Lemma \fref{MainDeformationLemma}: modulo the image, we have
\begin{align*}
[V \to V, \Li] &= [T_0 \to V, \Li_0] \\
&\equiv c_1(\Li) \bullet [T_1 \to V, \Li_1] \\
&\equiv c_1(\Li)^2 \bullet [T_2 \to V, \Li_2] \\
& \ \ \vdots
\end{align*}
and as the first Chern class $c_1(\Li)$ is nilpotent (by the definition of a precobordism theory), we see that $[V \to V, \Li]$ lies in the image.
\end{proof}

Next we are going to prove Theorem \fref{AlgebraicCobordismOfBGm}. We will need the following lemmas.

\begin{lem}\label{ThirdFibreLemma}
Let $\pi: \Proj^n \to pt$ be the structure map, and let $T_i(\Proj^n) = T_i(\Proj^n, \OO(1))$ and $\Li_i(\Proj^n) = \Li_i(\Proj^n, \OO(1))$ be defined as in Lemma \fref{FirstFibreLemma}. Then
$$\pi_* \bigl( c_1(\OO(1))^n \bullet [T_n(\Proj^n) \to \Proj^n, \Li_n(\Proj^n)] \bigr) = [P_n \to pt, \Cal M_n]$$
where $P_n$ and $\Cal M_n$ are defined as in Lemma \fref{SecondFibreLemma}.
\end{lem}
\begin{proof}
Recall that $c_1(\OO(1))^n \in \B^{*,1}(X)$ is represented by the cycle associated to a morphism $pt \to \Proj^n$ so that $\OO(1) \vert_{pt}$ is trivial. Moreover, it is clear that the fibre over $pt$ of any $(T_i(\Proj^n), \Li_i(\Proj^n))$ is $(P_i, \Cal M_i)$, from which the claim follows in the case $i = n$.
\end{proof}

\begin{lem}\label{ProjectiveSpacesGenerate}
The group $\B^{*,1}(pt)$ is generated as an $\B^*(pt)$-module by the elements of form $[\Proj^i \to pt, \OO(1)]$, where $i=0,1,2,...$.
\end{lem}

\begin{proof}
Suppose $[V \to pt, \Li] \in \B^{*,1}(pt)$. By applying Lemma \fref{MainDeformationLemma} as in the proof of Proposition \fref{KunnethSurjectivity} to 
$$[V \to V, \Li] \in \B^{*,1}(V),$$
and then applying the Gysin pushforward morphism $\B^{*,1}(V) \to \B^{*,1}(pt)$, we see that it is enough to show that $[P_i \to pt, \Cal M_i]$ is expressible in the desired form. 

We will proceed by induction on $i$, the case $i=0$ being trivial. Suppose the claim is known to hold up to $i$. Consider the elements in the $\B^*(\Proj^{i+1})$-linear span of
$$[P_0 \to pt, \Cal M_0] \times 1_{\Proj^{i+1}}, [P_1 \to pt, \Cal M_1]\times 1_{\Proj^{i+1}}, ..., [P_i \to pt, \Cal M_i]\times 1_{\Proj^{i+1}}$$
in $\B^{*,1}(\Proj^{i+1})$, which we are going to call \emph{easy elements} for the rest of the proof. Note that any easy element pushes forward to an element expressible in the desired form by the inductive assumption.

Next we will apply Lemma \fref{MainDeformationLemma} to
$$[\Proj^{i+1} \to \Proj^{i+1}, \OO(1)] \in \B^{*,1}(\Proj^{i+1}).$$
It follows that, modulo easy elements, the above element is equivalent to
\begin{align*}
[\Proj^{i+1} \to \Proj^{i+1}, \OO(1)] &= [T_0(\Proj^{i+1}) \to \Proj^{i+1}, \Li_0(\Proj^{i+1})] \\
&\equiv c_1(\OO(1)) \bullet [T_1(\Proj^{i+1}) \to \Proj^{i+1}, \Li_1(\Proj^{i+1})] \\
& \ \ \vdots \\
&\equiv c_1(\OO(1))^{i+1} \bullet[T_{i+1}(\Proj^{i+1}) \to \Proj^{i+1}, \Li_{i+1}(\Proj^{i+1})]
\end{align*}
in the notation of Lemma \fref{ThirdFibreLemma}. Note that the first element pushes forward to $[\Proj^{i+1} \to pt, \OO(1)]$, and, by Lemma \fref{ThirdFibreLemma}, the last element above pushes forward to $[P_{i+1} \to pt, \Cal M_{i+1}] \in \B^{*,1} (pt)$. Hence,
$$[P_{i+1} \to pt, \Cal M_{i+1}] = [\Proj^{i+1} \to pt, \OO(1)] + \beta_i \bullet [\Proj^{i} \to pt, \OO(1)] + \cdots + \beta_0 \bullet 1 \in \B^{*,1}(pt),$$
where $\beta_i \in \B^*(pt)$, proving the claim.
\end{proof}

We are now ready to prove the main theorems.

\begin{proof}[Proof of Theorem \fref{AlgebraicCobordismOfBGm}]
The proof is just an imitation of the methods of \cite{LeeP}, but as the proof is very explicit in this special case, we are going to write it down. We are going to define a homomorphism
$$\psi_0: \B^{*,1}(pt) \to \B^{*}(pt)$$
so that
$$\psi_0(\beta_0 + \beta_1 \bullet [\Proj^1 \to pt, \OO(1)] + \cdots +\beta_n \bullet [\Proj^n \to pt, \OO(1)]) = \beta_0$$
where $\beta_i \in \B^*(pt)$ for all $i \geq 0$. Here $\B^*(pt)$ is identified with $\B^{*,0}(pt)$ when we consider the above bivariant product $\beta_i \bullet [\Proj^{i} \to pt, \OO(1)]$ for each $i$. The claim follows from the existence of such a $\psi_0$: indeed, by Lemma \fref{ProjectiveSpacesGenerate} the elements $[\Proj^i \to pt, \OO(1)]$, generate, so it is enough to show that 
$$\beta_0 + \beta_1 \bullet [\Proj^1 \to pt, \OO(1)] + \cdots +\beta_n \bullet [\Proj^n \to pt, \OO(1)] = 0 \in \B^{*,1}(pt)$$
implies $\beta_i = 0$ for all $i$. But this follows from the fact that $\psi_0$ is a homomorphism and from the fact that
$$\psi_0 \partial^i_{c_1} (\beta_0 + \beta_1 \bullet [\Proj^1 \to pt, \OO(1)] + \cdots +\beta_n \bullet [\Proj^n \to pt, \OO(1)]) = \beta_i.$$

It is easy to construct $\psi_0$ as a $\B^*(pt)$-valued series of maps $\Cal F \partial_{c_1}^i$, where $\Cal F$ is the forgetful map $\B^{*,1}(pt) \to \B^*(pt)$. Indeed, we define
$$\psi_0 := \Cal F - [\Proj^1 \to pt] \Cal F \partial_{c_1} + ([\Proj^1 \times \Proj^1 \to pt] - [\Proj^2 \to pt]) {\Cal F} \partial_{c_1}^2  + \cdots$$
which gives a well defined homomorphism $\B^{*,1}(pt) \to \B^{*}(pt)$ since  $\partial_{c_1}^i \alpha = 0$ for $i \gg 0$ for any $\alpha \in \B^{*,1}(pt)$.
\end{proof}

\begin{proof}[Proof of Theorem \fref{KunnethFormulaForBGm}]
The idea is the same as above. By Proposition \fref{KunnethSurjectivity} and Lemma \fref{ProjectiveSpacesGenerate}, we can conclude that $\B^{*,1}(X \to Y)$ is generated, as an Abelian group, by elements of the form $\beta \times [\Proj^i \to pt, \OO(1)]$, hence it is enough to show that 
$$\beta_0 + \beta_1 \times [\Proj^1 \to pt, \OO(1)] + \cdots +\beta_n \times [\Proj^n \to pt, \OO(1)] = 0 \in \B^{*,1}(X \to Y)$$
implies $\beta_i = 0$ for all $i$, when all $\beta_i$ are in $\B^*(X \to Y)$. But this follows from the fact that
$$\psi_0 \partial^i_{c_1} (\beta_0 + \beta_1 \times [\Proj^1 \to pt, \OO(1)] + \cdots +\beta_n \times [\Proj^n \to pt, \OO(1)]) = \beta_i$$
where $\psi_0: \B^{*,1}(X \to Y) \to \B^*(X \to Y)$ is defined as in the previous proof.
\end{proof}

\subsection{Precobordism of trivial projective bundles}\label{WPBFSect}
The purpose of this section is to understand the structure of the groups 
$$\B^*(\Proj^n \times X \to Y),$$
where $\B^*$ is a precobordism theory of quasi-projective derived schemes over a Noetherian ring $A$. The morphisms
$$c_1\bigl( \mathrm{pr}_1^* \OO(1) \bigr)^i \bullet \theta(\mathrm{pr}_2) \bullet: \B^*(\Proj^n \times X \to Y) \to \B^{* - n + i}(\Proj^n \times X \to Y)$$
for $i = 0..n$ give rise to a morphism
$$\mathscr Proj: \bigoplus_{i=0}^n \B^{* + n - i}(X \to Y) \to \B^{*}(\Proj^n \times X \to  Y).$$

\noindent The main theorem of the subsection is the following
\begin{thm}[Weak Projective Bundle Formula]\label{WPBF}
Let $\B^*$ be a precobordism theory over a Noetherian base ring $A$ of finite Krull dimension. Then the above map 
\begin{equation*}
\mathscr Proj: \bigoplus_{i=0}^n \B^{* + n - i}(X \to Y) \to \B^{*}(\Proj^n \times X \to Y)
\end{equation*}
is an isomorphism of $\B^*(pt)$-modules.
\end{thm}

\begin{rem} The above formula has the adjective ``weak'' and the ideal formula is ``projective bundle formula'', which is that for a vector bundle $\pi:E \to X$ of rank $n+1$ and for the projectivization $\mathbb P(\pi): \mathbb P(E) \to X$ the following isomorphism would hold:
\begin{equation*}
\mathscr Proj: \bigoplus_{i=0}^n \B^{*}(X \to Y) \xrightarrow {\cong} \B^{*}(\mathbb P(E) \xrightarrow {f \circ \mathbb P(\pi)} Y).
\end{equation*}
This formula will be proved in \cite{An4}.
\end{rem}

Before embarking on the proof, let us list some nice consequences. Using the above result, and the fact that $\B^*$ has strong orientations along smooth morphisms (cf. Proposition \ref{prop-strong-orient}), we have the following easy corollary.

\begin{cor}\label{RingWPBF}
Let $\B^*$ be a precobordism theory over a Noetherian base ring $A$ of finite Krull dimension. We have natural isomorphisms of rings
$$\B^*(\Proj^n \times X) \cong \B^*(X)[t]/(t^{n+1})$$
where $t := c_1 \bigl(\op{pr_1}^* \OO(1) \bigr) \in \B^1(\Proj^n \times X)$ is the pullback of the class $c_1(\OO(1)) \in \B^1(\Proj^n)$ of a hyperplane.
\end{cor}
\begin{proof}
Consider the derived Cartesian square 
$$\xymatrix{
\Proj^n \times X \ar[d]^{\mathrm{pr}_2} \ar[r]^{\mathrm{Id}} & \Proj^n \times X \ar[d]^{\mathrm{pr}_2}\\
X \ar[r]^{\mathrm{Id}} & X.
}
$$
The commutativity of the bivariant theory $\B^*$ implies that for all $\alpha \in \B^*(X)$
$$\theta(\mathrm{pr}_2) \bullet \alpha = \mathrm{pr}_2^*(\alpha) \bullet \theta(\mathrm{pr}_2).$$
This observation, together with Theorem \fref{WPBF}, implies that the morphism
$$\bigoplus_{i=0}^n \B^{* - i}(X) \to \B^{*}(\Proj^n \times X)$$
defined using the maps
$$t^i \bullet \mathrm{pr}_2^*: \B^{*}(X) \to \B^{* + i}(\Proj^n \times X)$$
for $i = 0..n$ is an isomorphism of $\B^*(pt)$-modules. Since $t^{n+1} = 0$, the claim follows.
\end{proof}

We can use Corollary \fref{RingWPBF} to show that the first Chern classes of line bundles are controlled by a formal group law. Indeed, consider the varieties $\Proj^n \times \Proj^m$ for various $n$ and $m$ and consider the class $c_1(\OO(1,1)) \in \B^*(\Proj^n \times \Proj^m) \cong \B^*(pt)[x,y]/(x^{n+1}, y^{m+1})$. The class is uniquely expressible as a sum
$$c_1(\OO(1,1)) = \sum_{i,j} a^{n,m}_{i,j} x^i y^j$$
where $a^{n,m}_{i,j} \in \B^*(pt)$, and moreover, by naturality of Chern classes in pullbacks, $a^{n,m}_{i,j}$ do not depend of $n$ and $m$ as long as $i \leq n$ and $j \leq m$ (denote this by $a_{ij}$). We then have the following standard result.

\begin{thm}\label{PrecobordismFGL}
Let $\B^*$ be a precobordism theory over a Noetherian base ring $A$ of finite Krull dimension. Then the series
\begin{equation}\label{FGLSeries}
F_{\B^*}(x,y) := \sum_{i,j} a_{ij} x^i y^j \in \B^*(pt)[[x,y]]
\end{equation}
is a formal group law. Moreover, given a quasi-projective derived $A$-scheme $X$, $\Li_1$ and $\Li_2$ globally generated line bundles on $X$, we have
\begin{equation}\label{PrecobordismFGLEq}
c_1(\Li_1 \otimes \Li_2) = F_{\B^*}(c_1(\Li_1), c_1(\Li_2)) \in \B^*(X).
\end{equation}
\end{thm}

\begin{proof}
The final claim follows from naturality of Chern classes whenever $\Li_1$ and $\Li_2$ are globally generated. Associativity follows from the fact that
\begin{align*}
F_{\B^*}\bigl(c_1(\Li_1), F_{\B^*}(c_1(\Li_2), c_1(\Li_3))\bigr) &= c_1(\Li_1 \otimes \Li_2 \otimes \Li_3) \\
&= F_{\B^*}\bigl(F_{\B^*}(c_1(\Li_1), c_1(\Li_2)), c_1(\Li_3)\bigr)
\end{align*}
(do the computations on $\Proj^n \times \Proj^m \times \Proj^k$ using the obvious choices of line bundles). Commutativity and the fact that $0$ is identity follow using similar argument. The fact that (\fref{PrecobordismFGLEq}) is true for $\Li_1, \Li_2$ arbitrary (and not just globally generated) follows from (\fref{GeometricFGL}) using an argument similar to the proof of Lemma \fref{NilpPrecobordismChernClass}. 
\end{proof}

As an immediate corollary, we get:

\begin{cor}\label{LazardMapsToPrecobordism}
Let $\B^*$ be a precobordism theory over a Noetherian base ring $A$ of finite Krull dimension. Then the formal group law of Theorem \fref{PrecobordismFGL} induces a homomorphism of rings
$$\Laz \to \B^*(pt)$$
where $\Laz$ is the Lazard ring.
\end{cor}

\begin{rem}
Theorem \fref{PrecobordismFGL} does hold even if $\Li_1$ and $\Li_2$ are not globally generated, but the proof is more complicated. Since we were mainly interested in deriving Corollary \fref{LazardMapsToPrecobordism}, for which the weaker version is sufficient, we leave the proof of the general version to \cite{An4}.
\end{rem}

We are going to prove Theorem \fref{WPBF} by embedding $\B^*(\Proj^n \times X \to Y)$ into $\B^{*,1}(X \to Y)$, whose structure was fully determined in Section \fref{StructOfPrecobOfLinesSect}. This is achieved by first showing that $\B^{*,1}(X \to Y)$ is isomorphic to another group $\B^*_{\Proj^\infty}(X \to Y)$, which morally corresponds to $\B^*(\Proj^\infty \times X \to Y)$, where one should think $\Proj^\infty$ as $B \mathbb{G}_m$. Proving this is very easy in the case of the universal precobordism $\underline{\Omega}^*$. However, in general $\B^*$ is defined as a quotient of $\mathcal{M}_+^*$ by a bivariant ideal (this is the case for example for $\Omega^*$ constructed in \cite{An}), making it convenient to describe how the groups $\B^*_{\Proj^\infty}(X \to Y)$ give rise to a bivariant theory. This task takes up most of Section \fref{BPinftyBglSect}. However, we claim that even though the proofs are quite long, they are also fundamentally very easy. Minor technical problems are also caused by the fact that not all line bundles are globally generated, and we introduce a third bivariant theory $\Omega^{*,1}_\mathrm{gl}$ to get around this issue. After all the necessary definitions and preliminary results in Section \fref{BPinftyBglSect}, we prove Theorem \fref{WPBF} in Section \fref{WPBFsubsubsect}.

\subsubsection{The theories $\B^{*,1}_\mathrm{gl}$ and $\B^*_{\Proj^\infty}$}\label{BPinftyBglSect}

The purpose of this section is to define and show the equivalence of two bivariant theories $\Omega^{*,1}_\mathrm{gl}$ and $\Omega^*_{\Proj^\infty}$, which play an auxiliary role in connecting the theory $\Omega^{*,1}$ to the precobordism of trivial projective bundles. The first of these theories is easier to define.

\begin{cons}[Precobordism with globally generated line bundles.]
Let $\mathcal{M}_{+, \mathrm{gl}}^{*,1}$ be the subtheory of Yokura's universal bivariant theory with line bundles $\mathcal{M}^{*,1}_+$ consisting of cycles 
$$[V \to X, \Li]$$
where $\Li$ is globally generated. As the trivial line bundle is globally generated, there is a Grothendieck transformation $\mathcal{M}_{+}^{*} \to \mathcal{M}_{+, \mathrm{gl}}^{*,1}$. If $\B^* \cong \mathcal{M}_{+}^{*} / \bI$ is a precobordism theory, we construct the corresponding \emph{precobordism with globally generated line bundles} $\B_\mathrm{gl}^{*,1}$ as the quotient theory $\mathcal{M}_{+, \mathrm{gl}}^{*,1} / \langle \bI \rangle_{\mathcal{M}_{+, \mathrm{gl}}^{*,1}}$. As global generation is stable under pullbacks and tensor products, we get a well defined bivariant theory $(\B^{*,1}_{\mathrm{gl}}, \bullet_\otimes)$. Moreover, it is clear that we have a comparison map $\B^{*,1}_{\mathrm{gl}} \to \B^{*,1}$.
\end{cons}

\begin{rem}\label{GlobGenHFib}
Let $X \to Y$ be a map of derived schemes. By the usual argument, we see that given a proper map $W \to \Proj^1 \times X$ so that the composition $W \to \Proj^1 \times Y$ is quasi-smooth, and a \emph{globally generated} line bundle $\Li$ on $W$, we get the equality
\begin{equation}\label{GlobGenHFibFormula}
[W_0 \to X, \Li \vert_{W_0}] = [W_\infty \to X, \Li \vert_{W_\infty}] \in \B^{*,1}_\mathrm{gl}(X \to Y)
\end{equation}
where $W_0$ and $W_\infty$ are the homotopy fibres over $0$ and $\infty$ respectively. Note that it is not clear that (\fref{GlobGenHFibFormula}) would hold if we only assumed $\Li \vert_{W_0}$ and $\Li \vert_{W_\infty}$ to be globally generated. This will cause a minor technical inconvenience later.
\end{rem}

Next we are going to define $\B^*_{\Proj^\infty}$. The idea is that $\B^*_{\Proj^\infty}(X \to Y)$ should be just the cobordism group $\B^*(\Proj^{\infty} \times X \to Y)$. The standard way of making sense of such a thing is by approximating $\Proj^\infty$ as the infinite union of finite dimensional projective spaces. Choose the standard linear embeddings 
\begin{equation*}
pt \stackrel{\iota_0}{\to} \Proj^1 \stackrel{\iota_1}{\to} \Proj^2 \stackrel{\iota_2}{\to} \Proj^3 \stackrel{\iota_3}{\to} \cdots,
\end{equation*}
which induce a sequence
\begin{equation}\label{ApproxSequence}
\B^*(X \xrightarrow f Y) \stackrel {(\iota_0 \times 1)_*} \longrightarrow \B^*(\Proj^1 \times X \xrightarrow{f \circ \op{pr}_2} Y) \stackrel{(\iota_1 \times 1)_*}\longrightarrow \B^*(\Proj^2 \times X \xrightarrow {f \circ \op{pr}_2} Y) \longrightarrow \cdots
\end{equation}
of bivariant cobordism groups. Here $\op{pr}_2: \Proj^k \times X \to X$ is the projection to the second factor. Noticing that for each linear embedding $\iota_k: \Proj^k \to \Proj^{k+1}$ we have the following commutative diagram
$$\xymatrix
{
\Proj^k  \times X \ar[d]_{\iota_k \times 1} \ar[rr]^{f \circ \op{pr}_2} && Y\\
\Proj^{k+1} \times X \ar[urr]_{f \circ \op{pr}_2} 
}
$$
the bivariant pushforward $(\iota_k \times 1)_*:\B^*(\Proj^k \times X  \xrightarrow {f \circ \op{pr}_2} Y)  \to \B^*(\Proj^{k+1} \times X  \xrightarrow {f \circ \op{pr}_2} Y)$ is well-defined.

\begin{defn}
We define $\B^*_{\Proj^\infty}(X \to Y)$ as the colimit of the sequence in (\fref{ApproxSequence}). Note that only the bivariant pushforward and pullback of $\B^*$ naturally induce operations on $\B^*_{\Proj^\infty}$; as of now, we don't have a bivariant product on $\B^*_{\Proj^\infty}$.
\end{defn}

The goal is to show that the groups $\B^{*,1}_\mathrm{gl}(X \to Y)$ and $\B^*_{\Proj^\infty}(X \to Y)$ are isomorphic. In the case of $\B^* = \underline{\Omega}^*$ this can be done immediately: we define morphisms $\underline{\Omega}^*_{\Proj^\infty}(X \to Y) \to \underline{\Omega}^{*,1}(X \to Y)$ and $\underline{\Omega}^{*,1}(X \to Y) \to \underline{\Omega}^*_{\Proj^\infty}(X \to Y)$ on the cycle level using formulas
$$[V \xrightarrow{f} \Proj^n \times X] \mapsto [V \xrightarrow{\mathrm{pr}_2 \circ f} X, f^* \OO(1)]$$
and
$$[V \xrightarrow{g} X, \Li] \mapsto [V \xrightarrow{\wtil g} \Proj^n \times X]$$
respectively, where $\wtil g$ is induced by $n+1$ global sections generating $\Li$. Showing that these maps are well defined is easy, since the only relations in $\underline{\Omega}^*$ come from double point cobordism formula (see Proposition \fref{DPCProp}). However, in general, the only thing we know about the relations of $\B^*$ is that they form a bivariant ideal. This forces us to employ another strategy: we have to show that $\B^*_{\Proj^\infty}$ is a bivariant theory! Even though the general proof is much longer, we believe it to be absolutely essential since $\underline{\Omega}^*$ does not have enough relations to be called \emph{bivariant algebraic cobordism}!

Let us consider the cobordism cycle 
\begin{equation*}
[V \stackrel{h}{\to} X, \Li ] \in \mathcal{M}_{+, \mathrm{gl}}^{*,1}(X \xrightarrow f Y).
\end{equation*}
Choosing global sections $s_0,...,s_n$ which generate $\Li$, we get a morphism 
$s: V \to \Proj^n$ and therefore also an element
\begin{equation*}
[V \stackrel{s \times h}{\to} \Proj^n \times X] \in \B^{*}(\Proj^n \times X \xrightarrow {f \circ \op{pr}_2} \times Y)
\end{equation*}
Thus we define 
$$\Phi_{\Proj^{\infty}}:\mathcal{M}_{+, \mathrm{gl}}^{*,1}(X \xrightarrow f Y) \to \B_{\Proj^{\infty}}(X \xrightarrow f Y)$$
by
$$ \Phi_{\Proj^{\infty}}([V \stackrel{h}{\to} X, \Li ]):=\iota_{n \infty}([V \stackrel{s \times h}{\to} \Proj^n \times X]).$$
The fact that the map $\Phi_{\Proj^\infty}$ is well defined follows from the next lemma:
 
\begin{lem} 
The image $\Phi_{\Proj^\infty}([V \to X, \Li]) \in \B_{\Proj^\infty}^*(X \to Y)$ does not depend on the choice of global sections of $\Li$.
\end{lem}
\begin{proof}
Suppose we have two sequences $s_0,...,s_n$ and $s'_0,...,s'_m$ of generating global sections of $\Li$. We need to show that they define the same element of $\B^*(\Proj^i \times X \to Y)$ for $i$ large enough. In fact, $i = n+m+1$ is good enough: the sections $x_0 s_0, ..., x_0 s_n, x_1 s'_0, ..., x_1 s'_m$ generate the line bundle $\Li(1)$ on $\Proj^1 \times V$, which induces a proper map
\begin{equation*}
\Proj^1 \times V \to \Proj^1 \times \Proj^{n+m+1} \times X,
\end{equation*}
so that the composition with target $\Proj^1 \times Y$ is quasi-smooth. 

On the fibres over $0$ and $\infty$ respectively, the maps $V \to \Proj^i$ are induced by the sequences sections $s_0,...,s_n,0,...,0$ and $0,...,0,s'_0,...,s'_m$ respectively. On the other hand, there is clearly a chain of $m+1$ cobordisms showing that the sequences $0,...,0,s'_0,...,s'_m$ and $s'_0,...,s'_m,0,...,0$ give rise to cobordant maps to $\Proj^{n+m+1} \times X$, showing at last that the elements in $\B^*_{\Proj^\infty}(X \to Y)$ determined by the two sequences of sections of $\Li$ agree.
\end{proof}

We have shown that there exists a well defined map $\mathcal{M}_{+,\mathrm{gl}}^{*,1} \to \B_{\Proj^\infty}^*$. It wouldn't be hard to show directly that this map descends the homotopy fibre relation, but unfortunately the relations we have for the theory $\B_{\mathrm{gl}}^{*,1}$ are defined in the terms of bivariant ideal generated by bivariant subset. Hence, the cleanest way of showing that the above map descends to a map $\B_{\mathrm{gl}}^{*,1} \to \B_{\Proj^\infty}^*$, is to consider its functorial behavior. We begin with an easy lemma.

\begin{lem}
The map $\mathcal{M}_{+,\mathrm{gl}}^{*,1} (X \xrightarrow f Y) \to \B_{\Proj^\infty}^* (X \xrightarrow f Y)$ commutes with bivariant pushforwards and pullbacks.
\end{lem}
\begin{proof}
Commuting with pushforwards is a tautology, as it does not change $V$ or $\Li$. We prove commuting with pullback. Suppose we have a cobordism cycle $[V \xrightarrow h X, \Li]$ and a map $g: Y' \to Y$, and form the homotopy Cartesian diagram
\begin{equation*}
\begin{tikzcd}
V' \arrow[]{r}{h'} \arrow[]{d}{g''} & X' \arrow[]{r}{f'} \arrow[]{d}{g'} & Y' \arrow[]{d}{g} \\
V \arrow[]{r}{h} & X \arrow[]{r}{f} & Y
\end{tikzcd}
\end{equation*}
Choose global sections $s_0, ..., s_n$ generating $\Li$, and choose their pullbacks $s'_0,...,s'_n$ to generate the pulled back bundle $\Li'$ on $V'$. Now the diagram
\begin{equation*}
\begin{tikzcd}
V' \arrow[]{r}{(s', h')} \arrow[]{d}{g''} & \mathbb P^n \times X' \arrow[]{r}{f' \circ \mathrm{pr}_2} \arrow[]{d}{\mathrm{Id}_{\Proj^n} \times g'} & Y' \arrow[]{d}{g} \\
V \arrow[]{r}{(s, h)} & \mathbb P^n \arrow[]{r}{f \circ \mathrm{pr}_s} \times X & Y
\end{tikzcd}
\end{equation*}
is homotopy Cartesian because 
$V' \stackrel{s'}{\to} \Proj^n = V' \to V \stackrel{s}{\to} \Proj^n$ (by the universal property). As the pullback 
\begin{equation*}
g^*([V \stackrel{h}{\to} X, \Li]) = [V' \stackrel{h'}{\to} X', \Li'] 
\end{equation*}
maps to
\begin{equation*}
[V' \stackrel{(s', h')}{\longrightarrow} \Proj^{n} \times X'] = g^*([V \stackrel{(s, h)}{\longrightarrow} \Proj^{n} \times X]),
\end{equation*}
thus we are done.
\end{proof}

Our next goal is to construct a bivariant product $\bullet$ on $\B^*_{\Proj^\infty}$ that makes the map $\mathcal{M}_{+,\mathrm{gl}}^{*,1} \to \B_{\Proj^\infty}^*$ a Grothendieck transformation. As the above map is a surjection, it is enough to show that the bivariant product $\bullet_\otimes$ descends this morphism, which is done in the next lemma.

\begin{lem}
Let us have the bivariant elements $\alpha \in \B^*(\Proj^n \times X \to Y)$ and $\beta \in \B^*(\Proj^m \times Y \to Z)$, and the homotopy Cartesian square
\begin{equation*}
\begin{tikzcd}
\mathbb P^n \times \mathbb P^m \times \arrow[]{r} \arrow[]{d} X & \mathbb P^m \times Y \arrow[]{d} \\
\mathbb P^n \times X \arrow[]{r} & Y
\end{tikzcd}
\end{equation*}




\noindent Moreover, let $i$ be the Segre embedding $\Proj^n \times \Proj^m \to \Proj^{(n+1)(m+1)-1}$. Now we may define a bivariant product 
$[\alpha] \bullet [\beta]$ on $\B_{\Proj^\infty}^*$ with $[\alpha] \in \B_{\Proj^\infty}^*(X \to Y)$ and $[\beta] \in \B_{\Proj^\infty}^*(Y \to Z)$ by the formula
\begin{align*}
[\alpha] \bullet [\beta] := [(i \times \op{id}_X)_*(p^*(\alpha) \bullet \beta)],
\end{align*}
where $i \times \op{id}_X:\Proj^n \times \Proj^m \times X \to \Proj^{(n+1)(m+1)-1} \times X$ and the product $\bullet$ on the right hand side is that of $\B^*$. With these definitions the map $(\mathcal{M}_{+, \mathrm{gl}}^{*,1}, \bullet_\otimes) \to (\B^*_{\Proj^\infty}, \bullet)$ is a Grothendieck transformation. 
\end{lem}
\begin{proof}
Suppose we have bivariant cycles $[V \to X, \Li_1]$ and $[W \to Y, \Li_2]$, and let us have sequences $s_0,...,s_n$ and $s'_0,...,s'_m$ of generating global sections of $\Li_1$ and $\Li_2$ respectively. Their bivariant product is $[V' \to X, \Li_1 \boxtimes \Li_2]$ where $V'$ is given by the usual homotopy pullback diagram

\begin{equation*}
\begin{tikzcd}
V' \arrow[]{d} \arrow[]{r} & X' \arrow[]{r} \arrow[d]{} & W \arrow[]{d} \\
V \arrow[]{r} & X \arrow[]{r} & Y \arrow[]{r} & Z. 
\end{tikzcd}
\end{equation*}



Consider now the product $\alpha \bullet \beta := [V \to \Proj^n \times X] \bullet [W \to \Proj^m \times Y]$ as defined in the statement. To compute $p^*(\alpha) \bullet \beta$, we note that the diagram
\begin{equation*}
\begin{tikzcd}
V' \arrow[]{r} \arrow[]{d} & \Proj^n \times X' \arrow[]{r} \arrow[]{d} & W \arrow[]{d} \\
\Proj^m \times V \arrow[]{r} \arrow[]{d} & \Proj^m \times \Proj^n \times X \arrow[]{r} \arrow[]{d} &  \Proj^m \times Y \arrow[]{r} \arrow[]{d} & Z \\
V \arrow[]{r} & \Proj^n \times X \arrow[]{r} & Y
\end{tikzcd}
\end{equation*}
is actually homotopy Cartesian. Hence $(i \times \op{id}_X)_*(p^*(\alpha) \bullet \beta)$ is just $[V' \to \Proj^{(n+1)(m+1)-1} \times X]$ and one uses, using the universal properties of projective spaces and properties of the Segre embedding, that the pullback of $\OO(1)$ from $\Proj^{(n+1)(m+1)-1} \times X$ to $V'$ is merely $\Li_1 \boxtimes \Li_2$. 

To finish the proof, we need only to show that $\bullet$ is well defined, i.e., that it does not depend on the choice of $n$, $m$. To show the independence from $m$, consider the commutative diagrams
\begin{equation*}
\begin{tikzcd}
\Proj^m \times \Proj^n \times X  \arrow[]{r}{i \times \mathrm{Id}_X} \arrow[]{d}{j \times \mathrm{Id}_{\Proj^n \times X}} &  \Proj^{(m+1)(n+1) - 1} \times X \arrow[]{d}{j' \times \mathrm{Id}_{X}} \\
\Proj^{m+1} \times \Proj^n \times X \arrow[]{r}{i' \times \mathrm{Id}_X} & \Proj^{(m+2)(n+1) - 1} \times X,
\end{tikzcd}
\end{equation*}
where $i,i'$ are Segre embeddings and $j,j'$ the obvious linear immersions, as well as
\begin{equation*}
\begin{tikzcd}
\Proj^{m} \times \Proj^n \times X \arrow[]{r} \arrow[]{d}{j \times \mathrm{Id}_{\Proj^n \times X}} & \Proj^{m} \times Y \arrow[]{d}{j'' \times \mathrm{Id}_Y} \\
\Proj^{m+1} \times \Proj^n \times X \arrow[]{r} \arrow[]{d} & \Proj^{m+1} \times Y \arrow[]{d}{p} \\
\Proj^n \times X \arrow[]{r} & Y.
\end{tikzcd}
\end{equation*}
Suppose $\beta \in \B^*(\Proj^m \times Y \to Z)$, and $\alpha \in \B^*(\Proj^n \times X \to Z)$. For the sake of simplicity, the map $i \times \op{id}_X$, $i' \times \op{id}_X$ and so on are written simply by $i$ and $i'$ and so on, deleting $\times \op{id}_X$. We can now use the bivariant projection formula to see that
\begin{align*}
i'_* (p^*(\alpha) \bullet j''_*(\beta)) &= i'_* j_* (j''^*p^*(\alpha) \bullet \beta) \\
&= j'_*( i_* (j''^*p^*(\alpha) \bullet \beta))
\end{align*}
This shows that it does not matter for which $m$ we let $\beta \in \B^*(\Proj^m \times Y \to Z)$: in the end the result will be the same in $\B_{\Proj^\infty}^*(X \to Z)$. The independence from $n$ follows in a fairly similar way from the axioms of bivariant theories, and the proof is left for the reader.
\end{proof}

We are finally in a situation where we can easily prove the following

\begin{lem}
The map $\mathcal{M}_{+,\mathrm{gl}}^{*,1} \to \B_{\Proj^\infty}^*$ descends to a map $\B^{*,1}_{\mathrm{gl}} \to \B_{\Proj^\infty}^*$.
\end{lem}
\begin{proof}
One can show, as in Proposition \fref{RelationsInCobordismOfBundles}, that the kernel of $\mathcal{M}_{+,\mathrm{gl}}^{*,1} \to \B^{*,1}_{\mathrm{gl}}$ consists of linear combinations of elements of form $f_*([\Li] \bullet_\otimes r)$, where $r$ is a relation for $\B^*$. 
The result now follows from the fact that the map $\mathcal{M}_{+,\mathrm{gl}}^{*,1} \to \B_{\Proj^\infty}^*$ is a map of bivariant theories.
\end{proof}

On the other hand, we have a map in the other direction:

\begin{lem}
We have a map $\B_{\Proj^\infty}^*(X \to Y) \to \B^{*,1}_{\mathrm{gl}}(X \to Y)$ defined on the level of cycles as 
\begin{equation*}
[V \to \Proj^n \times X ] \mapsto [V \to X, \Li]
\end{equation*}
where $\Li$ is the pullback of $\OO(1)$ from $\Proj^n$ to $V$ (more precisely, $\Li:= (\op{pr}_1 \circ h)^*\OO(1)$ where $h:V \to \Proj^n \times X$ and $\op{pr}_1:\Proj^n \times X  \to \Proj^n$).
\end{lem}
\begin{proof}
We need to show that the map $\B^{*} (\Proj^n \times X \to Y) \to \B^{*,1}_\mathrm{gl}(X \to Y)$, as in the statement above, is well defined. But clearly the map is just the composition of the standard inclusion
\begin{equation*}
\B^*(\Proj^n \times X \to Y) \to \B_{\mathrm{gl}}^*(\Proj^n \times X \to Y),
\end{equation*}
the multiplication map
\begin{equation*}
[\Proj^n \times X, \OO (1)] \bullet -: \B_{\mathrm{gl}}^*(\Proj^n \times X \to Y) \to \B_{\mathrm{gl}}^*(\Proj^n \times X \to Y)
\end{equation*}
and the pushforward
\begin{equation*}
(\pi \times \op{id}_X)_* :  \B_{\mathrm{gl}}^*(\Proj^n \times X \to  Y) \to  \B_{\mathrm{gl}}^*(X \to Y)
\end{equation*}
and is therefore well defined ($\pi$ is the projection $\Proj^n \to pt$).
\end{proof}

We are now ready to prove the main theorem of this subsection:

\begin{thm}
The Grothendieck transformations described above give an isomorphism between the bivariant theories $\B^{*,1}_\mathrm{gl}$ and $\B_{\Proj^\infty}^*$. 
\end{thm}
\begin{proof}
The other direction follows from the fact that the image of $[V \to X, \Li]$ in $\B_{\Proj^\infty}^*$ does not depend on the chosen generating global sections, and the other direction follows from the fact that if we choose global sections $s_0,...,s_n$ generating $\Li$ and therefore defining a map $V \to \Proj^n$, then the pullback of $\OO(1)$ along this map is $\Li$. 
\end{proof}
\subsubsection{Proof of the weak projective bundle formula}\label{WPBFsubsubsect}

Now that we understand the theories $\B^*_\mathrm{gl}$ and $\B^*_{\Proj^\infty}$, we are ready to begin the proof of the weak projective bundle formula. We will need the following variations of results of Section \fref{StructOfPrecobOfLinesSect}.

\begin{lem}[cf. Theorem \fref{AlgebraicCobordismOfBGm}]\label{GGCobOfBGm}
Let $\B^*$ be a precobordism theory. Now the group $\B^{*,1}_\mathrm{gl}(pt)$ admits an $\B^*(pt)$-linear basis $([\Proj^i \to pt, \OO(1)])_{i=0}^\infty$. 
\end{lem}

\begin{lem}[cf. Theorem \fref{KunnethFormulaForBGm}]\label{GGKunneth}
Let $\B^{*}$ be a precobordism theory. The natural (cross product) map
\begin{equation*}
\B^{*,1}_\mathrm{gl} (pt \to pt) \otimes_{\B^*(pt)} \B^{*}(X \to Y) \to \B^{*,1}_\mathrm{gl}(X \to Y)
\end{equation*}
defined by
$$([V \to pt, \Li], [W \to X]) \mapsto [V \to pt, \Li] \times [W \to X] = [V \times W \to X, \mathrm{pr}_1^*\Li]$$
is an isomorphism. Above, the map $\mathrm{pr}_1$ is the natural projection $V \times W \to V$.
\end{lem}

The results are proven essentially the same way as in Section \fref{StructOfPrecobOfLinesSect}, but we need to make a minor modification to Construction \fref{FirstTowerConstruction}. Suppose we have a derived scheme $X$, a globally generated line bundle $\Li$ on $X$,  and a virtual Cartier divisor $D$ in the linear system of $\Li$. Recall that $W$ was defined as the derived blow up of $\Proj^1 \times X$ at $\{\infty\} \times D$, and that the line bundle $\wtil \Li$ was defined as $\Li (-E)$. Recalling the Remark \fref{GlobGenHFib}, the only thing stopping us from proving a globally generated analogue of Lemma \fref{MainDeformationLemma} is that $\wtil\Li$ might fail to be globally generated (after all, global generation is stable under tensor product, pullback, and $\OO(1)$ is globally generated for $\Proj(E)$ if $E$ is a globally generated vector bundle). Let us define a line bundle $\wtil \Li '$ on $W$ to be $\Li(1 - E)$. As 
$$\wtil \Li \vert_{W_0} \simeq \wtil \Li' \vert_{W_0}$$
and
$$\wtil \Li \vert_{W_\infty} \simeq \wtil \Li' \vert_{W_\infty},$$
one can replace $\wtil \Li$ with $\wtil \Li'$ in the arguments of Section \fref{StructOfPrecobOfLinesSect} to get the same results. If we can show that $\wtil \Li'$ is globally generated, we will also obtain proofs for Lemmas \fref{GGCobOfBGm} and \fref{GGKunneth}. This is taken care of by the following lemma.

\begin{lem}
Let $X$, $\Li$, $D$, $W$ and $\wtil\Li'$ be as above. Then $\wtil \Li'$ is a globally generated line bundle on $W$.
\end{lem}
\begin{proof}
Let $s_0,...,s_n$ be a generating sequence of global sections of $\Li$ so that $s_0$ cuts out $D$. It is clear that the sections $x_0 s_0,...,x_0 s_n, x_1 s_0, ..., x_1 s_n$ generate the line bundle $\Li (1)$ on $\Proj^1 \times X$. We claim that the sections 
\begin{equation}\label{LineGen}
x_0 s_0 + x_1 s_0, x_0 s_0 + x_1 s_1, ..., x_0 s_0 + x_1 s_n, x_0 s_0
\end{equation}
generate the line bundle $\Li (1)$ outside $\{\infty\} \times D$: indeed, outside $\{\infty\} \times X$ they restrict to sections
$$t s_0 + s_0, t s_0 + s_1, ..., t s_0 + s_n, t s_0$$
on $\A^1 \times X$, which clearly generate (here $\A^1 = \Spec(k[t])$), and on $\{\infty\} \times X$ all of them restrict to $s_0$, which by assumption cuts out $D$. We claim that the strict transformations $\wtil D_0,..., \wtil D_{n+1}$ of the divisors $D_0,...,D_{n+1}$ given as vanishing loci of the sections (\fref{LineGen}) have empty total intersection in $W$. This implies that $\wtil \Li' = \Li(1 - E)$ is globally generated, as $\wtil D_i$ are in the linear system of $\wtil \Li'$ (the centre is contained in each $D_i$).

By \fref{LinearizationInBlowUp}, we can identify the inclusions $\wtil D_i\vert_E \hookrightarrow E$, $i=0..n$ (we do not need $\wtil D_{n+1}$) as the projectivization of the natural inclusion
\begin{equation*}
\mathcal{N}_{\infty \times D / D_i} \to \mathcal{N}_{\infty \times D/\Proj^1 \times X}
\end{equation*}
which arises as the dual of the right end of the exact sequence
\begin{equation*}
\begin{tikzcd}
0 \arrow[]{r} \arrow[]{d} & \Li^\vee(-1)\vert_{\infty \times D} \arrow[]{r} \arrow[]{d}{\simeq} & \mathcal{N}^\vee_{\infty \times D / \Proj^1 \times X} \arrow[]{r} \arrow[]{d}{\simeq} & \mathcal{N}^\vee_{\infty \times D / D_i} \arrow[]{r} \arrow[]{d}{\simeq} & 0 \arrow[]{d} \\
0 \arrow[]{r} & \Li^\vee \arrow[]{r} & \Li^\vee \oplus \OO^\vee \arrow[]{r} & \OO^\vee \arrow[]{r} & 0
\end{tikzcd}
\end{equation*}
and where the left hand side map $\Li^\vee \to \Li^\vee \oplus \OO^\vee$ is given by $(1, s_i^\vee)$. Hence, dually, the inclusion $\mathcal{N}_{\infty \times D / D'_i} \to \mathcal{N}_{\infty \times D/\Proj^1 \times X}$ is merely the kernel of $(1,s_i)^T: \Li \oplus \OO \to \Li$, which is the trivial line bundle generated by the section $e_1 - s_i$ (where $e_1$ is the generator of $\OO$).

This shows that the divisors $\Proj(\mathcal{N}_{\infty \times D / D_i})$ have an empty intersection inside the exceptional divisor $\Proj(\mathcal{N}_{\infty \times D / \Proj^1 \times X})$. Indeed, as $s_0$ vanishes on $\infty \times D$, $i=0$ gives the zero section of the projective bundle, and as $s_1,...,s_n$ generate $\Li$ on $\infty \times D$ the intersection must be empty. This finishes the proof of the lemma, as well as the proofs of Lemmas \fref{GGCobOfBGm} and \fref{GGKunneth}.
\end{proof}

It is now easy to prove the following result.

\begin{prop}
The comparison map $\B^{*,1}_\mathrm{gl} \to \B^{*,1}$ is an isomorphism.
\end{prop}
\begin{proof}
By Lemma \fref{GGCobOfBGm} the map $\B^{*,1}_\mathrm{gl}(pt) \to \B^{*,1}(pt)$ is an isomorphism, and Lemma \fref{GGKunneth} extends this for all bivariant groups $\B^{*,1}_\mathrm{gl}(X \to Y) \to \B^{*,1}(X \to Y)$.
\end{proof}

Our next result implies that the natural maps
\begin{equation*}
\B^*(\Proj^n \times X \to Y) \to \B_{\Proj^\infty}(X \to Y) 
\end{equation*}
are injective.

\begin{lem}
The map
\begin{equation*}
(i \times \op{id}_X)_*: \B^* (\Proj^n \times X \to Y) \to \B^* (\Proj^{n+1} \times X \to Y)
\end{equation*}
is an injection. Here $i=\iota_n:\Proj^n \to \Proj^{n+1}$ is the linear embedding.
\end{lem}
\begin{proof}
By blowing up $\Proj^{n+1}$ at the point $[0:...:0:1]$, we obtain a space $H_n$, which is a $\Proj^1$ bundle over $\Proj^n$. We have the following commutative diagram

\begin{equation*}
\xymatrix{
& H_n \ar[d]^{\pi}\\
\mathbb P^n \ar[ur]^s \ar[r]_i & \mathbb P^{n+1}
}
\end{equation*}


As the image of $i$ lands into the the area where $\pi$ maps isomorphically, one verifies on the level of cycles that $(s \times \op{id}_X)_*$ and $(\pi \times \op{id}_X) ^* (i \times \op{id}_X)_*$ are the same map 
\begin{equation*}
\B^*(\Proj^n \times X \to Y) \to \B^*(H_n \times X \to Y).
\end{equation*}
However, the inclusion $s$ has one sided inverse given by the projection $H_n \to \Proj^n$, and therefore the map $(s \times \op{id}_X)_*$ is injective. This proves also that $(i \times \op{id}_X)_*$ is injective, so we are done.
\end{proof}

\begin{proof}[Proof of the Weak Projective Bundle Formula \fref{WPBF}.]
Consider composing the map
$$\mathscr Proj: \bigoplus_{i=0}^n \B^{* + n - i}(X \to Y) \to \B^{*}(\Proj^n \times X \to  Y)$$
with the injection
$$\iota: \B^*(\Proj^n \times X \to Y) \hookrightarrow \B^*_{\Proj^\infty}(X \to Y) \cong \B^{*,1}(X \to Y)$$
which is defined on the level of cycles by
$$[V \xrightarrow{f} \Proj^{n} \times X] \mapsto [V \xrightarrow{\mathrm{pr}_2 \circ f} X, f^* \mathrm{pr}_1^*\OO(1)].$$
By unwinding the definitions, we see that $\iota \circ \Cal Proj$ is described by the formula
$$(\alpha_0, \alpha_1, ..., \alpha_n) \mapsto \sum_{i=0}^n [\Proj^{n-i}, \OO(1)] \bullet \alpha_i$$
which shows the injectivity by Theorems \fref{AlgebraicCobordismOfBGm} and \fref{KunnethFormulaForBGm}. To check the surjectivity of $\Cal Proj$, it is enough to check the image of $\iota$ is contained in the image of $\iota \circ \Cal Proj$. But if
$$\beta = \sum_{i=0}^N [\Proj^{N-i}, \OO(1)] \bullet \alpha_i \in \B^{*,1}(X \to Y)$$
is in the image of $\iota$, then $\partial_{c_1}^{n+1} \beta = 0$. On the other hand, since $\partial_{c_1}$ is $\B^*$-linear by Proposition \fref{diff}, 
$$\partial_{c_1}^{n+1} \beta = \sum_{i=n+1}^N [\Proj^{N - i - n - 1}, \OO(1)] \bullet \alpha_i \in \B^{*,1}(X \to Y)$$
so in fact
$$\beta = \sum_{i=0}^n [\Proj^{N-i}, \OO(1)] \bullet \alpha_i \in \B^{*,1}(X \to Y),$$
which is in the image of $\iota \circ \Cal Proj$.
\end{proof}

\end{document}